\documentclass[aap,MSNbibl,seceqn,citesort,dvips]{arximspdf}
\usepackage{subenv}

\doi{10.1214/10-AAP754}
\volume{22}
\issue{3}
\pubyear{2012}
\firstpage{881}
\lastpage{930}

\makeatletter

\newcommand{\dist}{\stackrel{\mathcal{D}}{\sim}}
\newcommand{\tr}{\operatorname{trace}}
\newcommand{\bbR}{\mathbb{R}}
\newcommand{\EE}{\mathbb{E}}
\newcommand{\Normal}{\mathrm{N}}
\newcommand{\h}{\mathcal{H}}
\newcommand{\eqdef}{\stackrel{\mathrm{def}}{=}}
\newcommand{\ccdot}{\cdot}
\newcommand{\weak}{\Rightarrow}
\newcommand{\longweak}{\Longrightarrow}
\newcommand{\B}{B_s}

\newtheorem{theorem}{Theorem}[section]
\newtheorem{lemma}[theorem]{Lemma}
\newproclaim{remark}[theorem]{Remark}

\newtheorem{corollary}[theorem]{Corollary}
\newtheorem{prop}[theorem]{Proposition}
\newproclaim{assumptions}[theorem]{Assumptions}

\makeatother

\begin{document}
\begin{frontmatter}

\title{Diffusion limits of the random walk Metropolis algorithm in high dimensions}
\runtitle{Diffusion limits of the random walk Metropolis algorithm}

\begin{aug}
\author[A]{\fnms{Jonathan C.} \snm{Mattingly}\corref{}\thanksref{t1}\ead[label=e1]{jonm@math.duke.edu}},
\author[B]{\fnms{Natesh S.} \snm{Pillai}\ead[label=e2]{pillai@stat.harvard.edu}}
and~%
\author[C]{\fnms{Andrew~M.}~\snm{Stuart}\thanksref{t2}\ead[label=e3]{a.m.stuart@warwick.ac.uk}}

\runauthor{J. C. Mattingly, N. S. Pillai and A. M. Stuart}

\affiliation{Duke University, Harvard University and Warwick University}

\address[A]{J. C. Mattingly\\
Department of Mathematics\\
Center for Theoretical\\
\quad and Mathematical Sciences\\
Center for Nonlinear and Complex Systems\\
\quad and Department of Statistical Sciences\\
Duke University\\
Durham, North Carolina 27708-0251\\
USA\\
\printead{e1}}
\address[B]{N. S. Pillai\\
Department of Statistics\\
Harvard University\\
Cambridge, Massachusetts 02138\\
USA\\
\printead{e2}}
\address[C]{A. M. Stuart\\
Mathematics Institute\\
Warwick University\\
CV4 7AL\\
United Kingdom\\
\printead{e3}}
\end{aug}

\thankstext{t1}{Supported by NSF Grants DMS-04-49910 and DMS-08-54879.}

\thankstext{t2}{Supported by EPSRC and ERC.}

\received{\smonth{3} \syear{2010}}
\revised{\smonth{11} \syear{2010}}

%
\begin{abstract}
Diffusion limits of MCMC methods in high dimensions provide a~useful
theoretical tool for studying computational complexity. In particular,
they lead directly to precise estimates of the number of steps required
to explore the target measure, in stationarity, as a~function of the
dimension of the state space. However, to date such results have mainly
been proved for target measures with a~product structure, severely
limiting their applicability. The purpose of this paper is to study
diffusion limits for a~class of naturally occurring high-dimensional
measures found from the approximation of measures on a~Hilbert space
which are absolutely continuous with respect to a~Gaussian reference
measure. The diffusion limit of a~random walk Metropolis algorithm to
an infinite-dimensional Hilbert space valued SDE (or SPDE) is proved,
facilitating understanding of the computational complexity of the
algorithm.
\end{abstract}

%

\end{frontmatter}

\section{Introduction}
\label{secintro}
Metropolis--Hastings methods~\cite{Metretal53,Hast70} form a~widely
used class of MCMC methods~\cite{Liu08,CaseRobe04} for sampling
from complex probability distributions. It is, therefore, of
considerable interest to develop mathematical analyses which explain
the structure inherent in these algorithms, especially structure which
is pertinent to understanding the computational complexity of the
algorithm. Quantifying computational complexity of an MCMC method is
most naturally undertaken by studying the behavior of the method on a~family of probability distributions indexed by a~parameter and
studying the cost of the algorithm as a~function of that parameter. In
this paper we will study the cost as a~function of dimension for
algorithms applied to a~family of probability distributions found from
finite-dimensional approximation of a~measure on an
infinite-dimensional space.\vadjust{\goodbreak}

Our interest is focused on Metropolis--Hastings
MCMC methods~\cite{CaseRobe04}.
We study the simplest of these, the random walk Metropolis algorithm
(RWM). Let~$\pi$ be a~target distribution on $\mathbb{R}^N$. To
sample from~$\pi$, the RWM algorithm creates a~$\pi$-reversible
Markov chain $\{x^n\}_{n=0}^{\infty} $ which moves from a~current
state $x^0$ to a~new state
$x^1$ via proposing a~candidate $y$, using a~symmetric Markov
transition kernel such as a~random walk, and accepting $y$ with
probability $\alpha(x^0,y)$, where $\alpha(x,y) = 1 \wedge
\frac{\pi(y)}{\pi(x)}$. Although\vspace*{1pt} the proposal is somewhat naive,
within the class of all Metropolis--Hastings algorithms, the RWM is
still used in many applications because of its simplicity. The only
computational cost involved in calculating the acceptance
probabilities is the relative ratio of densities
$\frac{\pi(y)}{\pi(x)}$, as compared to, say, the Langevin algorithm (MALA)
where one needs to evaluate the gradient of $\log\pi$.

A pioneering paper in the analysis of complexity for MCMC methods in
high dimensions is~\cite{Robeetal97}. This paper studied the
behavior of random walk Metropolis methods when applied to target
distributions with density
%
%
\begin{equation}
\label{eqtarget1}
\pi^{N}(x)=\prod_{i=1}^N f(x_i),\vspace*{-2pt}
\end{equation}
where $f(x)$ is a~one-dimensional probability density function.
The authors considered a~proposal of the form
\begin{eqnarray*}
y &=& x + \sqrt{\delta} \rho,\\[-2pt]
\rho&\dist& \Normal(0, \mathrm{I}_{N}),\vspace*{-2pt}
\end{eqnarray*}
and the objective was to study the complexity of the
algorithm as a~function of the dimension~$N$ of the state
space. It was shown that choosing the proposal variance
$\delta$ to scale as $\delta=2\ell^2 \lambda^2N^{-1}$
with\setcounter{footnote}{2}\footnote{If $f$ is the p.d.f. of a~Gaussian on
$\bbR$, then $\lambda$ is its standard deviation.}
$\lambda^{-2} =\int(\frac{f'}{f})^2 f \,dx$
($\ell>0$ is a~parameter which we will discuss later)
leads to an average acceptance
probability of order $1$ with respect to dimension~$N$.
Furthermore, with this choice of
scaling, individual components of the resulting Markov
chain converge to the solution of a~stochastic differential
equation (SDE).
To state this, we define a~continuous interpolant
%
%
\begin{equation} \label{eqnMCMCe}
z^N(t)= (Nt - k) x^{k+1} + (k+1 - Nt) x^{k},\qquad k \leq
Nt < k+1.\vspace*{-2pt}
\end{equation}
%
Then~\cite{Robeetal97} shows that, when the
Markov chain is started in stationarity,
$z^{N} \Rightarrow z$ as $N \to\infty$
in $C([0,T];\bbR)$ where $z$ solves the {SDE}\footnote{Our $h(\cdot
)$ and $\ell$
are different from the $h_{\mathrm{old}}$ and $\ell_{\mathrm{old}}$
used in~\cite{Robeetal97}. However, they can be recovered from
the identities $\ell^2_{\mathrm{old}} = 2\lambda^2 \ell^2$, $h_{\mathrm
{old}}(\ell
_{\mathrm{old}}) = 2\lambda^2 h(\ell)$.}
%
%
\begin{eqnarray}
\label{eqsde}
\frac{dz}{dt} &=& \lambda^2 h(\ell) [\log f(z)]'+\sqrt{2\lambda^2
h(\ell) }\,\frac{dW}{dt},\vadjust{\goodbreak}\\
\label{eqnhl}
h(\ell) &=& 2\ell^2 \Phi\biggl(-{\ell\over\sqrt{2}}\biggr).
\end{eqnarray}
Here $\Phi$ denotes the CDF of a~standard normal distribution,
``$\Rightarrow$'' denotes weak convergence and
$C([0,T],\bbR)$ denotes the Banach space of real-valued continuous
functions defined on the interval $[0,T]$ endowed
with the usual supremum norm.
Note that the
invariant measure of the SDE~(\ref{eqsde}) has the density $f$ with
respect to the Lebesgue measure. This weak convergence result leads to
the interpretation that, started in stationarity and applied to target
measures of the form~(\ref{eqtarget1}), the RWM algorithm will take
on the order of~$N$ steps to explore the invariant
measure. Furthermore, it may be shown that the value of $\ell$ which
maximizes $h(\ell)$ and, therefore, maximizes the speed of convergence
of the limiting diffusion, leads to a~universal acceptance
probability, for random walk Metropolis algorithms applied to targets
(\ref{eqtarget1}), of approximately $0.234$.

These ideas have been generalized to other proposals, such as the MALA
algorithm in~\cite{RobeRose98}. For Langevin proposals, the scaling of
$\delta$ which achieves order 1, acceptance probabilities is
$\delta\propto N^{-1/3}$
and the choice of the constant of proportionality
which maximizes the speed of the limiting SDE results
from an acceptance probability of approximately 0.574.
Note, in particular, that this method will take
on the order of $N^{1/3}$ steps to explore the
invariant distribution. This quantifies the advantage of
using information about the gradient of log~$\pi$ in
the proposal; RWM algorithms, which do not use
this information, take on the order of~$N$ steps.

The work by Roberts and co-workers was among the first to develop a~mathematical theory of Metropolis--Hastings methods in high dimension
and does so in a~fashion which leads to clear criteria which
practitioners can use to optimize algorithmic performance, for
instance, by tuning the acceptance probabilities to $0.234$ (RWM) or
$0.574$ (MALA). Yet it is open to the criticism that, from a~practitioner's
perspective, target measures of the form~(\ref{eqtarget1}) are too
limited a~class of probability distributions to be useful and, in any
case, can be tackled by sampling a~single one-dimensional target
because of the product structure. There have been papers which
generalize this work to target measures which retain the product
structure inherent in~(\ref{eqtarget1}), but are no longer i.i.d.
(see~\cite{Beda07,RobeRose01}),
%
%
\begin{equation}
\label{eqtarget2}
\pi^{N}_0(x)=\prod_{i=1}^N \lambda_i^{-1}f(\lambda_i^{-1}x_i).
\end{equation}
However, the same criticism may be applied to
this scenario as well.

Despite the apparent simplicity of target measures of the
form~(\ref{eqtarget1}) and~(\ref{eqtarget2}),
the intuition obtained from the
study of Metropolis--Hastings methods applied to these
models with product structure
is, in fact, extremely valuable.
The two key results which need to be transferred to a~more general nonproduct measure
setting are (i) the scaling of the
proposal variance with~$N$ in order to ensure
order one acceptance probabilities; (ii)
the derivation of diffusion limits
for the RWM algorithm with a~time-scale factor which
can be maximized over all acceptance probabilities.
There is some work
concerning scaling limits for MCMC methods
applied to target measures which are not of
product form; the paper~\cite{Beda09} studies
hierarchical target distributions; the paper
\cite{Breyetal04} studies target measures
which arise in nonlinear regression and have a~mean
field structure and the paper~\cite{BreyRobe00} studies
target densities which are Gibbs measures.
We add further to this literature on scaling
limits for measures with nonproduct form
by adopting the framework studied in
\cite{Besketal08,Besketal0101,BeskStua07}.
There the authors consider a~target distribution~$\pi$
which lies in an infinite dimensional,
real separable Hilbert space
which is absolutely continuous with respect
to a~Gaussian measure $\pi_0$ with
mean zero and covariance operator~$C$ (see Section~\ref{secpre} for details).
The Radon--Nikodym derivative $\frac{d\pi}{d\pi_0}$ has the form
%
%
\begin{equation}\label{eqntargmeas}
\frac{d\pi}{d\pi_0} = M_{\Psi} \exp( -\Psi(x) )
\end{equation}
for a~real valued $\pi_0$-measurable
functional~$\Psi$ on the
Hilbert space and $M_{\Psi}$ a~normalizing constant. In Section~\ref{secpsin} we will specify
and discuss the precise assumptions on
$\Psi$ which we adopt in this paper.
This infinite-dimensional framework for the
target measures,
besides being able to capture a~huge number of useful models
arising in practice~\cite{HairStuaVoss10,Stua10},
also has an inherent mathematical structure which makes
it amenable to the derivation of diffusion limits in
infinite dimensions, while retaining links to the
product structure that has been widely studied.
We highlight two aspects of this mathematical
structure.

First, the theory of Gaussian measures naturally
generalizes from $\mathbb{R}^N$ to infinite-dimensional
Hilbert spaces. Let $(\h,\langle\cdot,\cdot\rangle,
\|\cdot\|)$ denote a~real separable
Hilbert space
with full measure under $\mu_0$ ($\Psi$
will be densely defined on $\h$). The covariance
operator $C\colon\h\mapsto\h$ is a~self-adjoint, positive\vspace*{1pt}
and trace
class operator on $\h$ with a~complete orthonormal eigenbasis
$\{\lambda_j^2, \phi_j\}$,
\[
C \phi_j = \lambda_j^2 \phi_j.
\]
Henceforth, we assume that the eigenvalues are
arranged in decreasing order and $\lambda_j > 0$.
Any function $x \in\mathcal{H}$ can be represented
in the orthonormal eigenbasis of~$C$ via the
expansion
%
%
\begin{equation}\label{eqneigenexp}
x= \sum_{j=1}^{\infty} x_j\phi_j,\qquad x_j \eqdef\langle x,\phi_j
\rangle.
\end{equation}
Throughout this paper we will often identify the
function $x$ with its coordinates
$\{x_j\}_{j=1}^{\infty} \in\ell^2$
in this eigenbasis,\vspace*{2pt} moving freely between the
two representations. Note, in particular, that
$C$ is diagonal with respect to the coordinates
in this eigenbasis.\vadjust{\goodbreak} By the Karhunen--Lo\'{e}ve~\cite{DaprZaby92} expansion,
a~realization $x$ from the Gaussian measure
$\pi_0$ can be expressed by allowing the~$x_j$ to
be independent random variables distributed
as $x_j \sim N(0,\lambda_j^2)$. Thus,\vspace*{1pt} in the
coordinates $\{x_j\}$, the prior has the product
structure~(\ref{eqtarget2}). For the random walk
algorithm studied in this paper we assume that the
eigenpairs $\{\lambda_j,\phi_j\}$ are known so that sampling
from $\pi_0$ is straightforward.

The measure~$\pi$ is absolutely continuous with respect
to $\pi_0$ and hence, any almost sure property under $\pi_0$
is also true under~$\pi$. For example, it is a~consequence
of the law of large numbers that, almost surely with
respect to $\pi_0$,
%
%
\begin{equation}
\label{eqLLN}
\frac{1}{N}\sum_{j=1}^N \frac{x_j^2}{\lambda_j^2} \to1
\qquad\mbox{as } N\to\infty.
\end{equation}
This also holds almost surely with respect to~$\pi$,
implying that a~typical draw from the target measure
$\pi$ must behave like a~typical draw from $\pi_0$
in the large $j$ coordinates.\footnote{For example, if
$\mu_0$ is the Gaussian measure
associated with Brownian motion on a~finite
interval, then~(\ref{eqLLN}) is an expression
for the variance scale in the quadratic variation,
and this is preserved
under changes of measure such as the Girsanov
formula.} This offers hope that ideas
from the product case are
applicable to measures~$\pi$ given by
(\ref{eqntargmeas}) as well.
However, the presence of~$\Psi$ prevents
use of the techniques from
previous work on this problem; the
fact that individual components of the Markov
chain converge to a~scalar SDE, as proved in~\cite{Robeetal97}, is a~direct consequence of the product structure inherent in
(\ref{eqtarget1}) or~(\ref{eqtarget2}). For target measures of the
form~(\ref{eqntargmeas}),
this structure is not present and individual components of the Markov
chain cannot be expected to converge to a~scalar SDE. However, it is
natural to expect convergence of the entire Markov chain to an
infinite-dimensional continuous time stochastic process and the
purpose of this paper is to carry out such a~program.

Thus, the second fact which makes the
target measure~(\ref{eqntargmeas}) attractive
from the point of view of establishing diffusion
limits is that fact that, as proved in a~series
of recent papers~\cite{Hairetal05,HairStuaVoss07},
it is invariant for Hilbert-space valued
SDEs (or stochastic PDES--SPDEs) with the form
%
%
\begin{equation}\label{eqnspde}
\frac{dz}{dt} = - h(\ell) \bigl(z + C \nabla\Psi(z) \bigr)+
\sqrt{2 h(\ell)} \,\frac{dW}{dt}, \qquad z(0)=z^0,
\end{equation}
where $W$ is a~Brownian motion (see~\cite{DaprZaby92})
in $\mathcal{H}$ with covariance operator~$C$. Thus, the above result
from SPDE theory
gives us a~natural candidate for the infinite-dimensional limit of an
MCMC method. We will prove such a~limit for a~RWM algorithm with
proposal covariance ${2\ell^2 \over N} C $. Moreover, we will show that
the time constant $h(\ell)$ is maximized for an average acceptance
probability of $0.234$, as obtained in~\cite{Robeetal97} in the
product case.

These measures~$\pi$ given by~(\ref{eqntargmeas}) have
a number of features which will enable us to develop
the ideas of diffusion limits for MCMC methods as
originally introduced in the i.i.d. product case.
Carrying out this program is worthwhile because
measures of the form given by~(\ref{eqntargmeas})
arise naturally in a~range of applications.
In particular, they arise in the context of nonparametric
regression in Bayesian statistics where the parameter
space is an infinite-dimensional function space.
The measure $\pi_0$ is the
prior and~$\Psi$ the log likelihood function.
Such Bayesian inverse problems are overviewed in
\cite{Stua10}. Another class of
problems leading to measures of the form~(\ref{eqntargmeas})
are conditioned diffusions (see~\cite{HairStuaVoss10}).


To sample from~$\pi$ numerically we need a~finite-dimensional
target measure. To this end, let $\Psi^N(\cdot)=\Psi(P^N\cdot)$
where $P^N$ denotes projection\footnote{Actually~$\Psi$ is only
densely defined on $\h$ but
the projection $P^N$ can also be defined on this dense subset.} (in
$\mathcal{H}$) onto
the first~$N$ eigenfunctions of~$C$.
Then consider the target measure $\pi^N$ with the form
%
%
\begin{equation}
\label{eqtarget3}
\frac{d\pi^{N}}{d\pi_0}(x) \propto\exp(-\Psi^N(x)) .
\end{equation}
This measure can be factored as the product of
two independent measures: it coincides
with $\pi_0$
on $\mathcal{H}\setminus P^N\mathcal{H}$ and has a~density
with respect to Lebesgue measure on $P^N\mathcal{H}$,
in the coordinates $\{x_j\}_{j=1}^N$.
In computational practice we implement a~random
walk method on ${\mathbb R}^N$ in the coordinate
system $\{x_j\}_{j=1}^N$, enabling
us to sample from $\pi^N$ in $P^N\mathcal{H}$.
However, in order to facilitate a~clean analysis,
it is beneficial to write this finite-dimensional random
walk method in $\mathcal{H}$, noting that the
coordinates $\{x_j\}_{j=N+1}^\infty$ in the\vspace*{2pt} representation
of functions sampled from $\pi^N$ do not then
change. We consider proposal distributions for the RWM
which exploit the covariance structure of $\pi_0$
and can be expressed in $\mathcal{H}$ as
%
%
\begin{equation}
\label{eqpropz}
y = x + {\sqrt{\frac{2\ell^2}{N}}} C^{1/2} \xi
\qquad\mbox{where } \xi= \sum_{j=1}^N \xi_j \phi_j
\quad\mbox{with } \xi_j \dist\Normal(0,1) \mbox{ i.i.d.}\hspace*{-40pt}
\end{equation}
Note that our proposal variance scales as $N^{-\gamma}$
with $\gamma=1$. The choice of~$\gamma$ in the proposal
variance affects the scale of the
proposal moves and identifying the optimal choice for~$\gamma$
is a~delicate exercise. The larger~$\gamma$ is, the
more ``localized'' the proposed move is
and, therefore, for the algorithm to explore the state
space rapidly,~$\gamma$ needs to be as small as possible.
However, if we take~$\gamma$
arbitrarily small, then the acceptance probability decreases to
zero very rapidly as a~function of~$N$. In fact, it was shown in
\cite{Besketal0101,BeskStua07,Besketal08} that, for a~variety of
Metropolis--Hastings proposals, there is $\gamma_c>0$ such that choice
of $\gamma<\gamma_c$
leads to\vadjust{\goodbreak} average acceptance probabilities which are
smaller than any inverse power of~$N$. Thus, in higher dimensions,
smaller values of~$\gamma$ lead to very poor mixing because of the
negligible acceptance probability. However, it turns out that at the
critical value $\gamma_c$, the acceptance probability is
$\mathcal{O}(1)$ as a~function of~$N$. In \cite
{Besketal0101,BeskStua07}, the value of $\gamma_c$ was
identified to be $1$ and $1/3$ for the RWM and MALA,
respectively. Finally, when using the scalings leading to $\mathcal{
O}(1)$ acceptance probabilities, it was also shown that the mean square
distance moved is maximized by choosing the acceptance probabilities
to be $0.234$ or $0.574$ as in the i.i.d. product case
(\ref{eqtarget1}). Guided by this intuition, we have
chosen $\gamma=\gamma_c= 1$ for our RWM proposal variance which,
as we will prove below, leads to
$\mathcal{O}(1)$ acceptance probabilities.

Summarizing the discussion so far,
our goal is to obtain an invariance
principle for the RWM Markov chain with
proposal~(\ref{eqpropz}) when applied to target
measures of the form~(\ref{eqntargmeas}). The
diffusion limit will be obtained in stationarity
and will be given by the SPDE~(\ref{eqnspde}).
We show that the continuous
time interpolant $z^N$ of the Markov chain $\{x^k\}$
defined by~(\ref{eqnMCMCe}) converges to
$z$ solving~(\ref{eqnspde}). This will show that, in
stationarity and properly scaled to achieve $\mathcal{O}(1)$ acceptance
probabilities, the random walk Metropolis algorithm takes $\mathcal{
O}(N)$ steps to explore the target distribution. From a~practical
point of view, the take home message of this work is that standard RWM
algorithms applied to approximations of target measures with the form
(\ref{eqntargmeas}) can be tuned to behave optimally by adjusting
the acceptance probability to be approximately $0.234$ in the case
where the proposal covariance is proportional to the covariance~$C$ in
the reference measure.
This will
lead to $\mathcal{O}(N)$ steps to explore the target measure in
stationarity. This extends the work in~\cite{Robeetal97}
and shows that the ideas developed there apply to nontrivial
high-dimensional targets arising in applications.
Although we only analyze the RWM proposal~(\ref{eqpropz}),
we believe
that our techniques can be applied to a~larger class of
Metropolis--Hastings methods, including the MALA algorithm,
and/or RWM methods with isotropic proposal variance.
In this latter case we expect to get
a different (nonpreconditioned)~$\pi$-invariant
SPDE as the limit when the dimension
goes to infinity (see
\cite{Hairetal05,HairStuaVoss07}
for analysis of these SPDEs) and a~different
(more severe) restriction on the scaling
of the proposal variance with~$N$; however,
we conjecture that the optimal acceptance
probability would not be changed.
The proposal that we study in this paper relies on
knowledge of the eigenstructure of the covariance
operator of the prior or reference measures $\pi_0$.
In some applications, this may be a~reasonable assumption,
for example, for conditioned diffusions
or for PDE inverse problems in simple geometries.
For others it may not, and then the isotropic
proposal covariance is more natural.

We analyze the RWM algorithm started at stationarity, and thus do not
attempt to answer the question of ``burn-in time'': the number of steps
required to reach stationarity and how the proposal scaling affects
the rate of convergence. These are important questions which we hope
to answer in a~future paper.\vadjust{\goodbreak} Furthermore, practitioners wishing to
sample from probability measures on function space with the form
(\ref{eqntargmeas}) should be aware that for some examples, new
generalizations of random walk Metropolis algorithms, defined on
function space, can be more efficient than the standard random walk
methods analyzed in this paper~\cite{Besketal08,BeskStua07};
whether or not they are more efficient depends on a~trade-off between
number of
steps to explore the measure (which is lower for the new generalized
methods) and cost per step (which can be higher,
but may not be).

There exist several methods in the literature to prove invariance
principles. For instance, because of the reversibility of the RWM
Markov chain, utilizing the abstract but powerful theory of Dirichlet
forms~\cite{MaRock92} is appealing. Another alternative is
to show the convergence of generators of the associated Markov
processes~\cite{EthiKurt86} as used in~\cite{Robeetal97}.
However, we chose a~more ``hands on''
approach using simple probabilistic tools, thus gaining more intuition
about the RWM algorithm in higher dimensions. We show that with the
correct choice of scaling, the one step transition for the RWM Markov
chain behaves nearly like an Euler scheme applied to
(\ref{eqnspde}). Since the noise enters~(\ref{eqnspde}) additively,
the induced It\^{o} map which takes
Wiener trajectories into solutions is
continuous in the supremum-in-time topology. This fact, which would
not be true if~(\ref{eqnspde}) had multiplicative noise, allows to
employ an argument simpler than the more general techniques often used
(see~\cite{EthiKurt86}). We first show that the martingale
increments converge weakly to a~Hilbert space-valued
Wiener process using a~martingale central limit theorem~\cite{Berg86}. Since weak
convergence is preserved under a~continuous map, the fact that the
It\^{o} map is continuous implies the RWM Markov chain converges to the
SPDE~(\ref{eqnspde}). Finally, we emphasize that
diffusion limits for the RWM proposal are necessarily
of weak convergence type. However,
strong convergence results are available for
the MALA algorithm, in fixed finite dimension
(see~\cite{BouVan09}).\looseness=-1

\subsection{Organization of the paper}
We start by setting up the notation that is used for
the remainder of the paper in Section~\ref{secrwm}. We then investigate
the mathematical structure of the RWM algorithm when applied to target
measures of the form~(\ref{eqtarget3}). Before presenting details, a~heuristic but detailed outline of the proof strategy is given for
communicating the main ideas. In Section~\ref{secmainthm} we state our
assumptions and give the proof of the main theorem at a~high
level, postponing proofs of some technical estimates.
In Section~\ref{secinvp} we prove the invariance
principle for the noise process. Section~\ref{secstec}
contains the proof of the drift and diffusion estimates.
All universal constants, unless otherwise stated, are
denoted by the letter $M$ whose precise value might
vary from one line to the next.

\section{Diffusion limits of the RWM algorithm}
\label{secrwm}

In this section we state the main theorem, set
it in context and explain the proof technique.
We first introduce an approximation of the
measure~$\pi$, namely $\pi^N$,
which is finite dimensional.
We then\vadjust{\goodbreak} state the main theorem concerning a~diffusion
limit of the algorithm and sketch the ideas of the proof
so that technical details in later sections can be
readily digested.

\subsection{Preliminaries} \label{secpre}

Recall that $\h$ is a~separable Hilbert space of real-valued functions
with inner-product and norm $\langle\cdot, \cdot\rangle$ and \mbox{$\|
\cdot\|$}. Let~$C$ be a~positive, trace class operator on $\h$. Let
$\{\phi_j,\lambda^2_j\}$ be the eigenfunctions and eigenvalues of~$C$,
respectively, so that
\[
C\phi_j = \lambda^2_j \phi_j,\qquad j \in\mathbb{N}.
\]
We assume a~normalization under which $\{\phi_j\}$
forms a~complete orthonormal basis in $\h$. We also assume
that the eigenvalues are arranged in decreasing order.
For every $x \in\h$ we have the representation
(\ref{eqneigenexp}).
Using this expansion, we define the Sobolev spaces $\h^r, r \in\bbR
$, with the inner-products and norms defined by
%
%
\begin{equation}\label{eqnSob}
\langle x,y \rangle_r\eqdef\sum_{j=1}^\infty j^{2r}x_jy_j,\qquad
\|x\|^2_r \eqdef\sum_{j=1}^\infty j^{2r} x_j^{2}.
\end{equation}
Notice that $\h^0 = \h$.
Furthermore,
$\h^r \subset\h\subset\h^{-r}$ for any $r >0$. For $r \in\mathbb{R}$,
let $B_r \dvtx\h\mapsto\h$ denote the operator which is
diagonal in the basis $\{\phi_j\}$ with diagonal entries
$j^{2r}$, that is,
\[
B_r \phi_j = j^{2r} \phi_j
\]
so that $B^{1/2}_r \phi_j = j^r \phi_j$. The operator $B_r$
lets us alternate between the Hilbert space $\h$ and the
Sobolev spaces $\h^r$ via the identities
%
%
\begin{equation}\label{eqnSob2}
\langle x,y \rangle_r=\langle B^{1/2}_r x,B^{1/2}_r y
\rangle,\qquad
\|x\|^2_r =\|B^{1/2}_r x\|^2.
\end{equation}
Let $\otimes$
denote the outer product operator in $\h$
defined by
%
%
\begin{equation}\label{eqnoutprod}
(x \otimes y) z \eqdef\langle y, z \rangle x\qquad\forall x,y,z \in\h.
\end{equation}
For an operator $L\dvtx\h^r \mapsto\h^l$, we denote the operator norm on
$\h$ by $\|\cdot\|_{\mathcal{L}(\mathcal{H}^{r},\mathcal{H}^{l})}$
defined by
\[
\|L\|_{\mathcal{L}(\mathcal{H}^{r},\mathcal{H}^{l})} \eqdef\sup_{\|
x\|_r=1} \|L x\|_l.
\]
For self-adjoint $L$ and $r=l=0$ this is, of course, the spectral
radius of $L$.
For a~positive, self-adjoint operator $D \dvtx\h\mapsto\h$, define its
trace as
\[
\tr(D) \eqdef\sum_{j=1}^\infty\langle\phi_j, D \phi_j \rangle.
\]
Since $\tr(D)$ does not depend on the orthonormal basis, 
an operator $D$ is said to be trace class if $\tr(D) < \infty$ for
some, and hence any, orthonormal basis $\{\phi_j\}$.

Let $\pi_0$ denote a~mean zero Gaussian measure on $\h$ with
covariance operator~$C$, that is, $\pi_0 \eqdef
\Normal(0,C)$. If $x \dist\pi_0$, then the~$x_j$
in~(\ref{eqneigenexp}) are independent $\Normal(0,\lambda_j^2)$
Gaussians and we may write (Karhunen--Lo\'{e}ve)
%
%
\begin{equation}\label{eqnKLexp}
x = \sum_{j=1}^\infty\lambda_j \rho_j \phi_j\qquad
\mbox{with } \rho_j \dist\Normal(0,1) \mbox{ i.i.d.} 
\end{equation}
Since $\|B^{-1/2} \phi_k\|_r = \|\phi_k\|=1$, we deduce that $\{
B^{-1/2}_r \phi_k \}$ form an orthonormal basis for $\h^r$ and, therefore,
we may write~(\ref{eqnKLexp}) as
%
%
\begin{equation}\label{eqnKLexprw}
x = \sum_{j=1}^\infty\lambda_j j^r \rho_j B^{-1/2}_r\phi_j
\qquad\mbox{with } \rho_j \dist\Normal(0,1) \mbox{ i.i.d.} 
\end{equation}
If $\Omega$ denotes the probability space for sequences
$\{\rho_j\}_{j \ge1}$, then
the sum converges in $L^2(\Omega;\h^r)$
as long as $\sum_{j=1}^\infty\lambda^2_j j^{2r} < \infty$.
Thus,\vspace*{1pt} under this condition,
the distribution induced by $\pi_0$ may be viewed
as that of a~centered Gaussian
measure on $\h^r$ with
covariance operator $C_r$ given by
%
%
\begin{equation} \label{eqncovop}
C_r = B^{1/2}_r C B^{1/2}_r.
\end{equation}
The assumption on summability is the usual
trace-class condition for Gaussian measures
on a~Hilbert space: $\tr(C_r) < \infty$. In what follows, we
freely alternate between the Gaussian measures $\Normal(0,C)$ on $\h$
and $\Normal(0,C_r)$ on $\h^r$, for values of $r$ for which
the trace-class property of $C_r$ holds.

Our goal is to sample from a~measure~$\pi$ on $\h$ given by
(\ref{eqntargmeas}),
\[
\frac{d\pi}{d\pi_0} = M_{\Psi} \exp( -\Psi(x) )
\]
with $\pi_0$ as constructed above.
Frequently in applications, the functional~$\Psi$ may not be defined
on all of $\h$,
but only on a~subset $\h^r \subset\h$ for some exponent $ r > 0$.
For instance, if $\h= L_2([0,1])$, the functional~$\Psi$ might only act
on continuous functions, in which case it is natural to define
$\Psi$ on some Sobolev space $\h^r[0,1]$ for $r > \frac12$.
Even though the Gaussian measure $\pi_0$ is defined on $\h$,
depending on the decay of the eigenvalues of~$C$,
there exists an entire range of values $r$ such that $\tr(C_r) <
\infty$ so that
the measure~$\pi_0$ has full support on $\h^r$, that is,
$\pi_0(\h^r) = 1$. From now onward we fix a~distinguished
exponent $s \geq0$ and assume that $\Psi\dvtx\h^s \mapsto\bbR$
and that the prior is chosen so that $\tr(C_s) < \infty$.
Then $\pi_0 \sim\Normal(0,C)$ on $\h$ and $\pi(\h^s)=1$; in
addition, we may view $\pi_0$ as a~Gaussian measure
$\Normal(0,C_s)$ on $\h^s$.
The precise connection between the exponent $s$ and the eigenvalues of~$C$
is given in Section~\ref{secpsin}.

In order to sample from~$\pi$ we first approximate it by a~finite-dimensional
measure. Recall that
%
%
\begin{equation} \label{eqnphat}
\widehat{\phi}_k \eqdef B^{-1/2}_s \phi_k
\end{equation}
form an orthonormal basis for $\h^s$.
For $N \in\mathbb{N}$, let $P^N \dvtx\h^s \mapsto X^N \subset\h^s
$ be the projection operator in $\h^s$ onto $X^N \eqdef{\operatorname{span}}\{
\widehat{\phi}_1,\widehat{\phi}_2,\ldots,\widehat{\phi}_N\}$,
that is,
\[
P^N x \eqdef\sum_{j=1}^N x_j \widehat{\phi}_j
\qquad\mbox{where } x_j =
\langle x, \widehat{\phi}_j \rangle_s, x \in\h^s .
\]
This shows that $X^N$ is isomorphic to $\mathbb{R}^N$.
Next, we approximate~$\Psi$
by $\Psi^N \dvtx X^N \mapsto
\mathbb{R}$ and attempt to sample from the following approximation
to~$\pi$, namely,
\[
\frac{d\pi^N}{d\pi_0}(x) \eqdef M_{\Psi^N}\exp( -\Psi^N(x))
\qquad\mbox{where }
\Psi^N(x) \eqdef\Psi(P^Nx).
\]
Note that $\nabla\Psi^N(x)=P^N\nabla\Psi(P^N x)$
and $\partial^2\Psi^N(x)=P^N \partial^2\Psi(P^N x) P^N$.
The constant $M_{\Psi^N}$ is chosen so that
$\pi^N(\h^s) = 1$. It may be shown that, for large~$N$, the measure
$\pi^N$ is close to the
measure~$\pi$ in the Hellinger metric (see~\cite{cottetal09}). Set
%
%
\begin{equation} \label{eqnCopN}
C^N\eqdef P^NCP^N,\qquad C^N_r \eqdef B^{1/2}_r C^N B^{1/2}_r .
\end{equation}
Notice that on $X^N$, $\pi^N$ has Lebesgue density\footnote{For ease
of notation we do not distinguish between a~measure
and its density.}
%
%
\begin{eqnarray} \label{eqnpitrunc}\qquad
\pi^N(x) &=& M_{\Psi^N} \exp\bigl( -\Psi^N( x ) -
\tfrac{1}{2}\langle P^Nx,C^{-1}(P^N x) \rangle\bigr),\qquad x
\in X^N \nonumber\\[-8pt]\\[-8pt]
&=&
M_{\Psi^N} \exp\bigl( -\Psi^N( x ) -
\tfrac{1}{2}\langle x,(C^N)^{-1} x \rangle\bigr)\nonumber
\end{eqnarray}
since $C^{N}$ is invertible on $X^N$ because the eigenvalues are
assumed to be strictly positive.
On $\h^s \setminus X^N$ we have that $\pi^N=\pi_0$.
Later we will impose natural assumptions on~$\Psi$
(and hence, on $\Psi^N$)
which are motivated by applications.

\subsection{The algorithm}
\label{ssecalg}

Our goal is now to sample from~(\ref{eqnpitrunc}) with $x \in
X^N$.
As explained in the \hyperref[secintro]{Introduction}, we use a~RWM proposal
with covariance operator $2\frac{\ell^2}{N} C$ on $\h$ given by
(\ref{eqpropz}).
The noise $\xi$ is finite dimensional and is independent of $x$.
Hence, even though the Markov chain evolves in
$\mathcal{H}^s$,~$x$ and $y$ in~(\ref{eqpropz}) differ only in the
first~$N$ coordinates when written in the eigenbasis
of~$C$; as a~consequence, the Markov chain does
not move at all in ${\mathcal H}^s\setminus P^N{\mathcal H}^s$
and can be implemented in $\bbR^N$. However
the analysis is cleaner when written in $\mathcal{H}^s$.
The acceptance probability also only depends
on the first~$N$ coordinates of $x$ and $y$ and has the
form
%
%
\begin{equation}\label{eqnaccprob}
\alpha(x,\xi) = 1 \wedge\exp(Q(x,\xi)),
\end{equation}
where
%
%
\begin{eqnarray}\label{eqnaccrem}
Q(x,\xi) &\eqdef&\tfrac{1}{2}\|C^{-1/2}(P^Nx)\|^2 -
\tfrac{1}{2}\|C^{-1/2}(P^Ny)\|^2\nonumber\\[-8pt]\\[-8pt]
&&{} + \Psi^N(x) -
\Psi^N(y) .\nonumber
\end{eqnarray}
The Markov chain for $\{x^k\}, k \geq0 $ is
then given by
%
%
\begin{equation}\label{eqnRWMfancy}\hspace*{32pt}
x^{k+1} = \gamma^{k+1} y^{k+1} + (1- \gamma^{k+1}) x^{k}
\quad\mbox{and}\quad
y^{k+1} = x^k + {\sqrt{\frac{2\ell^2}{N}}} C^{1/2} \xi^{k+1}
\end{equation}
with
\begin{eqnarray}
\gamma^{k+1} &\eqdef&\gamma(x^k,\xi^{k+1}) \dist
\operatorname{Bernoulli}(\alpha(x^k,\xi^{k+1}))\quad
\mbox{and} \quad
\xi^{k+1} =
\sum_{i=1}^N \xi_i^{k+1} \phi_i\nonumber\\
&&\eqntext{\mbox{where }\xi_i^{k+1} \dist
\Normal(0,1) \mbox{ i.i.d.}}
\end{eqnarray}
%
with some initial condition $x^0$.
The random variables $\xi^{k}$ and $x^0$
are independent of one another. Furthermore, conditional
on $\alpha(x^{k-1},\xi^{k})$,
the Bernoulli random variables $\gamma^k$
are chosen independently of all other sources of randomness.
This can be seen in the usual way by introducing
an i.i.d. sequence of uniform random variables $\operatorname{Unif}[0,1]$
and using these for each $k$ to construct the Bernoulli
random variable.

In summary, the Markov chain that we have described in
$\mathcal{H}^s$ is, when projected into coordinates
$\{x_j\}_{j=1}^N$, equivalent to a~standard random
walk Metropolis method for the Lebesgue density
(\ref{eqnpitrunc}) with proposal variance given by
$C^N$ on $\h$.
Recall that the target measure~$\pi$ in~(\ref{eqntargmeas}) is the
invariant measure of the SPDE~(\ref{eqnspde}). Our goal is to obtain
an invariance principle for the continuous interpolant
(\ref{eqnMCMCe}) of the Markov chain
$\{x^k\}$ started in stationarity:
to show weak convergence
of $z^{N}(t)$
to the solution $z(t)$ of the SPDE~(\ref{eqnspde}),
as the dimension $N \rightarrow\infty$.

In the rest of the section, we will give a~heuristic outline of our
main argument. The emphasis will be on the proof strategy and main
ideas. So we will not yet prove the error bounds and
use the symbol ``$\approx$'' to indicate so. Once the main skeleton
is outlined, we retrace our arguments and make them rigorous
in Sections~\ref{secmainthm},~\ref{secinvp}
and~\ref{secstec}.

\subsection{Main theorem and implications}

As mentioned earlier for fixed~$N$, the Markov chain evolves in $X^N
\subset\h^s$ and we prove the invariance principle for the Markov
chain in the Hilbert space $\h^s$ as~$N$ goes to
infinity.
Define the constant $\beta$,
%
%
\begin{equation} \label{eqnbeta}
\beta\eqdef2\Phi\bigl(-{\ell}/{\sqrt{2}}\bigr),
\end{equation}
where $\Phi$ denotes the CDF of the standard normal distribution. Note
that with this definition of $\beta$, the time scale $h(\ell)$
appearing in~(\ref{eqnspde}), and defined in~(\ref{eqnhl}), is
given by
$h(\ell) = \ell^2 \beta$.
The following is the main result of this article (it
is stated precisely, with conditions, as
Theorem~\ref{thmmain}):\vadjust{\goodbreak}
\begin{maintheorem*}
Let the initial condition $x^0$ of the
RWM algorithm be such that $x^0 \dist\pi^N$ and let
$z^N(t)$ be a~piecewise linear, continuous interpolant of the RWM algorithm (\ref
{eqnRWMfancy})
as defined in~(\ref{eqnMCMCe}). Then
$z^N(t)$ converges weakly in $C([0,T],\h^s)$ to the diffusion process
$z(t)$ given by~(\ref{eqnspde}) with \mbox{$z(0) \dist\pi$}.
\end{maintheorem*}

We will now explain the following two important implications of this result:
\begin{itemize}
\item it demonstrates that, in stationarity, the work required
to explore the invariant measure scales as $O(N)$;
\item it demonstrates that the speed at which the invariant measure is
explored, again in stationarity, is maximized by tuning the average
acceptance probability to $0.234$.
\end{itemize}

The first implication follows from~(\ref{eqnMCMCe}) since this shows
that $O(N)$ steps of the Markov chain~(\ref{eqnRWMfancy}) are
required for $z^N(t)$ to approximate $z(t)$ on a~time interval $[0,T]$ long
enough for $z(t)$ to have explored its invariant measure.
The second implication follows from~(\ref{eqnspde}) for $z(t)$
itself. 
The maximum of the time-scale $h(\ell)$ over the\vspace*{1pt} parameter
$\ell$
(see~\cite{Robeetal97}) occurs at a~universal
acceptance probability of $\widehat{\beta} =0.234$, to three decimal places.
Thus, remarkably, the optimal acceptance probability identified in
\cite{Robeetal97}
for product measures, is also optimal for the nonproduct measures studied
in this paper.

\subsection{Proof strategy}
\label{secconvSPDEOutline}

Let $\mathcal{F}_k$ denote the sigma algebra generated by $\{x^n,
\xi^n,\break \gamma^{k}$, $n \leq k\}$.
We denote the
conditional expectations $\mathbb{E}(\ccdot| \mathcal{F}_k)$ by
$\mathbb{E}_k(\ccdot)$. We first compute the one-step expected drift
of the Markov chain $\{x^k\}$. 
For notational convenience
let $x^0 = x$ and $\xi^1 = \xi$. We set $\xi^0 = 0$ and
$\gamma^0=0$.
Then, under the assumptions on $\Psi, \Psi^N$ given in
Section~\ref{secpsin}, we prove the following proposition
estimating the mean one-step drift and diffusion. The proof is
given in Sections~\ref{secdrifes} and~\ref{secdiffes}.
\begin{prop} \label{thmdrift-diffus} Let Assumptions
\ref{ass1} and~\ref{ass2} (below) hold.
Let $\{x^k\}$ be the RWM Markov chain with $x^0 = x \dist\pi^N$. Then
%
%
\begin{eqnarray}\label{thmDD-drift}\quad
N \mathbb{E}_0(x^{1} - x) &=&
-\ell^2\beta\bigl(P^N x + C^N \nabla\Psi^N( x)\bigr) + r^N,\\
\label{thmDD-diffus}
N \mathbb{E}_0[(x^{1} - x)\otimes(x^{1} - x)] &=&
2\ell^2 \beta C^N+ E^{N},
\end{eqnarray}
where the error terms $r^N$ and $E^N$ satisfy
$\mathbb{E}^{\pi^N} \|r^{N}\|_s^2 \rightarrow0$,
$\mathbb{E}^{\pi^N} \sum_{i=1}^N |\langle\phi_i,\break E^N \phi_i
\rangle_s| \rightarrow0$
and $\mathbb{E}^{\pi^N} |\langle\phi_i, E^N \phi_j \rangle_s|
\rightarrow0$ as $N \rightarrow\infty$, for any pair of indices
$i,j$ and
for $s$ appearing in Assumptions~\ref{ass1}.
\end{prop}

Thus the discrete time Markov chain $\{x^k\}$ obtained by the successive
accepted samples of the RWM algorithm has approximately
the expected drift and covariance structure of the SPDE~(\ref{eqnspde}).\vadjust{\goodbreak}
It is
also crucial to our subsequent argument involving the martingale
central limit theorem that the error terms $r^N$
and $E^N$ converge to zero in the Hilbert space $\h^s$
norm and inner-product as stated.


With this\vspace*{1pt} in hand, we need to establish the appropriate
invariance principle to show that the dynamics of the Markov chain
$\{x^k\}$, when seen as the values of a~continuous
time process on a~time mesh with steps of $O(1/N)$,
converges weakly to the law of the SPDE given in
(\ref{eqnspde}) on $C([0,T],\h^s)$. To this end we define,
for $k \geq0$,
%
%
\begin{eqnarray}
\label{eqnmeandrandgamm}
m^N(\ccdot) &\eqdef& P^N(\ccdot) + C^N\nabla\Psi^N (\ccdot)
\Gamma^{k+1,N} \nonumber\\[-8pt]\\[-8pt]
&\eqdef&{\sqrt{\frac{N}{2\ell^2\beta}}} \bigl(x^{k+1}
- x^k - \mathbb{E}_k(x^{k+1} - x^k)\bigr),\nonumber\\
\label{eqnRN}
r^{k+1,N} &\eqdef& N \mathbb{E}_k(x^{k+1} - x^{k}) +
\ell^2\beta\bigl(P^N x^k + C^N \nabla\Psi^N(x^k)\bigr), \\
\label{eqnEN}
E^{k+1,N} &\eqdef& N \mathbb{E}_k[(x^{k+1} - x^{k})\otimes
(x^{k+1} - x^{k})] -
2\ell^2 \beta C^N
\end{eqnarray}
with $E^{0,N}, \Gamma^{0,N}, r^{0,N} = 0$.
Notice that for fixed~$N$, $\{r^{k,N}\}_{k \geq1}, \{E^{k,N}\}_{k\geq
1}$ are, since $x^0 \sim\pi^N$,
stationary sequences.

By definition,
%
%
\begin{equation}
\label{eqnmexp11}
x^{k+1} = x^k + \mathbb{E}_k(x^{k+1} - x^k) +
{\sqrt{{\frac{2\ell^2\beta}{N}}}} \Gamma^{k+1,N} .
\end{equation}
From~(\ref{thmDD-drift}) in Proposition
\ref{thmdrift-diffus}, for large enough~$N$,
%
%
\begin{eqnarray}\label{eqnEscheme0}
x^{k+1} &\approx& x^k
- {\frac{\ell^2 \beta}{N}} \bigl(P^N x^k + C^N\nabla\Psi^N
(x^k)\bigr)
+
{\sqrt{\frac{2\ell^2\beta}{N}}} \Gamma^{k+1,N}\nonumber\\[-8pt]\\[-8pt]
&=& x^k - {\frac{\ell^2 \beta}{N}} m^N(x^k) +
{\sqrt{\frac{2\ell^2\beta}{N}}} \Gamma^{k+1,N}.\nonumber
\end{eqnarray}
From the definition of $\Gamma^{k,N}$ in~(\ref{eqnmeandrandgamm}),
and from~(\ref{thmDD-diffus})
in Proposition~\ref{thmdrift-diffus},
\[
\mathbb{E}_k(\Gamma^{k+1,N} ) = 0 \quad\mbox{and}\quad
\mathbb{E}_k (\Gamma^{k+1,N} \otimes\Gamma^{k+1,N} ) \approx C^N .
\]
%
Therefore, for large enough~$N$, equation~(\ref{eqnEscheme0})
``resembles'' the
Euler scheme for
simulating the finite-dimensional approximation of the SPDE~(\ref{eqnspde})
on~$\mathbb{R}^N$, with drift function $m^N(\ccdot)$ and
covariance operator $C^N$:
\[
x^{k+1} \approx x^k - h(\ell) m^N(x^k)\Delta t +
\sqrt{2h(\ell) \Delta t} \Gamma^{k+1,N}
\qquad\mbox{where } \Delta t \eqdef\frac{1}{N}.
\]
%
This is the key idea underlying our main result (Theorem \ref
{thmmain}): the Markov chain~(\ref{eqnRWMfancy})
looks like a~weak Euler approximation of~(\ref{eqnspde}).

Note that there is an important difference in analyzing the weak
convergence from the traditional Euler scheme. In our case, for any
fixed $N \in\mathbb{N}$, $\Gamma^{k,N} \in X^N$ is finite
dimensional, but clearly the dimension of $\Gamma^{k,N} $ grows with
$N$. Also,\vadjust{\goodbreak} the distribution of the initial condition $x(0) \dist\pi^N$
changes\vspace*{1pt} with~$N$, unlike the case of the traditional Euler scheme
where the distribution of $x(0)$ does not change with~$N$. Moreover,
for any fixed~$N$, the ``noise'' process $\{\Gamma^{k,N} \}$ are not
formed of
independent random variables. However, they are identically
distributed (a stationary sequence) because the Metropolis algorithm
preserves stationarity. To obtain an invariance principle, we first
use a~version of the martingale\vspace*{1pt} central limit theorem
(Proposition~\ref{thmBergmclt}) to show that the noise process
$\{\Gamma^{k,N} \}$, when rescaled and summed, converges weakly to a~Brownian motion on $C([0,T],\h^s)$ with covariance operator $C_s$,
for any $T = \mathcal{O}(1)$. We then use continuity of an appropriate
It\^{o} map
to deduce the desired result.

Before we proceed, we introduce some notation. Fix $T > 0$, and
define
%
%
\begin{equation}\label{eqndeltaconst}
\Delta t\eqdef1/N ,\qquad
t^k\eqdef k \Delta t ,\qquad
\eta^{k,N} \eqdef\sqrt{\Delta t} \sum_{l = 1}^{k} \Gamma^{l,N}
\end{equation}
and
%
%
\begin{equation}\label{eqnwnCproc}
W^N(t) \eqdef\eta^{\lfloor{Nt\rfloor},N} +
\frac{Nt -\lfloor Nt \rfloor}{\sqrt{N}}
\Gamma^{\lfloor{Nt\rfloor}+1,N},\qquad t \in[0,T] .
\end{equation}
Let $W(t), t \in[0,T]$ be an $\h^s$ valued Brownian motion with
covariance operator~$C_s$. Using a~martingale
central limit theorem, we will prove
the following proposition in Section~\ref{secinvp}.
\begin{prop}
\label{lembweakconv}%
Let Assumptions~\ref{ass1} (below) hold.
Let $x^0 \sim\pi^N$. The process $W^N(t)$ defined in (\ref
{eqnwnCproc}) converges
weakly to $W$ in $C([0,T],\h^s)$ as~$N$ tends to $\infty$, where $W$ is
a Brownian motion in time with covariance operator $C_s$ in $\h^s$
and $s$ is defined in Assumptions~\ref{ass1}.
Furthermore, the pair $(x^0,W^N(t))$ converges weakly to $(z^0,W)$ where
$z^0 \sim\pi$ and
Brownian motion~$W$ is independent of the initial condition $z^0$
almost surely.
\end{prop}

Using this invariance principle for the noise process and
the fact that
the noise process is additive (the diffusion coefficient is constant),
the invariance principle for the Markov chain follows from a~continuous mapping argument which we now outline.
For any $(z^0,W) \in\mathcal{H}^s \times C([0,T]; \mathcal{H}^s)$,
we define the It\^{o} map
$\Theta\colon\mathcal{H}^s \times C([0,T];\mathcal{H}^s)
\rightarrow C([0,T];\mathcal{H}^s)$ by
$\Theta\dvtx(z^0,W) \mapsto z$ where $z$ solves
%
%
\begin{equation} \label{eqnpcxN1jcm}
z(t) = z^0 - h(\ell) \int_0^t \bigl(z(s) + C
\nabla\Psi(z(s)) \bigr) \,ds +
\sqrt{2h(\ell)} W(t)
\end{equation}
for all $t \in[0,T]$ and $h(\ell) = \ell^2 \beta$ is as defined in
(\ref{eqnhl}). Thus $z=\Theta(z^0,W)$ solves the SPDE
(\ref{eqnspde}) with $h(\ell) = \ell^2 \beta$.
We will see in Lemma~\ref{lemcontmap} that
$\Theta$ is a~continuous map from $\mathcal{H}^s \times
C([0,T]; \mathcal{H}^s)$ into
$C([0,T]; \mathcal{H}^s)$.

We now define the piecewise constant interpolant of $x^k$,
%
%
\begin{equation}\label{ebarZdef}
\bar{z}^N(t) = x^k \qquad\mbox{for } t \in[t^k, t^{k+1}).\vadjust{\goodbreak}
\end{equation}
Set
%
%
\begin{equation}\label{eqnDN}
d^N(x) \eqdef N \mathbb{E}_0(x^{1} - x).
\end{equation}
Note that $d_N(x) \approx-h(\ell) m_N(x)$.
We can use $\bar{z}^N$ to construct a~continuous piecewise
linear interpolant of $x^k$ by defining
%
%
\begin{equation}\label{eqnprocessxN1}
z^N(t) = z^0 + \int_0^t d^N(\bar{z}^N(s)) \,ds +
\sqrt{2h(\ell)} W^N(t).
\end{equation}
Notice that $d^N(x)$ defined in~(\ref{eqnDN}) is a~function
which depends on arbitrary $x = x^0$
and averages out the randomness in $x^1$ conditional on fixing $x =
x^0$. We may then evaluate this
function at any $x\in\mathcal{H}^s$ and, in particular,
at $\bar{z}^N(s)$ as in~(\ref{eqnprocessxN1}).
Use of the stationarity of the sequence $x^k$,
together with equations~(\ref{eqnmexp11}),
(\ref{eqndeltaconst}) and~(\ref{eqnwnCproc}),
reveals that the definition~(\ref{eqnprocessxN1})
coincides with that given in~(\ref{eqnMCMCe}).
Using the closeness of $d^N$ and $-h(\ell) m^N$,
of $z^N$ and ${\bar z}^N$ and of $m^N$ and the
desired limiting
drift, we will see that there exists a~$\widehat W^N \weak W$ as $N\rightarrow\infty$,
such that
%
%
\begin{equation}\label{eqnprocessxN2}\qquad
z^N(t) = z^0 - h(\ell) \int_0^t \bigl(z^N(s) + C
\nabla\Psi(z^N(s)) \bigr) \,ds +
\sqrt{2h(\ell)} \widehat W^N(t),
\end{equation}
so that $z^N=\Theta(z^0,\widehat W^N)$.
By the continuity of $\Theta$ we will show, using
the continuous mapping theorem, that
%
%
\begin{equation}\label{eqnZweakConv}
z^N=\Theta(z^0,\widehat W^N) \quad\longweak\quad z=\Theta(z^0,W)
\qquad\mbox{as } N \rightarrow\infty.
\end{equation}
It will be important to show that the weak limit of
$(z^0,\widehat W^N)$, namely $(z^0,W)$, comprises of
two independent random variables $z^0$ (from the
stationary distribution) and $W$.

The weak convergence in~(\ref{eqnZweakConv}) is the principal result of
this article and is stated precisely in Theorem~\ref{thmmain}. To
summarize, we have argued that the RWM is well
approximated by an Euler approximation of~(\ref{eqnspde}). The Euler
approximation itself can be seen as an approximate solution of~(\ref{eqnspde})
with a~modified Brownian motion. As $N \rightarrow\infty$, all
approximation errors go to zero in the
appropriate sense and one deduces that the RWM algorithm converges to
the solution of~(\ref{eqnspde}).
%

\subsection{A framework for expected drift and diffusion}
We now turn to the question of how the RWM algorithm
produces the appropriate drift and covariance encapsulated
in Proposition~\ref{thmdrift-diffus}. This result,
which shows that the algorithm (approximately) performs
a noisy steepest ascent process, is
at the heart of why the Metropolis algorithm works.
In the rest of this section we set up a~framework
which will be used for deriving the expected drift and diffusion terms.


Recall the setup from Section~\ref{secrwm}. Starting from
(\ref{eqnaccrem}), after some algebra we obtain
%
%
\begin{equation}\label{eqnRx,y1}
Q(x,\xi) = -{\sqrt{\frac{2\ell^2}{N}}}\langle\zeta,\xi\rangle-
{\frac{\ell^2}{N}} \|\xi\|^2 - r(x,\xi),
\end{equation}
where we have defined
%
%
\begin{eqnarray}\label{eqnzetalt}
\zeta& \eqdef &C^{-1/2}(P^N x) + C^{1/2} \nabla\Psi^N(x),\\
r(x,\xi) & \eqdef&\Psi^N(y) - \Psi^N(x) - \langle\nabla
\Psi^N(x), P^Ny - P^N x\rangle.
\end{eqnarray}

\begin{remark}
If $x\dist\pi_0$ in $\h^s$, then the random variable $C^{-1/2}x$ is
not well defined in $\h^s$ because $C^{-1/2}$ is not a~trace class
operator. However, equation~(\ref{eqnzetalt}) is still well defined
because the operator $C^{-1/2}$ acts only in~$X^N$ for any
\textit{fixed}~$N$. 
\end{remark}

Notice that $C^{1/2}\zeta$ is approximately the
drift term in the SPDE~(\ref{eqnspde}) and this
plays a~key role in obtaining the
mean drift from the accept/reject mechanism; this
point is elaborated on in the arguments leading up to
(\ref{eqlambdaZeta}).
By~(\ref{eqnpsidiff}) and Assumptions
\ref{ass1},~\ref{ass2} on~$\Psi$ and $\Psi^N$ below, we will
obtain a~global
bound on the remainder term of the form
%
%
\begin{equation}\label{eqnrx,y,xi}
|r(x,\xi)|\leq M {\frac{\ell^2}{N}} \|C^{1/2}
\xi\|^2_s .
\end{equation}
Because of our assumptions on~$C$ in~(\ref{eqndecayeigen}), the
moments of $\|C^{1/2} \xi\|^2_s$ stay uniformly bounded as
$N\rightarrow
\infty$. Hence, we will neglect this term to explain the heuristic
ideas. Since
$\xi= \sum_{i=1}^N \xi_i\phi_i$ with $\xi_i \dist\Normal(0,1)$,
we find that
for fixed $x$,
%
%
\begin{equation}\label{eqnweakconvQ}
Q(x,\xi) \approx\Normal\biggl(-\ell^2, 2 \ell^2
{\frac{\|\zeta\|^2}{N}}\biggr)
\end{equation}
for large~$N$ (see Lemma~\ref{lemestRx}). Since $x \dist\pi$, we
have that $C^{-1/2}(P^Nx) =
\sum_{k=1}^N \rho_j \phi_j$, where $\rho_j$ are i.i.d.
$\Normal(0,1)$. Much as with the term $r(x,\xi)$ above,
the second term in expression~(\ref{eqnzetalt}) for $\zeta$ can
be seen as a~perturbation term which
is small in magnitude compared to the first term in
(\ref{eqnzetalt}) as $N \rightarrow\infty$. Thus,
as shown in Lemma~\ref{lemshiftmeas}, we have ${\|\zeta\|^2}/{N}
\rightarrow1$ for~$\pi$-a.e. $\zeta$
as $N\rightarrow\infty$. Returning to~(\ref{eqnweakconvQ}), this suggests
that it is reasonable for~$N$ sufficiently large to make the approximation
%
%
\begin{equation} \label{eqnweakconvQ2}
Q(x,\xi) \approx\Normal(-\ell^2, 2 \ell^2), \qquad\pi\mbox{-a.s.}
\end{equation}
Much of this section is concerned with understanding the behavior of
one step of the
RWM algorithm if we make the approximation in
(\ref{eqnweakconvQ2}). Once this is understood, we will retrace our
steps being more careful to control the approximation error leading to
(\ref{eqnweakconvQ2}).\vadjust{\goodbreak}

The following lemma concerning normal
random variables will be critical to identifying the source of the
observed drift. It gives us the relation between the constants
in the expected drift
and diffusion coefficients which ensures~$\pi$ invariance, as will be
seen later in this section.
\begin{lemma}\label{lemmagic} Let $Z_\ell\dist
\Normal(-\ell^2,2\ell^2)$. Then $\mathbb{P}(Z_\ell> 0) =
\mathbb{E}(e^{Z_\ell}1_{Z_{\ell}< 0}) =\break
\Phi(-{\ell}/{\sqrt{2}})$ and
%
%
\begin{equation}\label{eqnmagic}
\mathbb{E}(1 \wedge e^{Z_{\ell}}) =
2\Phi\bigl(-{\ell}/{\sqrt{2}}\bigr) = \beta.
\end{equation}
Furthermore, if $z \dist\Normal(0,1)$ then
%
%
\begin{equation}\label{lemnormmo}
\mathbb{E}[z(1 \wedge e^{az +b})] =
a\exp(a^2/2 + b) \Phi\biggl(-\frac{b}{|a|} - |a|\biggr)
\end{equation}
for any real constants $a$ and $b$.
\end{lemma}
\begin{pf}
A straightforward calculation. See Lemma 2 in~\cite{Besketal0101}.
\end{pf}

The calculations of the expected one step drift and diffusion
needed to prove Proposition~\ref{thmdrift-diffus} are long
and technical. In order to enhance the readability, in the next two
sections we outline
our proof strategy emphasizing the key calculations.

\subsection{Heuristic argument for the expected drift}\label{secdiftHuristic}
In this section, we will give heuristic arguments which underly
(\ref{thmDD-drift}) from Proposition~\ref{thmdrift-diffus}. Recall that
$\{\phi_1,\phi_2,\ldots\}$ is an orthonormal basis for $\h$. Let
$x^k_i, i \leq N$, denote the $i$th coordinate of $x^k$ and
$C^N$ denote the covariance operator on $X^N$, the span of
$\{\phi_1,\phi_2,\ldots,\phi_N\}$. Also recall that $\mathcal{F}_k$
denotes the sigma algebra generated by $\{x^n, \xi^n, \gamma^n, n
\leq k\}$
and the conditional expectations $\mathbb{E}(\cdot| \mathcal{F}_k)$
are denoted by $\mathbb{E}_k(\ccdot)$. Thus $\mathbb{E}_0(\cdot)$
denotes the expectation with respect to $\xi^1$ and $\gamma^1$ with~$x^0$
fixed. Also, for notational convenience, set $x^0 = x$ and $\xi^1
= \xi$. Letting $\mathbb{E}_0^\xi$ denote the expectation with
respect to~$\xi$, it follows that
%
%
\begin{eqnarray}\label{eqnexpQ}\quad
N \mathbb{E}_0(x^{1}_i - x^0_i) &=& N \mathbb{E}_0\bigl(\gamma^1
(y^1_i - x_i)\bigr)\nonumber\\
&=& N \mathbb{E}_0^\xi\Biggl(\alpha(x,\xi) {\sqrt{\frac{2\ell
^2}{N}}}(C^{1/2}\xi)_i\Biggr) \nonumber\\[-8pt]\\[-8pt]
&=& \lambda_i \sqrt{2\ell^2 N} \mathbb{E}_0^\xi
(\alpha(x,\xi) \xi_i) \nonumber\\
&=&\lambda_i \sqrt{2\ell^2 N}
\mathbb{E}_0^\xi\bigl(\bigl(1 \wedge
e^{Q(x,\xi)}\bigr)\xi_i\bigr) .\nonumber
\end{eqnarray}
To approximately evaluate~(\ref{eqnexpQ}) using Lemma~\ref{lemmagic},
it is easier to first factor $Q(x,\xi)$ into components involving
$\xi_i$ and those orthogonal (under $\mathbb{E}^\xi_0$) to them. To
this end we introduce the following
terms:
%
%
\begin{eqnarray}
\label{eqnR}
R(x,\xi) &\eqdef& -{\sqrt{\frac{2\ell^2}{N}}}\sum_{j=1}^N \zeta_j
\xi_j -
\frac{\ell^2}{N} \sum_{j=1}^N \xi_{j}^2, \\
\label{eqnRi}
R_i(x,\xi) &\eqdef&
-{\sqrt{\frac{2\ell^2}{N}}}\sum_{j=1, j \neq i}^N \zeta_j \xi_j -
\frac{\ell^2}{N}\sum_{j=1, j \neq i}^N \xi_{j}^2 .
\end{eqnarray}
Hence, for large~$N$ (see Lemma~\ref{lemdo1}),
%
%
\begin{eqnarray}\label{eqnQRiapprox}\quad
Q(x,\xi)&=& R(x,\xi) - r(x,\xi) = R_i(x,\xi) -
{\sqrt{\frac{2\ell^2}{N}}} \zeta_i\xi_i -{\frac{\ell^2}{N}}\xi_i^2
-r(x,\xi)\nonumber\\
&=& R_i(x,\xi) - {\sqrt{\frac{2\ell^2}{N}}} \zeta_i \xi_i +
O\biggl(\frac
{1}{N}\biggr)\\
&\approx& R_i(x,\xi) - {\sqrt{\frac{2\ell^2}{N}}}
\zeta_i\xi_i .\nonumber
\end{eqnarray}
The important observation here is that conditional on $x$, the random
variable $R_i(x,\xi)$ is independent of $\xi_i$. Hence, the expectation
$\mathbb{E}_0^\xi((1 \wedge e^{Q(x,\xi)})\xi_i)$ can be
computed by first computing it over $\xi_i$ and then over
$\xi\setminus\xi_i$. Let $\mathbb{E}^{\xi_i^{-}},
\mathbb{E}^{\xi_i}$ denote the expectation with respect to $\xi
\setminus\xi_i, \xi_i$, respectively. Using the relation
(\ref{eqnQRiapprox}), and applying~(\ref{lemnormmo}) with
$a =
-\sqrt{\frac{2\ell^2}{N}}\zeta_i$, $z = \xi_i$ and $b = R_i(x,\xi)$,
we obtain (see Lemma~\ref{lemdo2})
%
%
\begin{eqnarray}\label{eqnpenexes1}
&&\mathbb{E}_0^{\xi} \bigl(\bigl(1 \wedge e^{Q(x,\xi)}\bigr)\xi_i\bigr)\nonumber\\
&&\qquad
\approx
-{\sqrt{\frac{2\ell^2}{N}}}\zeta_i \mathbb{E}_0^{\xi_i^-}
e^{R_i(x,\xi) + ({\ell^2/N}) \zeta^2_i} \Phi
\Biggl({ \frac{-R_i(x, \xi)}
{{\sqrt{{2\ell^2}/{N}}}|\zeta_i|}} -
{\sqrt{\frac{2\ell^2}{N}}}|\zeta_i|\Biggr) \\
&&\qquad \approx
-{\sqrt{\frac{2\ell^2}{N}}}\zeta_i \mathbb{E}_0^{\xi_i^-}
e^{R_i(x,\xi)+{\ell^2}/{N} \zeta^2_i}
\Phi\biggl({ \frac{-R_i(x, \xi)}
{\sqrt{{2\ell^2}/{N}}|\zeta_i|}}\biggr).\nonumber
\end{eqnarray}
Now, again from the relation~(\ref{eqnQRiapprox}) and the
approximation $Q(x,\xi)$ encapsulated in~(\ref{eqnweakconvQ}),
it follows that for sufficiently large~$N$
%
%
\begin{equation}\label{eqnweakconvR}
R_i(x,\xi) \approx\Normal(-\ell^2,2\ell^2),\qquad \pi\mbox{-a.s.}
\end{equation}
Combining~(\ref{eqnpenexes1}) with the fact that, for large enough
$N$, $\Phi({-R_i(x, \xi)}/\break
{\sqrt{\frac{2\ell^2}{N}}|\zeta_i|}) \approx1_{R_i(x,\xi) < 0}$,
we see that Lemma~\ref{lemmagic} implies that (see Lemmas~\ref
{lemdo3}--\ref{lemdo4})
%
%
\begin{eqnarray}
\label{eqnRilimwbet1}
\mathbb{E}_0^{\xi_i^-}
e^{R_i(x,\xi)+{\ell^2}/{N} \zeta^2_i}
\Phi\biggl({ \frac{-R_i(x, \xi)}
{\sqrt{{2\ell^2}/{N}}|\zeta_i|}}\biggr) 
&\approx&\mathbb{E}_0^{\xi_i^-}
\bigl(e^{R_i(x,\xi)} 1_{R_i(x,\xi) < 0}\bigr)
\\
\label{eqnRilimwbet2}
&\approx&
\mathbb{E} e^{Z_\ell} 1_{Z_\ell< 0} = \beta/2,
\end{eqnarray}
where $Z_\ell\dist N(-\ell^2, 2\ell^2)$.
Hence, from~(\ref{eqnexpQ}),~(\ref{eqnpenexes1}) and
(\ref{eqnRilimwbet2}), we gather that for large~$N$,
\[
N \mathbb{E}_0(x^{1}_i - x^0_i) \approx-\ell^2\beta\lambda_i
\zeta_i .
\]
To identify the drift, observe that since $C^{-1/2}$ is self-adjoint
and $i \leq N$, we have $\lambda_i C^{-1/2}\phi_i = \phi_i$ and
%
%
\begin{eqnarray}\label{eqlambdaZeta}
\lambda_i \zeta_i
&=& \lambda_i \langle C^{-1/2} (P^Nx) +
C^{1/2}\nabla\Psi^N(x) , \phi_i \rangle\nonumber\\
&=& \lambda_i
\langle C^{-1/2} (P^Nx) +
C^{-1/2}C\nabla\Psi^N( x) , \phi_i \rangle\\
&=&\langle P^Nx + C^N\nabla\Psi^N( x) , \phi_i \rangle.\nonumber
\end{eqnarray}
Hence, for large enough~$N$, we deduce that (heuristically)
the expected drift in the $i$th coordinate after one step
of the Markov chain $\{x^k\}$ is well approximated by the
expression
\[
N \mathbb{E}_0(x^{1}_i - x^0_i) \approx-\ell^2\beta\bigl(P^N x + C^N
\nabla\Psi^N(x)\bigr)_i.
\]
This is an approximation of
the drift term that appears in the SPDE
(\ref{eqnspde}). Therefore, the above heuristic arguments
show how the
Metropolis algorithm achieves the ``change of measure'' by mapping
$\pi_0$ to~$\pi$. The above arguments can be made rigorous by
quantitatively controlling the errors made. In
Section~\ref{secstec}, we quantify the size of the neglected terms
and quantify the rate at which $Q$ is well approximated by a~Gaussian
distribution. Using these estimates, in
Section~\ref{secdrifes} we will retrace the arguments of this section
paying attention to the cumulative error, thereby proving
(\ref{thmDD-drift}) of Proposition~\ref{thmdrift-diffus}.


\subsection{Heuristic argument for the expected diffusion coefficient}\label{secdiffusionHuristic}
We now give the heuristic arguments for the expected diffusion
coefficient, after one step of the Markov chain $\{x^k\}$. The
arguments used here are much simpler than the drift calculations.
The strategy is the same as in the drift case except that now we
consider the
covariance between two coordinates $x^1_i$ and~$x^1_j$. For $1 \leq
i,j \leq N$,
%
%
\begin{eqnarray}\label{eqndecoptrc}
&&
N \mathbb{E}_0[(x_i^{1} - x_i^{0})(x_j^{1} - x_j^{0})]\nonumber\\
&&\qquad= N \mathbb{E}^\xi_0[(y_i^{1} - x_i)(y_j^{1} - x_j) \alpha
(x,\xi)]\nonumber\\[-8pt]\\[-8pt]
&&\qquad= N \mathbb{E}_0^\xi\bigl[(y_i^{1} - x_i)(y_j^{1} - x_j)
\bigl(1 \wedge\exp Q(x,\xi)\bigr)\bigr]\nonumber\\
&&\qquad= 2\ell^2 \mathbb{E}_0^{\xi}\bigl[(C^{1/2}\xi)_i (C^{1/2}\xi)_j
\bigl(1 \wedge\exp Q(x,\xi)\bigr)\bigr] .\nonumber
\end{eqnarray}
Now notice that
\[
\mathbb{E}_0^{\xi}[(C^{1/2}\xi)_i (C^{1/2}\xi)_j] =
\lambda_i\lambda_j \delta_{ij} ,
\]
where $\delta_{ij} = 1_{i=j}$.
Similar to the calculations used when evaluating the expected drift, we define
%
%
\begin{equation} \label{eqnRij}
R_{ij}(x,\xi) \eqdef
-{\sqrt{\frac{2\ell^2}{N}}}\sum_{k=1, k \neq{i,j}}^N \zeta_k \xi
_k -
{\frac{\ell^2}{N}}\sum_{k=1, k \neq{i,j}}^N \xi_{k}^2
\end{equation}
and observe that
\[
R(x,\xi) = R_{ij}(x,\xi) - {\sqrt{\frac{2\ell^2}{N}}} \zeta_i\xi
_i -{\frac{\ell^2}{N}} \xi_i^2 -
{\sqrt{\frac{2\ell^2}{N}}} \zeta_j\xi_j -{\frac{\ell^2}{N}} \xi
_j^2 .
\]
Hence, for sufficiently large~$N$, we have $Q(x,\xi)\approx
R_{ij}(x,\xi)$. By replacing $Q(x,\xi)$ in
(\ref{eqndecoptrc}) by $R_{ij}(x,\xi)$ we can take advantage of the
fact that $R_{ij}(x,\xi)$ is conditionally independent of
$\xi_i,\xi_j$. However,\vspace*{-1pt} the additional error term introduced is easy to
estimate because the function $f(x) \eqdef(1 \wedge e^x)$ is
1-Lipschitz. So, for large enough~$N$ (Lemma~\ref{lemdrifvar}),
%
%
\begin{eqnarray} \label{eqnvardiffbd}
&&\mathbb{E}_0^{\xi}\bigl[(C^{1/2}\xi)_i (C^{1/2}\xi)_j\bigl(1 \wedge
\exp Q(x,\xi)\bigr)\bigr]\nonumber\\
&&\qquad \approx
\mathbb{E}_0^{\xi}\bigl[(C^{1/2}\xi)_i (C^{1/2}\xi)_j
\bigl(1 \wedge\exp R_{ij}(x,\xi)\bigr)\bigr] \\
&&\qquad= \lambda_i\lambda_j \delta_{ij} \mathbb{E}_0^{\xi^-_{ij}}\bigl[
\bigl(1 \wedge\exp R_{ij}(x,\xi)\bigr)\bigr] .\nonumber
\end{eqnarray}
Again, as in the drift calculation, we have that
\[
R_{ij}(x,\xi) \longweak\Normal(-\ell^2,2\ell^2),\qquad \mbox{$\pi
$-a.s.}
\]
So by the dominated convergence theorem and Lemma~\ref{lemmagic},
%
%
\begin{equation}\label{eqnrconvbet}
\lim_{N\rightarrow\infty} \mathbb{E}^{\xi^-_{ij}}\bigl[
\bigl(1 \wedge\exp R_{ij}(x,\xi)\bigr)\bigr] = \beta.
\end{equation}
Therefore, for large~$N$,
\[
N \mathbb{E}_0[(x_i^{1} - x_i^{0})(x_j^{1} - x_j^{0})]
\approx
2\ell^2\beta\lambda_i \lambda_j \delta_{ij}=2\ell^2\beta\langle
\phi_i, C \phi_j\rangle
\]
or in other words,
\[
N \mathbb{E}_0[(x^{1} - x^{0})\otimes(x^{1} - x^{0})]
\approx
2\ell^2 \beta C^N .
\]
%
As with the drift calculations in the last section, these calculations
can be made rigorous by tracking the size of the neglected terms and
quantifying the rate at which $Q$ is approximated by the appropriate
Gaussian. We will substantiate these arguments Section~\ref{secdiffes}.


\section{Main theorem} \label{secmainthm}

In this section we state the assumptions we make on~$\pi_0$ and~$\Psi$ and then prove our main
theorem.

\subsection{\texorpdfstring{Assumptions on~$\Psi$ and~$C$}{Assumptions on Psi and C}}\label{secpsin}
%

The assumptions we make now concern (i) the rate of
decay of the standard deviations in the prior
or reference measure $\pi_0$ and (ii)
the properties of the Radon--Nikodym derivative (likelihood function).\vadjust{\goodbreak}
These assumptions are naturally
linked; in order for~$\pi$ to be well defined
we require that~$\Psi$ is $\pi_0$-measurable
and this can be achieved by ensuring that
$\Psi$ is continuous on a~space which has
full measure under $\pi_0$. In fact, in a~wide
range of applications,~$\Psi$ is Lipschitz
on such a~space~\cite{Stua10}.
In this paper we require, in addition,
that~$\Psi$ be twice differentiable
in order to define the diffusion limit. This, too,
may be established in many applications.
To avoid technicalities, we assume that $\Psi(x)$ is quadratically
bounded, with first derivative linearly bounded and second
derivative globally bounded. A~simple example of a~function
$\Psi$ satisfying the above
assumptions is $\Psi(x) = \|x\|_s^2$.
\begin{assumptions}
\label{ass1}
The operator~$C$ and functional~$\Psi$ satisfy the following:
\begin{longlist}[(1)]
\item[(1)]\textit{Decay of eigenvalues} $\lambda_i^2$ of~$C$:
There exist $M_{-}, M_+ \in(0,\infty)$ and $\kappa> \frac{1}{2}$
such that
%
%
\begin{equation}\label{eqndecayeigen}
M_{-} \leq i^\kappa\lambda_i \leq M_+\qquad \forall i \in\mathbb{Z}_+ .
\end{equation}
\item[(2)]\textit{Assumptions on}~$\Psi$:
There exist constants $M_i \in\mathbb{R}, i \leq4$
and $s \in[0,\break \kappa- 1/2)$ such that
%
%
\begin{eqnarray}
\label{eqnpsi1}
M_1 &\leq&\Psi(x) \leq M_2(1 + \|x\|_s^2) \qquad\forall x \in\h^s,\\
\label{eqnpsi2}
\| \nabla\Psi(x)\|_{-s} &\leq& M_3(1 + \|x\|_s) \qquad\forall x
\in\h^s,\\
\label{eqnpsi3}
\|\partial^2 \Psi(x)\|_{\mathcal{L}(\mathcal{H}^{s},\mathcal
{H}^{-s})} & \leq& M_4 \qquad\forall x
\in
\h^s.
\end{eqnarray}
\end{longlist}
\end{assumptions}

Notice also that the above assumptions on~$\Psi$
imply that for all $x,y \in\h^s$,
%
%
\begin{subequation}
\label{eqnpsidiff}
\begin{eqnarray}
\label{eqnpsidiffa}
|\Psi(x) - \Psi(y)| &\leq& M_5 ( 1 + \|x\|_s + \|y\|_s ) \|x
- y\|_s, \\
\label{eqnpsidiffb}
\Psi(y) &=& \Psi(x) + \langle\nabla\Psi(x), y-x \rangle+
\operatorname{rem}(x,y), \\
\label{eqnpsidiffc}
\operatorname{rem}(x,y) &\leq& M_6 \|x - y \|_s^2
\end{eqnarray}
\end{subequation}
for some constants $M_5, M_6 \in\mathbb{R}_+$.
\begin{remark}
The condition $\kappa>\frac12$ ensures
that the covariance operator for $\pi_0$ is trace class.
In fact, the $\h^r$ norm of a~realization of a~Gaussian measure
$\Normal(0,C)$ defined on $\h$ is almost surely finite
if and only if $r < \kappa- \frac{1}{2}$~\cite{DaprZaby92}.
Thus the choice of Sobolev space $\mathcal{H}^s$, with $s \in
[0,\kappa-\frac12)$ in which we
state the assumptions on~$\Psi$, is made to ensure
that the Radon--Nikodym derivative of~$\pi$ with respect
to $\pi_0$ is well defined. Indeed, under our
assumptions,~$\Psi$ is Lipschitz continuous on a~set
of full $\pi_0$ measure; it is hence $\pi_0$-measurable.
Weaker growth
assumptions on~$\Psi$, its Lipschitz constant
and second derivative
could be dealt with by use of stopping time
arguments.
\end{remark}

The following lemma will be used repeatedly.\vadjust{\goodbreak}
\begin{lemma} \label{lneeded}
Under Assumptions~\ref{ass1} it follows that, for
all $a \in\bbR$,
%
%
\begin{equation}
\label{eqequiv}
\|C^{a}x\| \asymp\|x\|_{-2\kappa a}.
\end{equation}
Furthermore, the function
$C\nabla\Psi\dvtx\h^s \to\h^s$ is
globally Lipschitz.
\end{lemma}
\begin{pf} The first result follows from the inequality
\[
\|C^{a}x\|^2=\sum_{j=1}^{\infty}\lambda_{j}^{4a}x_j^2
\le M_{+}\sum_{j=1}^{\infty}j^{-4a\kappa}x_j^2=M_+\|x\|_{-2\kappa a}^2,
\]
and a~similar lower bound, using~(\ref{eqndecayeigen}).
To prove the global Lipschitz property we first note that
%
%
\begin{eqnarray}\label{eqsure}
\nabla\Psi(u_1)-\nabla\Psi(u_2)&=&K(u_1-u_2)\nonumber\\[-8pt]\\[-8pt]
:\!&=&\int_0^1
\partial^2\Psi\bigl(tu_1
+(1-t)u_2\bigr)\, dt (u_1-u_2).\nonumber
\end{eqnarray}
Note that $\|K\|_{\mathcal{L}(\h^s,\h^{-s})} \le M_4$
by~(\ref{eqnpsi3}). Thus,
%
\begin{eqnarray*}
&&\bigl\|C\bigl(\nabla\Psi(u_1)-\nabla\Psi(u_2)\bigr)\bigr\|_{s}\\
&&\qquad\le M\|C^{1-s/2\kappa}K(u_1-u_2)\|\\
&&\qquad\le M\|C^{1-s/2\kappa}KC^{s/2k}C^{-s/2k}(u_1-u_2)\|\\
&&\qquad\le M\|C^{1-s/2\kappa}KC^{s/2k}\|_{\mathcal{L}(\h,\h)} \|u_1-u_2\|
_{s}\\
&&\qquad\le M\|C^{1-s/2\kappa}\|_{\mathcal{L}(\h^{-s},\h)}
\|K\|_{\mathcal{L}(\h^s,\h^{-s})}
\|C^{s/2k}\|_{\mathcal{L}(\h,\h^s)}\|u_1-u_2\|_{s}.
\end{eqnarray*}
The three linear operators are bounded between the
appropriate spaces, in the case of $C^{1-s/2\kappa}$ by using
the fact that $s<\kappa-\frac12$ implies $s<\kappa$.
\end{pf}

\subsection{Finite-dimensional approximation of the invariant distribution}

For simplicity we assume throughout this paper that
$\Psi^N(\cdot) = \Psi(P^N \cdot)$.
We note again that $\nabla\Psi^N(x)=P^N\nabla\Psi(P^N x)$
and $\partial^2\Psi^N(x)=P^N \partial^2\Psi(P^N x) P^N$.
Other approximations could be handled similarly.
The function $\Psi^N$ may
be shown to satisfy the following.
\begin{assumptions}[(Assumptions on $\Psi^N$)]
\label{ass2}
The functions $\Psi^N \dvtx X^N \mapsto\mathbb{R}$ satisfy the same
conditions imposed on~$\Psi$ given by
equations~(\ref{eqnpsi1}),~(\ref{eqnpsi2}) and~(\ref{eqnpsi3})
with the
same constants uniformly in~$N$.
\end{assumptions}

It is straightforward to show that the above assumptions on
$\Psi^N$ imply that the sequence of measures $\{\pi^N\}$ converges to
$\pi$
in the Hellinger metric (see~\cite{cottetal09}). Therefore, the
measures $\{\pi^N\}$ are good candidates for
\textit{finite-dimensional} approximations of~$\pi$.
Furthermore, the normalizing constants~$M_{\Psi^N}$ are uniformly
bounded and we use this fact to obtain uniform bounds
on moments of functionals in $\h$ under $\pi^N$.
\begin{lemma} \label{lemmchmeas} Under the Assumptions~\ref{ass2}
on $\Psi^N$,
\[
\sup_{N \in\mathbb{N}} M_{\Psi^N} < \infty
\]
and for any measurable functional $f\dvtx\h\mapsto\mathbb{R}$, and any
$p \geq1$,
%
%
\begin{equation}\label{eqnchmeasrn}
\sup_{N \in\mathbb{N}} \mathbb{E}^{\pi^N}|f(x)|^p \leq
M\mathbb{E}^{\pi_0}|f(x)|^p .
\end{equation}
\end{lemma}
\begin{pf}By definition,
\begin{eqnarray*}
M_{\Psi^N}^{-1} &=& \int_{\h} \exp\{- \Psi^N(x)\} \pi_0(dx)
\geq\int_{\h} \exp\{- M(1 + \|x\|_s^2)\} \pi_0(dx) \\
&\geq& e^{-2M}
\mathbb{P}^{\pi_0}(\|x\|_s \leq1)
\end{eqnarray*}
and therefore, if $\inf\{M_{\Psi^N}^{-1} \dvtx N \in\mathbb{N}\} > 0$,
then $\sup\{ M_{\Psi^N} \dvtx N \in\mathbb{N}\} < \infty$. Hence, for any
$f\dvtx\h\mapsto\mathbb{R}$,
\[
\sup_{N \in\mathbb{N}} \mathbb{E}^{\pi^N}|f(x)|^p \leq\sup_{N
\in
\mathbb{N}} M_{\Psi^N} \mathbb{E}^{\pi_0}\bigl(e^{-\Psi^N(x)}|f(x)|^p \bigr)
\leq M\mathbb{E}^{\pi_0}|f(x)|^p
\]
proving the lemma.
\end{pf}

The uniform estimate given in~(\ref{eqnchmeasrn}) will be used
repeatedly in
the sequel.


\subsection{Statement and proof of the main theorem}\label{secspde}
\label{secconv}
The assumptions made above allow us to fully state
the main result of this article, as outlined in
Section~\ref{secconvSPDEOutline}.
\begin{theorem}\label{thmmain} Let the\vspace*{-1pt} Assumptions
\ref{ass1},~\ref{ass2} hold.
Let the initial condition $x^0$ of the
RWM algorithm be such that $x^0 \dist\pi^N$ and let
$z^N(t)$ be a~piecewise linear, continuous interpolant of the RWM algorithm
(\ref{eqnRWMfancy})
as defined in~(\ref{eqnMCMCe}). Then
$z^N(t)$ converges weakly in $C([0,T],\h^s)$ to the diffusion process
$z(t)$ given by~(\ref{eqnspde}) with $z(0) \dist\pi$.
\end{theorem}

Throughout the remainder of the paper we assume
that Assumptions~\ref{ass1}, \ref{ass2}
hold, without explicitly stating this fact.
The proof of Theorem~\ref{thmmain} is given
below and relies on Proposition~\ref{thmdrift-diffus}
stated above and proved in Section~\ref{secstec},
Proposition~\ref{lembweakconv} stated above
and proved in Section~\ref{secinvp} and Lemma
\ref{lemcontmap}
which we now state and then prove at the end of this section.
\begin{lemma} \label{lemcontmap} Fix any $T>0$, any
$z^0 \in\mathcal{H}^s$ and any $W \in
C([0,T],\h^s)$.
Then the integral equation~(\ref{eqnpcxN1jcm}) has a~unique
solution $z \in C([0,T],\h^s)$. Furthermore, $z=\Theta(z^0,W)$
where
$\Theta\colon\mathcal{H}^s \times C([0,T];\mathcal{H}^s)
\rightarrow C([0,T];\mathcal{H}^s)$
as defined in~(\ref{eqnpcxN1jcm}) is continuous.
\end{lemma}
\begin{pf*}{Proof of Theorem~\ref{thmmain}}
We begin by tracking the error in
the Euler approximation argument. As before,\vadjust{\goodbreak} let $x^0 \dist\pi^N$ and
assume $x(0) = x^0$. Returning to~(\ref{eqnmexp11}), using the
definitions from~(\ref{eqnmeandrandgamm}) and
Proposition~\ref{thmdrift-diffus}, produces
%
%
\begin{eqnarray}
\label{eqnmexp1}
x^{k+1} &=& x^k + \mathbb{E}_k(x^{k+1} - x^k) +
{\sqrt{\frac{2\ell^2\beta}{N}}} \Gamma^{k+1,N},\\
\label{eqnmexp1d}
x^{k+1} &=& x^k + \frac{1}{N} d^N(x^k) +
{\sqrt{\frac{2\ell^2\beta}{N}}} \Gamma^{k+1,N}\\
\label{eqnEscheme}
&=& x^k - \frac{\ell^2 \beta}{N} m^N(x^k) +
{\sqrt{\frac{2\ell^2\beta}{N}}} \Gamma^{k+1,N}+ \frac
{r^{k+1,N}}{N},
\end{eqnarray}
where $d^N(\cdot)$ is defined as in~(\ref{eqnDN}) and $r^{k+1,N}$ as
in~(\ref{eqnRN}).
By construction, $ \mathbb{E}_k(\Gamma^{k+1,N} ) = 0$ and
%
%
\begin{eqnarray}\label{eqncovoper}
&&
\mathbb{E}_k (\Gamma^{k+1,N} \otimes\Gamma^{k+1,N} ) \nonumber\\
&&\qquad=
\frac{N}{2\ell^2\beta}\bigl[ \mathbb{E}_k \bigl((x^{k+1} -
x^k)\otimes(x^{k+1} - x^k)\bigr) \nonumber\\[-8pt]\\[-8pt]
&&\qquad\quad\hspace*{25.7pt}{} - \mathbb{E}_k(x^{k+1} - x^k)\otimes
\mathbb{E}_k(x^{k+1} - x^k)\bigr]\nonumber\\
&&\qquad= C^N + \frac{1}{2\ell^2\beta} E^{k+1,N} - \frac{N}{2\ell^2\beta}
[ \mathbb{E}_k(x^{k+1} - x^k)\otimes\mathbb{E}_k(x^{k+1} -
x^k)],
\nonumber
\end{eqnarray}
where $E^{k+1,N}$ is as given in~(\ref{eqnEN}).

Recall $t^k$ given by~(\ref{eqndeltaconst}) and $W^N$,
the linear interpolant of a~correctly scaled sum of the
$\Gamma^{k,N}$, given by~(\ref{eqnwnCproc}).
We now define $\widehat W^N$
so that~(\ref{eqnprocessxN2}) holds as stated
and hence, $\Theta(\widehat W^N)=z^N$. Define
\begin{eqnarray*}
r^N_1(t) &\eqdef& r^{k+1,N} \qquad\mbox{for } t \in[t^k,
t^{k+1}),\\
r^N_2(s) &\eqdef& \ell^2\beta\bigl(z^N(s) + C\nabla\Psi(z^N(s)) -
m^N(\bar{z}^N(s)) \bigr),
\end{eqnarray*}
where $r^{k+1,N}(\cdot)$ is given by~(\ref{eqnRN}), $m^N$ is
from~(\ref{eqnmeandrandgamm}), $\bar{z}^N$
from~(\ref{ebarZdef}) and~$z^N$ from~(\ref{eqnprocessxN1}). If
\[
\widehat W^N(t)\eqdef W^N(t)+\bigl(1/\sqrt{2\ell^2 \beta}\bigr)e^N(t)
\]
with
$ e^N(t)=\int_0^t (r^N_1(u) + r^N_2(u)) \,du $,
then~(\ref{eqnprocessxN2}) holds. To see this,
observe from~(\ref{eqnprocessxN1}) that
\begin{eqnarray*}
z^N(t) &=& z^0 + \int_0^t d^N(\bar{z}^N(u)) \,du + \sqrt{2 \ell
^2\beta} W^N(t)\\ 
&=& z^0 - \ell^2\beta\int_0^t m^N(\bar{z}^N(u)) \,du + \int_0^t
r_1^N(s) \,ds
+ \sqrt{2 \ell^2\beta} W^N(t)\\
&=& z^0 - \ell^2\beta\int_0^t \bigl(z^N(u) + C\nabla\Psi
(z^N(u))\bigr)\,du +
\int_{0}^t \bigl(r_1^N(s) + r_2^N(s)\bigr) \,ds
\\
&&{} + \sqrt{2 \ell^2\beta} W^N(t)\\
&=& z^0 - \ell^2\beta\int_0^t \bigl(z^N(u) + C\nabla\Psi
(z^N(u))\bigr)\,du +
\sqrt{2 \ell^2\beta} \widehat W^N(t)
\end{eqnarray*}
and hence, with this definition of $\widehat W^N$,
(\ref{eqnprocessxN2}) holds.

Furthermore, we claim that
%
%
\begin{equation} \label{eqnensup0}
\lim_{N \rightarrow\infty} \mathbb{E}^{\pi^N}\Bigl(\sup_{t \in
[0,T]} \|e^N(t)\|_s^2\Bigr) = 0.
\end{equation}
To prove this, notice that
\[
\sup_{t \in[0,T]} \|e^N(t)\|_s^2 \leq M\biggl(\sup_{t \in[0,T]}
\int_{0}^t \| r_1^N(u)\|_s^2 \,du+ \sup_{t \in[0,T]} \int_{0}^t
\|r_2^N(u)\|_s^2\,du\biggr).
\]
Also
\begin{eqnarray*}
&&\mathbb{E}^{\pi^N}\sup_{t \in[0,T]}\int_{0}^t \|r_1^N(u)\|_s^2 \,du
\\
&&\qquad\leq\mathbb{E}^{\pi^N}\int_{0}^T \|r_1^N(u)\|_s^2 \,du
\leq M\frac{1}{N}\mathbb{E}^{\pi^N}\sum_{k=1}^N \|r^{k,N}\|_s^2 \\
&&\qquad= M
\mathbb{E}^{\pi^N} \|r^{1,N}\|_s^2
\stackrel{N\rightarrow\infty}{\longrightarrow} 0,
\end{eqnarray*}
where we used stationarity of $r^{k,N}$ and~(\ref{thmDD-drift}) from
Proposition~\ref{thmdrift-diffus} in the last step. We now estimate the
second term similarly to complete the proof. Recall
that the function $z \mapsto z+C\nabla\Psi(z)$ is Lipschitz
on $\h^s$ by Lemma~\ref{lneeded}. Note also that
$C^N\nabla\Psi^N(\cdot)=CP^N\nabla\Psi(P^N\cdot)$.
Thus,
\begin{eqnarray*}
\|r_2^N(u)\|_s &\le& M\|z^N(u)-P^N\bar{z}^N(u)\|_s
+\|C(I-P^N)\nabla\Psi(P^N\bar{z}^N(u))\|_s\\
&\leq& M\bigl(\|z^N(u)-\bar{z}^N(u)\|_s +\|(I-P^N)\bar{z}^N(u)\|_s\bigr)
\\
&&{}+\|(I-P^N)C\nabla\Psi(P^N\bar{z}^N(u))\|_s .
\end{eqnarray*}
But for any $u \in[t^k,t^{k+1})$, we
have
\[
\|z^N(u)-\bar{z}^N(u)\|_s \le\|x^{k+1}-x^{k}\|_s \leq\|
y^{k+1}-x^{k}\|_s .
\]
This follows\vspace*{1pt} from the fact that
$\bar{z}^N(u)=x^k$ and $ z^N(u)=\frac{1}{\Delta
t}((u-t^k)x^{k+1}+(t^{k+1}-u)x^{k})$,
because $x^{k+1} - x^k = \gamma^{k+1}(y^{k+1} - x^k)$ and
$|\gamma^{k+1}|\le1$.
For $u \in[t^k,t^{k+1})$, we also have
\[
\|(P^N-I)\bar{z}^N(u)\|_s=\|(P^N-I)x^k\|_s=\|(P^N-I)x^0\|_s,
\]
because $x^k$ is not updated in $\mathcal{H}^s\setminus X^N$,
and
\[
\|(P^N-I)C\nabla\Psi(P^N\bar{z}^N(u))\|_s=
\|(P^N-I)C\nabla\Psi(P^N x^k)\|_s.
\]
Hence, we have by stationarity that, for all $u \in[0,T]$,
\begin{eqnarray*}
\mathbb{E}^\pi\|r^N_2(u)\|^2_s &\leq& M \mathbb{E}^\pi\|y^1 - x^0\|^2_s
\\
&&{}+ M\mathbb{E}^\pi\bigl(\|(P^N - I)x^0\|^2_s
+\|(P^N-I)C\nabla\Psi(P^N x^0)\|^2_s\bigr) .
\end{eqnarray*}
Equation~(\ref{eqnRWMfancy}) shows that $\mathbb{E}^\pi\|y^1 - x^0\|
^2_s \leq M N^{-1}$. The definition
of $P^N$ gives $\mathbb{E}^\pi\|(P^N - I)x\|^2_s
\leq N^{-(r-s)}\mathbb{E}^\pi\|x\|^2_r$ for any
$r \in(s, \kappa-1/2)$. Note that
$\mathbb{E}^\pi\|x^0\|^2_r$ is finite for $r \in(s, \kappa-1/2)$ by
Lemma~\ref{lemmchmeas} and
the properties of $\pi_0$. Similarly,
we have that for $r \le2\kappa-s<\kappa+\frac12$,
\begin{eqnarray*}
\EE\|C\nabla\Psi(P^N x^0)\|_r^2 &\le&
M\EE\bigl\|C^{1-(r+s)/2\kappa}\bigr\|_{\mathcal{L}(\h,\h)}\|\nabla\Psi
(P^Nx^0)\|_{-s}^2\\
&\le& M\EE(1+\|x^0\|_{s}^2).
\end{eqnarray*}
Hence, we deduce that $\mathbb{E}^{\pi^N}\|r_2^N(u)\|^2_s \to0$
uniformly for $u \in[0,T]$. It follows that
\[
\mathbb{E}^{\pi^N}\sup_{t \in[0,T]}\int_{0}^t \|r_2^N(u)\|_s^2 \,du
\leq\mathbb{E}^{\pi^N}\int_{0}^T \|r_2^N(u)\|_s^2 \,du
\leq\int_0^T \mathbb{E}^{\pi^N}\|r_2^N(u)\|_s^2 \,du \rightarrow0
\]
and we have proved the claim concerning $e^N$ made in~(\ref{eqnensup0}).

The proof\vspace*{1pt} concludes with a~straightforward
application of the
continuous mapping theorem.
Let\vspace*{-1pt} ${\widehat W}^N=W^N+ \frac{1}{\sqrt{2\ell^2 \beta}} e^N$. Let
$\Omega$ denote the probability
space generating the Markov chain in stationarity. We
have shown that
$e^N \to0$ in $L^2(\Omega;C([0,T],\h^s))$ and
by Proposition~\ref{lembweakconv}, $W^N$ converges
weakly to $W$ a~Brownian motion with covariance operator $C_s$ in
$C([0,T],\h^s)$. Furthermore,\vspace*{1pt} we also have that $W$ is independent of $z^0$.
Thus $(z^0,{\widehat W}^N)$ converges weakly to
$(z^0,W)$ in $\h^s \times C([0,T],\h^s)$, with $z^0$
and $W$ independent.
Notice that $z^N = \Theta(z^0,{\widehat W}^N)$, where $\Theta$ is
defined as in
Lemma~\ref{lemcontmap}.
Since $\Theta$ is a~continuous map by Lemma~\ref{lemcontmap},
we deduce from the continuous mapping theorem
that the process $z^N$ converges weakly in $C([0,T],\h^s)$ to $z$ with
law given by $\Theta(z^0,W)$. Since $W$ is independent of $z^0$, this
is precisely the law of the SPDE given by~(\ref{eqnspde}).
\end{pf*}
\begin{pf*}{Proof of Lemma~\ref{lemcontmap}}
Consider the mapping $z^{(n)}\mapsto z^{(n+1)}$ defined
by
\[
z^{(n+1)}(t) = z^0 - h(\ell) \int_0^t \bigl(z^{(n)}(s) + C
\nabla\Psi\bigl(z^{(n)}(s)\bigr) \bigr) \,ds +
\sqrt{2h(\ell)} W(t)
\]
for arbitrary $z^0 \in\mathcal{H}$ and $W \in C([0,T];\mathcal{H}^s)$.
Recall from Lemma~\ref{lneeded} that
$z \mapsto z+C\nabla\Psi(z)$
is globally Lipschitz on $\h^s$.
It is then a~straightforward
application of the contraction mapping
theorem to show that this mapping has a~unique
fixed point in $C([0,T]; \h^s)$, for $T$ sufficiently
small.
Repeated application of the same idea extends
this existence and uniqueness result to
arbitrary time-intervals. Let $z_i$
solve~(\ref{eqnpcxN1jcm}) with $(z^0,W)=(w_i,W_i), i=1,2$.
Subtracting the two equations and using the fact that
$z \mapsto z+C\nabla\Psi(z)$
is globally Lipschitz on $\h^s$ gives
\begin{eqnarray*}
\|z_1(t)-z_2(t)\|_s &\le& \|w_1-w_2\|_s+M\int_0^t \|z_1(s)-z_2(s)\|_s
\,ds\\
&&{} + \sqrt{2 \ell^2\beta} \|W_1(t)-W_2(t)\|_s .
\end{eqnarray*}
Thus,
\begin{eqnarray*}
\sup_{0 \le t \le T}\|z_1(t)-z_2(t)\|_s &\le&\|w_1-w_2\|_s+
M\int_0^T \sup_{0 \le\tau\le s}\|z_1(\tau)-z_2(\tau)\|_s \,ds\\
&&{} +
\sqrt{2 \ell^2\beta}\sup_{0 \le t \le T}\|W_1(t)-W_2(t)\|_s .
\end{eqnarray*}
The Gronwall lemma gives continuity in the desired spaces.
\end{pf*}

\section{\texorpdfstring{Weak convergence of the noise process: Proof of Proposition \protect\ref{lembweakconv}}
{Weak convergence of the noise process: Proof of Proposition 2.2}}
\label{secinvp}
\mbox{}

\noindent Throughout, we make the standing Assumptions~\ref{ass1},
\ref{ass2} without explicit mention.
The proof of Proposition~\ref{lembweakconv} uses the following result
concerning triangular martingale increment arrays. The result is
similar to the
classical results on triangular arrays of independent increments.

Let\vspace*{2pt} $k_N\dvtx[0,T] \rightarrow\mathbb{Z}_+$ be a~sequence of
nondecreasing, right-continuous
functions indexed by~$N$ with
$k_N(0) = 0$ and $k_N(T) \geq1$. Let\vspace*{1pt} $\{M^{k,N},\break
\mathcal{F}^{k,N}\}_{0 \leq k \leq k_N(T)}$ be an $ \h^s$ valued
martingale difference array. That is,\vspace*{1pt} for $k=1,\ldots,k_N(T)$,
we have
$\mathbb{E}(M^{k,N}|\mathcal{F}^{k-1,N}) = 0$,
$\mathbb{E}(\|M^{k,N}\|_{s}^2|\mathcal{F}^{k-1,N}) < \infty$ almost
surely, and $\mathcal{F}^{k-1,N} \subset\mathcal{F}^{k,N}$. We will
make use of the following result.

%
\begin{prop}[(\cite{Berg86}, Proposition 5.1)] \label{thmBergmclt} Let
$S\dvtx\h^s \rightarrow\h^s$ be a~self-adjoint, positive definite,
operator with finite trace. Assume that,
for all $x \in\h^s, \epsilon>0$ and
$t \in[0,T]$, the following limits hold in probability:
%
%
\begin{eqnarray}\label{eqntracmart}
\lim_{N \rightarrow\infty} \sum_{k=1}^{k_N(T)} \mathbb{E}(\|
M^{k,N}\|_s^2| \mathcal{F}^{k-1,N})
&=& T \tr(S),
\\
%
%
\label{eqntraceont}
\lim_{N \rightarrow\infty} \sum_{k=1}^{k_N(t)} \mathbb{E}(\langle
M^{k,N},x \rangle_s ^2| \mathcal{F}^{k-1,N}) &=& t\langle Sx,x
\rangle_s,
\\
%
%
\label{eqnLindtyp}
\lim_{N \rightarrow\infty} \sum_{k=1}^{k_N(T)} \mathbb{E}\bigl(\langle
M^{k,N},x \rangle^2_s 1_ {| \langle M^{k,N},x \rangle_s| \geq
\epsilon}|\mathcal{F}^{k-1,N}\bigr) &=& 0.
\end{eqnarray}
Define a~continuous time process $W^N$ by
$W^N(t) = \sum_{k=1}^{k_N(t)}M^{k,N}$ if $k_N(t) \ge1$
and $k_N(t) > \lim_{r \rightarrow0_+} k_N(t -r)$, and by
linear interpolation otherwise.
Then the sequence of\vadjust{\goodbreak} random variables $W^N$
converges weakly in $ C([0,T],\h^s)$ to an $\h^s$ valued
Brownian motion $W$, with $W(0) = 0$, $\mathbb{E}(W(T)) = 0$,
and with covariance operator $S$.
\end{prop}
\begin{remark}
\label{remcwrem}
The first two hypotheses of the above theorem
ensure the weak convergence of finite-dimensional
distributions of $W^N(t)$ using
the martingale central limit theorem in $\mathbb{R}^N$;
the last hypothesis
is needed to verify the tightness of the family $\{W^N(\ccdot)\}$.
As noted in~\cite{ChenWhit98},
the second hypothesis [equation~(\ref{eqntraceont})]
of Proposition~\ref{thmBergmclt} is implied by
%
%
\begin{equation}\label{eqntraceontmod}
\lim_{N \rightarrow\infty} \sum_{k=1}^{k_N(t)} \mathbb{E}(\langle
M^{k,N}, e_n
\rangle_s\langle M^{k,N}, e_m \rangle_s | \mathcal{F}^{k-1,N})
= t \langle Se_n, e_m \rangle_s
\end{equation}
in probability, where $\{e_n\}$ is any orthonormal basis for $\h^s$.
The third
hypothesis in~(\ref{eqnLindtyp}) is implied by
the Lindeberg type condition,
%
%
\begin{equation}\label{eqnLindtypmod}
\lim_{N \rightarrow\infty} \sum_{k=1}^{k_N(T)} \mathbb{E}\bigl(\|
M^{k,N}\| ^2 _s 1_
{ \| M^{k,N} \|_s \geq\epsilon}|\mathcal{F}^{k-1,N}\bigr)
= 0
\end{equation}
in probability, for any fixed $\epsilon>0$.
\end{remark}

Using Proposition~\ref{thmBergmclt}
we now give the proof of Proposition~\ref{lembweakconv}.
%
%
\begin{pf*}{Proof of Proposition~\ref{lembweakconv}}
We apply Proposition~\ref{thmBergmclt} with $k_N(t) \eqdef\lfloor Nt
\rfloor$, $M^{k,N} \eqdef\frac{1}{\sqrt{N}} \Gamma^{k,N}$ and $S
\eqdef C_s$;
the resulting definition\vspace*{-1pt} of $W^N(t)$ from
Proposition~\ref{thmBergmclt} coincides with that
given in~(\ref{eqnwnCproc}).
We set\vspace*{1pt} $\mathcal{F}^{k,N}$ to be the sigma algebra generated by $\{
x^j,\xi^j\}_{j \leq k}$
with $x^0 \sim\pi^N$. Since
the chain is stationary, the noise process $\{\Gamma^{k,N}, 1\leq k
\leq N\}$ is identically distributed, and so are the errors $r^{k,N}$ and
$E^{k,N}$ from~(\ref{eqnRN}) and~(\ref{eqnEN}), respectively. We
now verify
the three hypotheses required to apply Proposition~\ref{thmBergmclt}.
We generalize the notation $\mathbb{E}^\xi_0(\cdot)$
from Section~\ref{secdiftHuristic} and set
$\mathbb{E}^\xi(\cdot| \mathcal{F}^{k,N})=\mathbb{E}^\xi_k(\cdot)$.
\begin{itemize}
\item Condition~(\ref{eqntracmart}). It is enough to show that
\[
\lim_{N\rightarrow\infty} \mathbb{E}^{\pi^N}\Biggl| \frac
{1}{N}\sum_{k=1}^{\lfloor NT
\rfloor}\mathbb{E}^\xi_{k-1}(\|\Gamma^{k,N}\|_s^2 ) - \tr(C_s)
\Biggr| = 0
\]
and condition~(\ref{eqntracmart}) will follow from Markov's inequality.
By~(\ref{eqncovoper}) and~(\ref{eqnSob2}),
%
%
\begin{eqnarray}
\mathbb{E}_0^{\xi}(\|\Gamma^{1,N}\|_s^2) &=&
\sum_{j=1}^N\mathbb{E}^{\xi}_0 (\|\B^{1/2} \Gamma^{1,N}\|
^2) =
\sum_{j=1}^N \mathbb{E}^{\xi}_0 \langle\Gamma^{1,N}, \B^{
1/2} \phi_j \rangle^2 \nonumber\\
\label{eqnmcveric1}
&=&
\sum_{j=1}^N \mathbb{E}^{\xi}_0 \langle\B^{1/2} \phi_j,
\Gamma^{1,N}
\otimes\Gamma^{1,N} \B^{1/2} \phi_j \rangle\\
\label{eqnGammaNexp}
&=& \tr(C^N_s) + \frac{1}{2\ell^2 \beta} \sum_{j=1}^N \langle\phi
_j,E^{1,N} \phi_j \rangle_s\nonumber\\[-8pt]\\[-8pt]
&&{}- {N \over2\ell^2 \beta} \|\mathbb{E}_0(x^{1} - x^0) \|^2_s
.\nonumber
\end{eqnarray}
By Proposition~\ref{thmdrift-diffus} it follows that $\mathbb
{E}^{\pi^N}| {\sum_{j=1}^N} \langle\phi_j,E^{1,N} \phi_j \rangle_s
| \rightarrow0$. For the third
term, notice that
by Proposition~\ref{thmdrift-diffus}~(\ref{thmDD-drift}) we have
%
%
\begin{eqnarray}\label{eqnMartsecverref}
\mathbb{E}^{\pi^N}\frac{N}{2\ell^2\beta}\|\mathbb{E}_0(x^{1} -
x^0)\|^2_s &\leq&
M \frac{1}{N} \mathbb{E}^{\pi^N} \bigl(\|m^N(x^0)\|^2_s + \|r^{1,N}\|
^2_s\bigr)\nonumber\\
&\leq& M \frac{1}{N} \bigl(\mathbb{E}^{\pi^N} (1+ \|x^0\|_s)^2 +
\mathbb{E}^{\pi^N}\|r^{1,N}\|^2_s\bigr) \\
&\rightarrow&0,\nonumber
\end{eqnarray}
where the second inequality follows from the fact that $C\nabla\Psi$
is globally Lipschitz in $\h^s$.
Also $\{E^{k,N}\}$ is a~stationary sequence. Therefore,
\begin{eqnarray*}
&&\mathbb{E}^{\pi^N}\Biggl| \frac{1}{N}\sum_{k=1}^{\lfloor NT
\rfloor}\mathbb{E}^\xi_{k-1}(\|\Gamma^{k,N}\|^2_s ) - T\tr(C^N_s)\Biggr|\\
&&\qquad\leq
M\mathbb{E}^{\pi^N} \Biggl(\Biggl| \sum_{j=1}^N \langle\phi
_j,E^{1,N} \phi_j \rangle_s \Biggr| + \frac{N}{2\ell^2\beta}\|\mathbb
{E}_0(x^{1} -
x^0)\|^2_s \Biggr) \\
&&\qquad\quad{} + \tr(C^N_s)\biggl|\frac{ \lfloor NT \rfloor}{N} - T\biggr|
\rightarrow0 .
\end{eqnarray*}
Condition~(\ref{eqntracmart}) now follows from the fact that
\[
\lim_{N\rightarrow\infty} |{\tr}(C_s) - \tr(C^N_s)| = 0.
\]

%
\item
Condition~(\ref{eqntraceont}). By Remark~\ref{remcwrem}, it is enough
to verify~(\ref{eqntraceontmod}).
To show
(\ref{eqntraceontmod}), using stationarity and similar arguments used
in verifying condition~(\ref{eqntracmart}), it suffices to show that
%
%
\begin{equation}
\lim_{N \rightarrow\infty} \mathbb{E}^{\pi^N}|\mathbb
{E}^{\xi}_0 (\langle\Gamma^{1,N},
\widehat{\phi}_n \rangle_s\langle\Gamma^{1,N}, \widehat{\phi}_m
\rangle_s) - \langle\widehat{\phi}_n , C^N_s
\widehat{\phi}_m \rangle_s
| = 0, 
\end{equation}
where $\{\widehat{\phi}_k\}$ is as defined in~(\ref{eqnphat}). We have
\begin{eqnarray*}
&&
\mathbb{E}^{\pi^N}|\mathbb{E}^{\xi}_0 (\langle\Gamma^{1,N},
\widehat{\phi}_n \rangle_s\langle\Gamma^{1,N}, \widehat{\phi}_m
\rangle_s) - \langle\widehat{\phi}_n , C^N_s
\widehat{\phi}_m \rangle_s| \\
&&\qquad=
n^{-s} m^{-s}\mathbb{E}^{\pi^N}|\mathbb{E}^{\xi
}_0 (\langle\Gamma^{1,N},
{\phi}_n \rangle_s\langle\Gamma^{1,N}, {\phi}_m \rangle_s) -
\langle{\phi}_n , C^N_s
{\phi}_m \rangle_s
|
\end{eqnarray*}
and therefore, it is enough to show that
%
%
\begin{equation}\label{eqn2mcltverfcond}
\lim_{N \rightarrow\infty}\mathbb{E}^{\pi^N}|\mathbb{E}^{\xi
}_0 (\langle\Gamma^{1,N},
{\phi}_n \rangle_s\langle\Gamma^{1,N}, {\phi}_m \rangle_s) -
\langle{\phi}_n , C^N_s
{\phi}_m \rangle_s
| = 0 .
\end{equation}
Indeed we have
\begin{eqnarray*}
\langle\Gamma^{1,N}, \phi_n \rangle_s \langle\Gamma^{1,N},
\phi_m \rangle_s &=&\langle\Gamma^{1,N}, \B\phi_n \rangle
\langle\Gamma^{1,N},
\B\phi_m \rangle\\
&=&\langle\B\phi_n , \Gamma^{1,N}\otimes
\Gamma^{1,N}
\B\phi_m \rangle\\
&=&\langle\phi_n , \B^{1/2} \Gamma^{1,N}\otimes\Gamma
^{1,N}
\B^{1/2} \phi_m \rangle_s
\end{eqnarray*}
and from~(\ref{eqncovoper}) and
Proposition~\ref{thmdrift-diffus} we obtain
\begin{eqnarray*} 
&&\langle\phi_n , \B^{1/2} \Gamma^{1,N}\otimes\Gamma^{1,N}
\B^{1/2} \phi_m \rangle_s -
\langle\phi_n , C^N_s \phi_m \rangle_s \\
&&\qquad
= \langle\phi_n , \B^{1/2} \Gamma^{1,N}\otimes\Gamma
^{1,N}
\B^{1/2} \phi_m \rangle_s -
\langle\phi_n , \B^{1/2} C^N \B^{1/2} \phi_m
\rangle_s\\
&&\qquad= n^s m^s \langle\phi_n, E^{1,N} \phi_m \rangle_s
- {N \over2\ell^2 \beta}
\mathbb{E}_0(\langle x^1 - x^0, \phi_n \rangle_s ) \mathbb
{E}_0(\langle x^1 - x^0, \phi_m \rangle_s) .
\end{eqnarray*}
From Proposition~\ref{thmdrift-diffus}, it follows that $\lim_{N
\rightarrow
\infty} \mathbb{E}^{\pi^N}|\langle\phi_n, E^{1,N} \phi_m \rangle
_s| = 0$. Also notice that
\begin{eqnarray*}
&&N^2 [\mathbb{E}^{\pi^N}| \mathbb{E}_0(\langle x^1 - x^0, \phi
_n \rangle_s )
\mathbb{E}_0(\langle x^1 - x^0, \phi_m \rangle_s)| ]^2\\
&&\qquad\leq M
\mathbb{E}^{\pi^N}\bigl(N\|\mathbb{E}_0(x^1 - x^0)
\|^2_s \|\phi_n\|^2_s \bigr) \mathbb{E}^{\pi^N} \bigl( N\|\mathbb
{E}_0(x^1 - x^0) \|^2_s \|\phi_m\|^2_s\bigr) \\
&&\qquad\rightarrow0
\end{eqnarray*}
by the calculation done in
(\ref{eqnMartsecverref}). Thus~(\ref{eqn2mcltverfcond})
holds and since $|\langle\phi_n , C_s \phi_m \rangle_s- \langle
\phi_n , C^N_s \phi_m \rangle_s| \rightarrow0$, equation (\ref
{eqntraceont}) follows from
Markov's inequality.

\item
Condition~(\ref{eqnLindtyp}). From Remark~\ref{remcwrem} it follows that
verifying~(\ref{eqnLindtypmod}) suffices
to establish~(\ref{eqnLindtyp}).

To verify~(\ref{eqnLindtypmod}), notice that
for any $\epsilon>0$,
\begin{eqnarray*}
&&\mathbb{E}^{\pi^N}\Biggl|\frac{1}{N}\sum_{k=1}^{ \lfloor NT \rfloor
} \mathbb{E}^{\xi}_{k-1}\bigl(\|\Gamma^{k,N}\|^2_s 1_{
\{\|\Gamma^{k,N}\|^2_s \geq\epsilon N \}}\bigr)\Biggr|\\
&&\qquad \leq
\frac{ \lfloor NT \rfloor}{N} \mathbb{E}^{\pi^N}\bigl(\|\Gamma^{1,N}\|
^2_s 1_{ \{\|\Gamma^{1,N}\|^2_s \geq
\epsilon N \}}\bigr) \rightarrow0
\end{eqnarray*}
by the dominated convergence theorem since
\[
\lim_{N \rightarrow
\infty} \mathbb{E}^{\pi^N} \|\Gamma^{1,N} \|^2_s = \tr(C_s) <
\infty.
\]
Thus~(\ref{eqnLindtypmod}) is verified.
\end{itemize}
Thus we have verified all three hypotheses of
Proposition~\ref{thmBergmclt}, proving that $W^N(t)$ converges weakly to
$W(t)$ in $C([0,T];\h^s)$.

Recall that $X^R \subset\h^s$ denotes the $R$-dimensional
subspace $P^R \h^s$.
To prove the second claim of Proposition~\ref{lembweakconv}, we need
to show that $(x^0, W^N(t))$\vspace*{1pt} converges
weakly to $(z^0, W(t))$ in $(\h^s, C([0,T];\h^s) )$ as $N
\rightarrow\infty$ where $z^0 \sim\pi$
and $z^0$ is independent of the limiting
noise $W$. For showing this, it is enough to show that
for any $R \in\mathbb{N}$, the pair $(x^0, P^R W^N(t))$ converges
weakly to $(z^0, Z_R)$ for every $t > 0$, where
$Z_R$ is a~Gaussian random variable on $X^R$ with mean zero, covariance
$t P^R C_s P^R$
and independent of $z^0$. We will prove this statement as the corollary
of the following lemma.
\begin{lemma} \label{lemmartconvlem}
Let $x^0 \sim\pi^N$ and let $\{\theta^{k,N}\}$ be any stationary
martingale sequence adapted to the filtration
$\{\mathcal{F}^{k,N}\}$ and furthermore, assume that there exists a~stationary sequence
$\{U^{k,N}\}$ such that for all $k \geq1 $ and any $u \in X^R$:
\begin{longlist}[(1)]
\item[(1)] $\mathbb{E}^\xi_{k-1}| \langle u, P^R\theta^{k,N} \rangle_s|^2
= \langle u, P^R C_s u \rangle_s + U^{k,N},
\lim_{N \rightarrow\infty}\mathbb{E}^{\pi^N}|U^{1,N}|= 0$.
\item[(2)] $\mathbb{E}^\xi_{k-1}\| \theta^{k,N}\|_s^3 \leq M$.
\end{longlist}
%
Then for any $\mathbf{t} \in\h^s$, $u \in X^R$, $R \in\mathbb{N}$
and $t > 0$,
%
%
\begin{eqnarray} \label{eqnUNapp}
&&\lim_{N \rightarrow\infty} \mathbb{E}^{\pi^N}\bigl(e^{i\langle\mathbf{t},
x^0 \rangle_s + (i /{\sqrt{N}})\sum_{k=1}^{\lfloor Nt
\rfloor} \langle u, P^R \theta^{k,N} \rangle_s} \bigr)\nonumber\\[-8pt]\\[-8pt]
&&\qquad =
\mathbb{E}^{\pi}\bigl(e^{i\langle\mathbf{t}, z^0 \rangle_s - (t/
{2}) \langle u, P^R C_s u \rangle_s} \bigr).\nonumber
\end{eqnarray}
\end{lemma}

Note: Here and in Corollary~\ref{corx0wind}, $i = \sqrt{-1}$.
\begin{pf*}{Proof of Lemma~\ref{lemmartconvlem}}
We show~(\ref{eqnUNapp}) for $t=1$, since the calculations are nearly
identical for an arbitrary $t$
with minor notational changes.
Indeed, we have
\begin{eqnarray*}
&&\mathbb{E}^{\pi^N}\bigl(e^{i\langle\mathbf{t}, x^0 \rangle_s + (i/
{\sqrt{N}})\sum_{k=1}^{ N } \langle u, P^R \theta^{k,N}
\rangle_s} \bigr) \\
&&\qquad= \mathbb{E}^{\pi^N}
\bigl(\mathbb{E}^\xi_{N-1}\bigl(e^{i\langle\mathbf{t}, x^0 \rangle_s + (i
/{\sqrt{N}})\sum_{k=1}^{ N } \langle u, P^R \theta^{k,N}
\rangle_s} \bigr)\bigr).
\end{eqnarray*}
By Taylor's expansion,
%
%
\begin{eqnarray}
\label{eqnUNVN}\qquad
&&\mathbb{E}^{\pi^N}
\bigl(\mathbb{E}^\xi_{N-1}\bigl(e^{i\langle\mathbf{t}, x^0 \rangle_s + (i
{\sqrt{N}})\sum_{k=1}^{ N } \langle u, P^R \theta^{k,N}
\rangle_s} \bigr)\bigr)\nonumber\\
&&\qquad=
\mathbb{E} \biggl[e^{i\langle\mathbf{t}, x^0 \rangle
_s +
(i /{\sqrt{N}})\sum_{k=1}^{N-1} \langle u, P^R \theta^{k,N}
\rangle_s} \nonumber\\[-8pt]\\[-8pt]
&&\qquad\quad\hspace*{10pt}{}\times\biggl(1 - \frac{1}{2N}\mathbb{E}^\xi_{N-1}| \langle u, P^R
\theta^{N,N} \rangle_s|^2\nonumber\\
&&\qquad\quad\hspace*{69.4pt}{} + M\biggl(\frac{1}{N^{3/2}}V^N \wedge2\biggr)
\biggr)\biggr],\nonumber
\end{eqnarray}
where $|V^N| \leq\mathbb{E}^\xi_{N-1}| \langle u, P^R \theta^{N,N}
\rangle_s|^3 \leq M$, since by assumption $\mathbb{E}^\xi_{N-1}\|
\theta^{N,N}\|_s^3 \leq M$. We also have that
\begin{eqnarray*}
\mathbb{E}^\xi_{N-1}| \langle u, P^R\theta^{N,N} \rangle_s|^2
& = &\langle u, P^R C_s u \rangle_s+ U^{N,N} ,\\
\lim_{N\rightarrow\infty}\mathbb{E}^{\pi^N}|U^{N,N}|&=& 0 .
\end{eqnarray*}
Thus from~(\ref{eqnUNVN}) we deduce that
%
%
\begin{eqnarray}\label{eqnSNuN}\qquad
&&\hspace*{1pt}\mathbb{E}^{\pi^N} \bigl(e^{i\langle\mathbf{t}, x^0 \rangle_s + (i
/{\sqrt{N}})\sum_{k=1}^N \langle u, P^R \theta^{k,N} \rangle_s}\bigr) \nonumber\\
&&\hspace*{1pt}\qquad=
\mathbb{E}^{\pi^N} \biggl[e^{i\langle\mathbf{t}, x^0 \rangle_s + (i
/{\sqrt{N}})\sum_{k=1}^{N-1} \langle u, P^R \theta^{k,N}
\rangle_s} \biggl(1 - \frac{1}{2N} \langle u, P^RC_s u \rangle_s
\biggr)\biggr] \nonumber\\
&&\hspace*{1pt}\qquad\quad{}+ S^N,\\
&&|S^N| \leq M \mathbb{E}^{\pi^N} \biggl(\frac{1}{2N} |U^{N,N}| +
\frac{1}{N^{3/2}} |V^N|\biggr) \nonumber\\
&&\hphantom{|S^N|}= M\frac{1}{N}\mathbb{E}^{\pi^N} \biggl( |U^{N,N}|
+ \frac{1}{\sqrt{N}}\biggr) .\nonumber
\end{eqnarray}
Proceeding recursively we obtain
\begin{eqnarray*}
&&\mathbb{E}^{\pi^N} \bigl(e^{i\langle\mathbf{t}, x^0 \rangle_s + (i
/{\sqrt{N}})\sum_{k=1}^N \langle u, P^R \theta^{k,N} \rangle
_s} \bigr)\\
&&\qquad =
\mathbb{E}^{\pi^N} \biggl[e^{i\langle\mathbf{t}, x^0 \rangle_s}
\biggl(1 - \frac{1}{2N} \langle u, P^R C_s u \rangle_s \biggr)^N
\biggr] + \sum_{k=1}^NS^k.
\end{eqnarray*}
By the stationarity of $\{U^{k,N}\}$ and the fact that $\mathbb{E}^\pi
|U^{k,N}| \rightarrow0$ as
$N \rightarrow\infty$, from~(\ref{eqnSNuN}) it follows that
\[
\sum_{k=1}^N |S^k| \leq M \sum_{k=1}^N\frac{1}{N}\biggl(\mathbb{E}^{\pi
^N} |U^k| + \frac{1}{\sqrt{N}}\biggr) \leq M \biggl(\mathbb{E}^{\pi^N} |U^1| +
\frac{1}{\sqrt{N}}\biggr) \rightarrow0.
\]
Thus we have shown that
\[
\mathbb{E}^{\pi^N} \biggl[e^{i\langle\mathbf{t}, x^0
\rangle_s} \biggl(1 - \frac{1}{2N} \langle u, P^RC_s u
\rangle_s \biggr)^N\biggr] = \mathbb{E}^{\pi^N} \bigl[e^{i\langle
\mathbf{t}, x^0 \rangle_s - ({1/2}) \langle u, P^RC_s u \rangle
_s} \bigr] +o(1),
\]
and
the result follows from the fact that $\mathbb{E}^{\pi^N}
[e^{i\langle\mathbf{t}, x^0 \rangle_s} ] \rightarrow\mathbb
{E}^{\pi} [e^{i\langle\mathbf{t}, z^0 \rangle_s} ]$,
finishing the proof
of Lemma~\ref{lemmartconvlem}.
\end{pf*}

As a~corollary of Lemma~\ref{lemmartconvlem}, we obtain the following.
\begin{corollary} \label{corx0wind}The pair $(x^0, W^N)$ converges
weakly to $(z^0, W)$
in $C([0, T];\h^s)$ where $W$ is a~Brownian motion with covariance
operator $C_s$
and is independent of $z^0$ almost surely.
\end{corollary}
\begin{pf}
As mentioned before, it is enough to show that for any $\mathbf{t} \in
\h^s$, $u \in X^R$, $R \in\mathbb{N}$ and $t > 0$,
%
%
\begin{eqnarray}\label{corchareq}
&&\lim_{N \rightarrow\infty} \mathbb{E}^{\pi^N}\bigl(e^{i\langle\mathbf{t},
x^0 \rangle_s + (i /{\sqrt{N}})\sum_{k=1}^{\lfloor Nt
\rfloor} \langle u, P^R \Gamma^{k,N} \rangle_s} \bigr) \nonumber\\[-8pt]\\[-8pt]
&&\qquad=
\mathbb{E}^{\pi}\bigl(e^{i\langle\mathbf{t}, z^0 \rangle_s - (t/{2})
\langle u, P^R C_s u \rangle_s} \bigr).\nonumber
\end{eqnarray}
Now we verify the conditions of Lemma~\ref{lemmartconvlem} to show
(\ref{corchareq}). To verify the first hypothesis of Lemma \ref
{lemmartconvlem}, notice that from Proposition~\ref{thmdrift-diffus}
we obtain that for $k \geq1$,
\begin{eqnarray*}
\mathbb{E}^\xi_{k-1}| \langle u, P^R \Gamma^{k,N} \rangle_s|^2
&=& \mathbb{E}^\xi_{k-1} \langle\B u, P^R \Gamma^{k,N} \otimes
\Gamma^{k,N} \B u \rangle\\
&=& \langle u, P^R C_s u \rangle_s+
U^{k,N}, \\
|U^{k,N}|&\leq&\frac{1}{2\ell^2 \beta} M \sum_{l,j=1}^{R \wedge N}
u_l u_j |\langle
\phi_l,P^ME^{k,N} \phi_j \rangle_s|\\
&&{} + {N \over2\ell^2 \beta} \|
\mathbb{E}^\xi_{k-1}(x^{{k}} - x^{k-1}) \|^2_s \| u \|_s^2 \\
&&{} + |\langle u, P^R C^N_s u \rangle_s - \langle u, P^R C_s u
\rangle_s|,
\end{eqnarray*}
where $\{E^{k,N}\}$ is as defined in~(\ref{eqnEN}). Because $\{\Gamma
^{k,N}\}$ is stationary, we deduce that $\{U^{k,N}\}$ is stationary.
From Proposition~\ref{thmdrift-diffus} we obtain
\[
\lim_{N \rightarrow\infty} \sum_{l,j =1}^{R \wedge N}\mathbb
{E}^{\pi^N}|\langle
\phi_l,P^ME^{k,N} \phi_j \rangle_s| = 0
\]
and $\mathbb{E}^{\pi^N}{N \over2\ell^2 \beta} \|\mathbb{E}^\xi
_{k-1}(x^{{k}} - x^{k-1}) \|_s^2 \rightarrow0$ by the calculation in
(\ref{eqnMartsecverref}). Thus we have shown that $\mathbb{E}^\pi
|U^{1,N}| \rightarrow0$ as
$N \rightarrow\infty$.
The second hypothesis of Lem\-ma~\ref{lemmartconvlem} is easily
verified since $\mathbb{E}^\xi_{k-1} \|\Gamma^{k,N} \|^3_s \leq M
\mathbb{E}^\xi_{k-1}\|C^{1/2} \xi^k\|^3_s \leq M$.
Thus the corollary follows from Lemma~\ref{lemmartconvlem}.
\end{pf}

Thus we have shown that $(x^0, W^N)$ converges weakly to $(z^0,W)$ where
$W$ is a~Brownian motion in $\h^s$ with covariance operator $C_s$, and
by the above corollary we
see that $W$ is independent of $x^0$ almost surely, proving the two
claims made in Proposition~\ref{lembweakconv} and the proof is complete.
\end{pf*}

\section{\texorpdfstring{Mean drift and diffusion: Proof of Proposition \protect\ref{thmdrift-diffus}}
{Mean drift and diffusion: Proof of Proposition 2.1}}\label{secstec}

To prove this key proposition we make the standing
Assumptions~\ref{ass1},~\ref{ass2}
from Section~\ref{secpsin} without
explicit statement of this fact within the individual
lemmas.
We start with several preliminary bounds and then
consider the drift and diffusion terms, respectively.

\subsection{Preliminary estimates}

Recall the definitions of $R(x,\xi)$, $R_i(x,\xi)$ and
$R_{ij}(x,\xi)$ from equations~(\ref{eqnR}),~(\ref{eqnRi}) and
(\ref{eqnRij}), respectively. These quantities were introduced so
that the term in the exponential
of the acceptance probability $Q(x,\xi)$ could be replaced with
$R_i(x,\xi)$ and $R_{ij}(x,\xi)$ to take advantage of the fact that,
conditional on $x$, $R_i(x,\xi)$ is independent of $\xi_i$ and
$R_{ij}(x,\xi)$ is independent of $\xi_i,\xi_j$. In the next lemma, we
estimate the additional error due to this replacement of $Q(x,\xi)$.
Recall that $\mathbb{E}_0^{\xi}$ denotes expectation
with respect to $\xi=\xi_0$ as in Section~\ref{ssecalg}.\vadjust{\goodbreak}
\begin{lemma}\label{lemestRx}
%
%
\begin{eqnarray}
\label{eqnRiest}
\mathbb{E}_0^{\xi} |Q(x,\xi) - R_i(x,\xi)|^2 &\leq& \frac{M}{N} (1
+ |\zeta_i|^2) , \\
\label{eqnRijest}
\mathbb{E}_0^\xi\bigl(Q(x,\xi) - R_{ij}(x,\xi)\bigr)^2 &\leq&
\frac{M}{N} (1 + |\zeta_i|^2 + |\zeta_j|^2) .
\end{eqnarray}
\end{lemma}
\begin{pf}
Since $\xi_j$ are i.i.d. $\mathrm{N}(0,1)$, using~(\ref{eqnSob}) and
(\ref{eqndecayeigen}), we obtain that
%
%
\begin{equation}\label{eqnRremest}
\mathbb{E}\|C^{1/2}\xi\|_s^4 \leq3(\mathbb{E}\|C^{1/2}\xi\|_s^2)^2
\leq M \Biggl(\sum_{j=1}^{\infty} j^{2s - 2k}\Biggr)^2 < \infty
\end{equation}
since $s<k-\frac12$.

Starting from~(\ref{eqnQRiapprox}), the estimates in~(\ref{eqnrx,y,xi})
and~(\ref{eqnRremest}) imply that
\begin{eqnarray*}
\mathbb{E}_0^{\xi} |Q(x,\xi) - R_i(x,\xi)|^2 &\leq& M
\biggl(\mathbb{E}^\xi_0|r(x,\xi)|^2 + \frac{1}{N} \mathbb{E}^\xi_0
\zeta_i^2 \xi_i^2 +
\frac{1}{N^2}\mathbb{E}\xi_i^4\biggr) \\
&\leq& M \biggl(\frac{1}{N^2}\mathbb{E}\|C^{1/2}\xi\|_s^4 + \frac{1}{N}
\zeta_i^2 + \frac{3}{N^2}\biggr) \\
&\leq& M \frac{1}{N}
(1 + \zeta_i^2 )
\end{eqnarray*}
verifying the first part of the lemma.
A very similar argument for the second part finishes the proof.
\end{pf}

The random variables $R(x,\xi)$, $R_i(x,\xi)$ and $R_{ij}(x,\xi)$
are approximately Gaussian random variables. Indeed it can be readily
seen that
\[
R(x,\xi) \approx\Normal\biggl(-\ell^2,2\frac{\ell^2}{N}\|\zeta\|^2\biggr) .
\]
The next lemma contains a~crucial observation.
We show that the sequence
of random variables $\{\frac{\|\zeta\|^2}{N}\}$ converges to $1$
almost surely under both $\pi_0$ and~$\pi$. Thus $R(x,\xi)$ converges
almost surely to $Z_\ell\eqdef\Normal(-\ell^2,2\ell^2)$ and
thus the expected acceptance probability $\mathbb{E}\alpha(x,\xi) =
1\wedge e^{Q(x,\xi)}$
converges to $\beta= \mathbb{E}(1 \wedge e^{Z_\ell})$.
\begin{lemma}\label{lemshiftmeas}As $N \rightarrow\infty$ we have
%
%
\begin{equation}\label{eqnalzeta0-1}
\frac{1}{N} \|\zeta\|^2 \rightarrow1, \qquad\pi_0 \mbox{-a.s.}
\quad\mbox{and}\quad
\frac{1}{N} \|\zeta\|^2 \rightarrow1, \qquad\pi\mbox{-a.s.}
\end{equation}
Furthermore, for any $m \in\mathbb{N}$, $\alpha\geq2$, $s < \kappa
- {1 \over2}$
and for any $ c \geq0$,
%
%
\begin{equation}\label{eqnalzeta2-3}
\limsup_{N \in\mathbb{N}}
\mathbb{E}^{\pi^N} \sum_{j=1}^N \lambda_j^{\alpha} j^{2s}
|\zeta_j|^m e^{(c/{N}) \|\zeta\|^2} < \infty.
\end{equation}
Finally, we have
%
%
\begin{equation}\label{eqnalzeta4}
\lim_{N \to\infty} \mathbb{E}^{\pi^N}\biggl(\biggl| 1 -
\frac{1}{N}\|\zeta\|^2 \biggr|^2\biggr) = 0 .
\end{equation}
\end{lemma}
%
%
\begin{pf}
The proof proceeds by
showing the conclusions first in the case when $x \dist\pi_0$; this
is easier
because the finite-dimensional distributions are Gaussian and by
Fernique's theorem
$x$ has exponential moments. Next we notice that the almost sure properties
are preserved under the change of measure~$\pi$. To show the
convergence of moments, we use our hypothesis that the
Radon--Nikodym\vspace*{1pt} derivative $\frac{d\pi^N}{d\pi_0}$ is bounded
from above
independently of~$N$, as shown in Lemma~\ref{lemmchmeas}, equation~(\ref{eqnchmeasrn}).

Indeed, first let $x \dist\pi_0$.
Recall that $\zeta=C^{-1/2} (P^Nx)+C^{1/2}\nabla\Psi^N(x)$ and
%
%
\begin{equation}
\|\nabla\Psi^N( x) \|_{-s}\le M_3(1+\|x\|_s).
\end{equation}
Using~(\ref{eqequiv}) and the fact that $s< \kappa- \frac{1}{2}$
so that $-\kappa<-s$, we deduce that
\begin{eqnarray*}
\|C^{1/2} \nabla\Psi^N(x) \| & \asymp&\|\nabla\Psi^N(x) \|_{-\kappa
}\\
& \le&\|\nabla\Psi^N(x) \|_{-s}\\
& \le &M (1+\|x\|_s)
\end{eqnarray*}
uniformly in~$N$. Also, since $x$ is Gaussian under $\pi_0$,
from~(\ref{eqnKLexp}), we may write $C^{-1/2}(P^Nx) = \sum
_{k=1}^N \rho_k \phi_k$, where
$\rho_k$
are i.i.d. $\Normal(0,1)$.
Note that
%
%
\begin{eqnarray}\label{eqnerchmain}\qquad
\frac{1}{N}\|\zeta\|^2& = & \frac{1}{N}\|C^{-1/2}(P^Nx)+C^{1/2}\nabla\Psi^N(x)\|^2\nonumber\\
&=&\frac{1}{N}\bigl(\|C^{-1/2}(P^Nx)\|^2+2\langle C^{-1/2}(P^Nx),
C^{1/2}\nabla\Psi^N(x)
\rangle\nonumber\\
&&\hspace*{151pt}{}+\|C^{1/2}\nabla\Psi^N(x)\|^2\bigr)\\
&=&\frac{1}{N}\bigl(\|C^{-1/2}(P^Nx)\|^2+2\langle P^Nx, \nabla
\Psi^N(x)
\rangle+\|C^{1/2}\nabla\Psi^N(x)\|^2\bigr)\nonumber\\
&=&\frac{1}{N} \sum_{k=1}^N \rho_k^2+\gamma,\nonumber
\end{eqnarray}
where
%
%
\begin{eqnarray}\label{eqnechormeas}
|\gamma| &\le&\frac{1}{N}\bigl(2\|x\|_s\|\nabla\Psi^N(x)\|_{-s}+
\|C^{1/2}\nabla\Psi^N(x)\|^2\bigr)\nonumber\\[-8pt]\\[-8pt]
&\leq&
\frac{M}{N}\bigl(2\|x\|_s(1 + \|x\|_s) + (1 +
\|x\|_s)^2\bigr) .\nonumber
\end{eqnarray}
Under $\pi_0$, we have $\|x\|_{s} < \infty$ a.s., for
$s<\kappa-\frac12$ and hence, by
(\ref{eqnechormeas}), we conclude that $|\gamma| \rightarrow0$
almost surely as $N\rightarrow\infty$. 
Now, by the strong law of large numbers, $\frac{1}{N} \sum_{k=1}^N
\rho_k^2 \rightarrow1$ almost\vadjust{\goodbreak} surely. Hence, from
(\ref{eqnerchmain}) we obtain that under $\pi_0$, $\lim_{N
\rightarrow\infty} \frac{1}{N} \|\zeta\|^2 = 1$ almost surely,
proving the first equation in~(\ref{eqnalzeta0-1}). Now the second
equation in~(\ref{eqnalzeta0-1}) follows by noting that almost sure
limits are preserved under a~(absolutely continuous) change of
measure.

Next, notice that by~(\ref{eqnerchmain}) and
the Cauchy--Schwarz inequality, for any $c>0$,
\begin{eqnarray*}
\bigl(\mathbb{E}^{\pi_0} e^{(c /{N}) \|\zeta\|^2}\bigr)^2 &\leq&
\bigl(\mathbb{E}^{\pi_0} e^{(2{c}/{N})\sum\rho_k^2}\bigr)
(\mathbb{E}^{\pi_0} e^{2c\gamma})\\
& \leq&\bigl(\mathbb
{E}^{\pi_0}
e^{(2{c}/{N})\sum\rho_k^2}\bigr)
\bigl(\mathbb{E}^{\pi_0}e^{({M / N})\|x\|_s^2}\bigr) .
\end{eqnarray*}
Using the fact that $\sum_{k=1}^N \rho^2_k$ has chi-squared
distribution with~$N$ degrees of freedom gives
%
%
\begin{equation}\label{eqnlimsupenest}
\bigl(\mathbb{E}^{\pi_0} e^{(c /{N}) \|\zeta\|^2}\bigr)^2
\leq Me^{-({N}/{2})\log(1
- {4c}/{N})} \bigl(\mathbb{E}^{\pi_0} e^{({M/N})
\|x\|^2_s}\bigr) \leq M,
\end{equation}
where the last inequality follows from Fernique's theorem since
$\mathbb{E}^{\pi_0} e^{({M/N})\|x\|_s^2} < \infty$ 
for sufficiently large~$N$. Hence,
by applying Lemma~\ref{lemmchmeas}, equation
(\ref{eqnchmeasrn}), it follows that $\limsup_{N \rightarrow
\infty}\mathbb{E}^{\pi^N} e^{(c /{N}) \|\zeta\|^2} < \infty$.
Notice that we also have the bound
\[
|\zeta_k|^m \leq M\bigl(|\rho_k|^m + |\lambda_k|^m(1 + \|x\|^{m}_s)\bigr) .
\]
Since $s < k-1/2$, we have that $ \sum_{j=1}^\infty\lambda^2_j
j^{2s} < \infty$ and
therefore, it follows that for $\alpha\geq2$,
%
%
\begin{equation} \label{eqnlamzetacl}
\limsup_{N \rightarrow\infty} \sum_{k=1}^N (\mathbb{E}^{\pi^N}
\lambda_k^{2\alpha} j^{2s} |\zeta_k|^{2m})^{1/2} < \infty.
\end{equation}
Hence the claim in
(\ref{eqnalzeta2-3}) follows from applying Cauchy--Schwarz combined
with~(\ref{eqnlimsupenest}) and~(\ref{eqnlamzetacl}).
Similarly,\vspace*{1pt} a~straightforward calculation yields
that $\mathbb{E}^{\pi^0}(| 1 - \frac{1}{N}\|\zeta\|^2 |^2)
\leq\frac{M}{N}$. Hence, again by Lemma~\ref{lemmchmeas},
\[
\lim_{N\rightarrow\infty} \mathbb{E}^{\pi^N}\biggl(\biggl| 1 -
\frac{1}{N}\|\zeta\|^2 \biggr|^2\biggr)=0
\]
%
proving the last claim and the proof is complete.
\end{pf}

Recall that $Q(x,\xi) = R(x,\xi) - r(x,\xi)$. Thus, from (\ref
{eqnrx,y,xi}) and Lemma~\ref{lemestRx} it follows that
$R_i(x,\xi)$
and $R_{ij}(x,\xi)$ also are approximately Gaussian.
Therefore, the conclusion of Lemma~\ref{lemshiftmeas} leads to the
reasoning that,
for any fixed realization of $x \dist\pi$, the random variables
$R(x,\xi), R_i(x,\xi)$
and $ R_{ij}(x,\xi)$ all
converge to the same weak limit $Z_\ell\sim\Normal(-\ell^2,2\ell
^2)$ as the dimension of the noise
$\xi$ goes to $\infty$.
In the
rest of this subsection, we rigorize this argument
by deriving a~Berry--Essen
bound for the weak convergence
of $R(x,\xi)$ to $Z_\ell$. 


For this purpose,
it is natural and convenient to obtain these bounds in the Wasserstein metric.
Recall that the Wasserstein\vadjust{\goodbreak} distance
between two random variables $\operatorname{Wass}(X,Y)$ is defined by
\[
\operatorname{Wass}(X,Y) \eqdef\sup_{f \in\mathrm{Lip}_1} \mathbb{E}\bigl(f(X) - f(Y)\bigr),
\]
where $\mathrm{Lip}_1$ is the class of 1-Lipschitz functions.
%
The following lemma 
gives a~bound for
the Wasserstein distance between $R(x,\xi)$ and $Z_{\ell}$.
\begin{lemma}\label{lemwass2} Almost surely
with respect to $x \sim\pi$,
%
%
\begin{eqnarray}
\label{eqnWasswz}\qquad
\operatorname{Wass}(R(x,\xi),Z_{\ell}) &\leq&
M\Biggl(\frac{1}{N^{3/2}}\sum_{j=1}^{N}|\zeta_j|^3 + \biggl|1 - \frac
{\|\zeta\|^2}{N} \biggr|
+ \frac{1}{\sqrt{N}} \Biggr) , \\
\qquad
\label{eqnestst1}\qquad
\operatorname{Wass}(R(x,\xi),R_i(x,\xi)) &\leq& \frac{M}{\sqrt{N}}
(|\zeta_i| + 1) .
\end{eqnarray}
\end{lemma}

\begin{pf}
Define the Gaussian random variable
$G \eqdef-\sqrt{\frac{2\ell^2}{N}}\sum_{k=1}^N\zeta_k \xi_k -
\ell^2$.
For any 1-Lipschitz function $f$,
\[
\bigl|\mathbb{E}^{\xi}\bigl(f(G) - f(R(x,\xi))\bigr)\bigr| \leq\ell^2 \mathbb{E}^{\xi
} \Biggl|1 - \frac{1}{N} \sum_{k=1}^N \xi_k^2 \Biggr| < M\frac
{1}{\sqrt{N}}
\]
implying that $\operatorname{Wass}(G, R(x,\xi)) \leq M\frac{1}{\sqrt{N}}$.
Now, from classical Berry--Esseen estimates (see~\cite{Stro93}), we
have that
\[
\operatorname{Wass}(G,Z_{\ell}) \leq M\frac{1}{N^{3/2}}\sum_{j=1}^{N}
|\zeta_j| ^3 +
M\biggl|1 - \frac{\|\zeta\|^2}{N} \biggr| .
\]
Hence the proof of the first claim follows from
the triangle inequality. To see the second claim,
notice that for any 1-Lipschitz function $f$ we have
\[
\mathbb{E}^\xi_0|f(R(x,\xi)) - f(R_i(x,\xi))| \leq
\mathbb{E}^\xi_0|R(x,\xi) - R_i(x,\xi)| \leq M\frac{1}{\sqrt
{N}}(1 +|\zeta_i|)
\]
and the proof is complete.
\end{pf}

Hence, from equations~(\ref{eqnestst1}) and~(\ref{eqnWasswz}),
we obtain
%
%
\begin{eqnarray}\label{eqnWasswiz}
&&
\operatorname{Wass}(R_i(x,\xi) , Z_\ell) \nonumber\\[-8pt]\\[-8pt]
&&\qquad\leq M \Biggl( \frac{1}{\sqrt
{N}}(|\zeta_i| +1) +
\frac{1}{N^{3/2}}\sum_{j=1}^{N}|\zeta_j|^3 + \biggl|1 - \frac{\|\zeta\|
^2}{N} \biggr| \Biggr).\nonumber
\end{eqnarray}
We conclude this section with
the following observation which will be used later.
Recall the Kolmogorov--Smirnov (KS) distance between two random
variables $(W,Z)$:
%
%
\begin{equation}\label{eqnksdis}
\operatorname{KS}(W,Z) \eqdef\sup_{t \in\mathbb{R}}| \mathbb{P}(W \leq t)
- \mathbb{P}(Z \leq t)| .\vadjust{\goodbreak}
\end{equation}

\begin{lemma}\label{lemkswass}
If a~random variable $Z$ has a~density with respect to the Lebesgue measure,
bounded by a~constant $M$, then
%
%
\begin{equation}\label{eqnkswassbound}
\mathrm{KS}(W,Z) \leq\sqrt{4M \operatorname{Wass}(W,Z)} .
\end{equation}
\end{lemma}

We could not find the reference for the above in any published literature,
so we include a~short proof here which was taken from the unpublished
lecture notes~\cite{Chat07}.
\begin{pf*}{Proof of Lemma~\ref{lemkswass}}
Fix $t \in\mathbb{R}$ and $\epsilon> 0$. Define two functions $g_1$
and $g_2$ as
$g_1 (y) = 1$ for $ y \in( -\infty, t)$, $g_1(y)= 0$ for
$ y \in[t + \epsilon, \infty)$
and linear interpolation in between. Similarly, define
$g_2 (y) = 1$, for $y \in(-\infty, t- \epsilon]$, $g_2(y)= 0$, for
$y \in[t, \infty)$
and linear interpolation in between. Then $g_1$ and $g_2$ form upper
and lower envelopes for the function $1_{(-\infty,t]}(y)$. So
\[
\mathbb{P} (W \leq t) - \mathbb{P} (Z \leq t) \leq\mathbb{E} g_1 (W
) -
\mathbb{E} g_1 (Z ) +\mathbb{E} g_1 (Z ) - \mathbb{P} (Z \leq T ).
\]
Since $g_1$ is $\frac{1}{\epsilon}$-Lipschitz, we have $\mathbb{E}
g_1 (W )- \mathbb{E} g_1 (Z ) \leq
\frac{1}{\epsilon} \operatorname{Wass}(W, Z )$ and $\mathbb{E} g_1 (Z )
-\mathbb{P} (Z \leq t)
\leq M \epsilon$ since Z has density bounded by $M$.
Similarly, using the function $g_2$, it follows that the same bound
holds for the difference $ \mathbb{P} (Z \leq t) - \mathbb{P} (W \leq
t)$. Optimizing over
$\epsilon$ yields the required bound.
\end{pf*}

\subsection{\texorpdfstring{Rigorous estimates for the drift: Proof of Proposition~\protect\ref{thmdrift-diffus}, equation~(\protect\ref{thmDD-drift})}
{Rigorous estimates for the drift: Proof of Proposition 2.1, equation (2.14)}}
\label{secdrifes}

In the following series of lemmas we retrace the
arguments from Section~\ref{secdiftHuristic} while deriving explicit bounds
for the error terms. Lemma~\ref{lemremest} at the end of the section
gives control of the error terms. 

The following lemma shows that $Q(x,\xi)$ is well approximated by
$R_i(x,\xi) - \sqrt{\frac{2\ell^2}{N}} \zeta_i \xi_i$, as
indicated in~(\ref{eqnQRiapprox}).
\begin{lemma}\label{lemdo1}
\begin{eqnarray*}
N \mathbb{E}_0(x^{1}_i - x_i) &=&
\lambda_i \sqrt{2\ell^2 N} \mathbb{E}_0^\xi
\bigl(\bigl(1 \wedge e^{R_i(x,\xi) - \sqrt{{2\ell^2}/{N}}\zeta_i\xi
_i}\bigr)\xi_i\bigr)
+ \omega_0(i),\\
|\omega_0(i)| &\leq& \frac{M}{\sqrt{N}}
\lambda_i .
\end{eqnarray*}
\end{lemma}
\begin{pf}
We have
\begin{eqnarray*}
N \mathbb{E}_0(x^{1}_i - x^0_i) &=&
N \mathbb{E}_0\bigl(\gamma^0 (y^0_i - x_i)\bigr) = N \mathbb
{E}_0^\xi\Biggl(\alpha(x,\xi) {\sqrt{\frac{2\ell
^2}{N}}}(C^{1/2}\xi
)_i\Biggr) \\
&=& \lambda_i \sqrt{2\ell^2 N} \mathbb{E}_0^\xi
(\alpha(x,\xi) \xi_i) =\lambda_i \sqrt{2\ell^2 N}
\mathbb{E}_0^\xi\bigl(\bigl(1 \wedge e^{Q(x,\xi)}\bigr)\xi_i\bigr) .
\end{eqnarray*}
Now we observe that
\[
\mathbb{E}_0^\xi\bigl(\bigl(1 \wedge e^{Q(x,\xi)}\bigr)\xi_i\bigr)
= \mathbb{E}_0^\xi\bigl(\bigl(1 \wedge e^{R_i(x,\xi) - \sqrt{
{2\ell^2}/{N}}\xi_i\zeta_i}\bigr)\xi_i\bigr)
+ \frac{\omega_0(i)}{\lambda_i {\sqrt{2\ell^2N}}} .
\]
By~(\ref{eqnrx,y,xi}) and~(\ref{eqnQRiapprox}),
%
%
\begin{equation}\label{eqnQRxest}
\Biggl|Q(x,\xi) - R_i(x,\xi) + {\sqrt{\frac{2\ell^2}{N}}} \zeta_i \xi
_i\Biggr|^2 \leq
\frac{M}{N^2} (|\xi_i|^4 + \|C^{1/2}\xi\|_s^4).
\end{equation}
Noticing that the map $y \mapsto1\wedge e^y$ is
Lipschitz, we obtain
\begin{eqnarray*}
|\omega_0(i)| &\leq&
M \lambda_i {\sqrt{N}} \mathbb{E}_0^{\xi}\bigl|\bigl(\bigl(1 \wedge
e^{Q(x,\xi)}\bigr)
- \bigl(1 \wedge e^{R_i(x,\xi) - \sqrt{{2\ell^2}/{N}}\xi_i\zeta
_i}\bigr)\bigr) \xi_i\bigr|\\
&\leq& M \lambda_i {\sqrt{N}} \Biggl[\mathbb{E}_0^{\xi}
\Biggl|Q(x,\xi
) -
R_i(x,\xi) + {\sqrt{\frac{2\ell^2}{N}}}\xi_i\zeta_i\Biggr|^2
\Biggr]^{1/2}[
\mathbb{E}^{\xi}_0 (\xi_i)^2]^{1/2} \\
&\leq&\frac{M}{\sqrt
{N}}\lambda_i,
\end{eqnarray*}
where the last inequality follows from~(\ref{eqnQRxest})
and the proof is complete.
\end{pf}

The next lemma takes advantage of the fact that $R_i(x,\xi)$ is independent
of $\xi_i$ conditional on $x$. Thus, using the identity~(\ref{lemnormmo}),
we obtain the bound for the approximation made in~(\ref{eqnpenexes1}).
\begin{lemma} \label{lemdo2}
%
%
\begin{eqnarray}\label{eqnpenexes}\quad
&&\hspace*{11pt}\mathbb{E}_0^\xi
\bigl(\bigl(1 \wedge e^{R_i(x,\xi) - {\sqrt{{2\ell^2}/{N}}}\zeta
_i\xi_i}\bigr)\xi_i\bigr) \nonumber\\
&&\hspace*{11pt}\qquad=
-{\sqrt{\frac{2\ell^2}{N}}}\zeta_i \mathbb{E}_0^{\xi_i^-}
e^{R_i(x,\xi) +
{\ell^2 {\zeta_i}^2}/{N}} \Phi\biggl(\frac{-R_i(x,\xi)}
{\sqrt{{2\ell^2}/{N}}|\zeta_i|} \biggr) + \omega_1(i),
\\
&&|\omega_1(i)| \leq
M|\zeta_i|^2\frac{1}{N}e^{({\ell^2}/{N})\|\zeta\|^2}
.\nonumber
\end{eqnarray}
\end{lemma}
\begin{pf}
Applying~(\ref{lemnormmo}) with $a =
-\sqrt{\frac{2\ell^2}{N}}\zeta_i$, $z = \xi_i$ and $b = R_i(x,\xi)$,
we obtain the identity
%
%
\begin{eqnarray}\label{eqntempoid}\qquad
&&\mathbb{E}_0^\xi\bigl(\bigl(1 \wedge e^{R_i(x,\xi) - \sqrt{{2\ell
^2}/{N}}\xi_i\zeta_i}\bigr)\xi_i\bigr)\nonumber\\[-8pt]\\[-8pt]
&&\qquad=
-{\sqrt{\frac{2\ell^2}{N}}}\zeta_i \mathbb{E}_0^{\xi_i^-}
e^{R_i(x,\xi) + ({\ell^2/N} )\zeta^2_i} \Phi
\Biggl({ \frac{-R_i(x, \xi)}
{{\sqrt{{2\ell^2}/{N}}}|\zeta_i|}} -
{\sqrt{\frac{2\ell^2}{N}}}|\zeta_i|\Biggr) .\nonumber
\end{eqnarray}
Now we observe that
%
%
\begin{eqnarray}\label{eqnesteR}
&&\mathbb{E}_0^{\xi_i^-} e^{R_i(x,\xi) +
{\ell^2 {\zeta_i}^2}/{N}} \nonumber\\
&&\qquad=
\mathbb{E}^{\xi_i^-}_0 \bigl(e^{-\sqrt{{2\ell^2}/{N}}\sum
_{j=1,j\neq i}^N \zeta_j \xi_j -
({\ell^2}/{N})\sum_{j=1, j \neq i}^N {\xi_j}^2 + ({\ell
^2}/{N}){\zeta_i}^2}\bigr) \\
&&\qquad\leq\mathbb{E}^{\xi_i^-}_0 \bigl(e^{-\sqrt{{2\ell^2}/{N}}\sum
_{j=1,j\neq i}^N \zeta_j\xi_j +
({\ell^2}/{N}){\zeta_i}^2}\bigr) = e^{({\ell^2}/{N})\|\zeta\|^2} .
\nonumber
\end{eqnarray}
Since $\Phi$ is globally Lipschitz, it follows that
%
%
\begin{eqnarray}\label{eqnesto11}
&&\hspace*{11pt}\mathbb{E}_0^{\xi_i^-}
e^{R_i(x,\xi) +
{\ell^2 {\zeta_i}^2}/{N}} \Phi\Biggl(-\frac{R_i(x,\xi)}
{\sqrt{ {2\ell^2}/{N}}|\zeta_i|} -
\sqrt{ \frac{2\ell^2}{N}}|\zeta_i|\Biggr)\nonumber\\[-2pt]
&&\hspace*{11pt}\qquad= \mathbb{E}_0^{\xi_i^-} e^{R_i(x,\xi) +
{\ell^2 {\zeta_i}^2}/{N}} \Phi\biggl(\frac{-R_i(x,\xi)}
{\sqrt{ {2\ell^2}/{N}}|\zeta_i|}\biggr) + \omega_1(i), \\[-2pt]
&&|\omega_1(i)| \leq M|\zeta_i|\frac{1}{\sqrt{N}}
\mathbb{E}_0^{\xi_i^-} e^{R_i(x,\xi) +
{\ell^2 {\zeta_i}^2}/{N}}
\leq M|\zeta_i|\frac{1}{\sqrt{N}}e^{({\ell^2}/{N})\|\zeta\|
^2,} \nonumber
\end{eqnarray}
where the last estimate follows from~(\ref{eqnesteR}).
The lemma follows from~(\ref{eqntempoid}) and~(\ref{eqnesteR}).\vspace*{-2pt}
\end{pf}

The next few lemmas are technical and give quantitative bounds for the
approximations
in~(\ref{eqnRilimwbet1}) and~(\ref{eqnRilimwbet2}).\vspace*{-2pt}
\begin{lemma} \label{lemdo3}
\begin{eqnarray*}
&&\hspace*{11pt}\mathbb{E}_0^{\xi_i^-} e^{R_i(x,\xi) +
{\ell^2 {\zeta_i}^2}/{N}} \Phi\biggl(\frac{-R_i(x,\xi)}
{\sqrt{{2\ell^2}/{N}}|\zeta_i|}\biggr)\\[-2pt]
&&\hspace*{11pt}\qquad= \mathbb{E}_0^{\xi_i^-} e^{R_i(x,\xi) +
{\ell^2 {\zeta_i}^2}/{N}}1_{{R}_i(x,\xi) < 0} +
\omega_2(i),\\[-2pt]
&&|\omega_2(i)| \leq Me^{({2\ell^2}/{N})\|\zeta\|^2}(|\zeta_i| +
1)\biggl[\mathbb{E}^{\xi}_0\frac{1}
{(1 +|R(x,\xi)|\sqrt{N})^2}\biggr]^{1/4} .\vspace*{-2pt}
\end{eqnarray*}
\end{lemma}
\begin{pf}
We first prove the following lemma needed for the proof.
\begin{lemma}\label{lemmodMillsratio} Let $\phi(\cdot)$ and
$\Phi(\cdot)$ denote the pdf and CDF of the standard normal
distribution, respectively. Then we have:
\begin{longlist}[(1)]
\item[(1)] for any $x \in\bbR$, $|\Phi(-x) - 1_{x<0}| = |1 - \Phi(|x|)|$.
\item[(2)] for any $x>0$ and $\epsilon\geq0$, $1
- \Phi(x) \leq\frac{1+\epsilon}{x + \epsilon}$.\vspace*{-2pt}
\end{longlist}
\end{lemma}
\begin{pf}
For the first claim, notice that if $x > 0$, $|\Phi(-x) - 1_{x<0}| =\break
|\Phi(-x)|
=  |1 - \Phi(|x|)|$. If $x <0$, $|\Phi(-x) - 1_{x<0}| = |1 - \Phi
(|x|)|$ and the claim
follows.

For the second claim,
\[
1 - \Phi(x) = \int_{x}^{\infty} \phi(u) \,du\leq\int_{x}^{\infty}
\frac{u + \epsilon}{x + \epsilon} \phi(u) \,du\leq\frac{\phi
(x) + \epsilon}{x+ \epsilon} \leq\frac{1+\epsilon}{x + \epsilon}
\]
since $\int_{-\infty}^{\infty} \phi(u) \,du = 1$.\vspace*{-2pt}
\end{pf}
%

We now proceed to the proof of Lemma~\ref{lemdo3}. By
Cauchy--Schwarz and
an estimate similar to~(\ref{eqnesteR}), 
%
%
\begin{eqnarray}\label{eqnestcomp}
|\omega_2(i)| &\leq& \mathbb{E}^{\xi_i^-}_0\biggl[e^{R_i(x,\xi) +
({\ell^2}/{N}){\zeta_i}^2} \biggl|
1_{R_i(x,\xi) <0}- \Phi\biggl( \frac{-R_i(x,\xi)}{{\sqrt
{{2\ell^2}/{N}}}|\zeta_i|}\biggr)\biggr| \biggr] \nonumber\\
&\leq& \bigl[\mathbb{E}^{\xi_i^-}_0e^{2R_i(x,\xi) +
(2{\ell^2}/{N}){\zeta_i}^2}\bigr]^{1/2}\biggl[ \mathbb{E}^{\xi
_i^-}_0\biggl|
1_{R_i(x,\xi) <0}- \Phi\biggl(\frac{-R_i(x,\xi)}
{\sqrt{{2\ell^2}/{N}}|\zeta_i|}\biggr)\biggr|^2 \biggr]^{1/2}
\nonumber\\[-8pt]\\[-8pt]
&\leq& Me^{({2\ell^2}/{N})\|\zeta\|^2} \mathbb{E}^{\xi_i^-}_0
\biggl[\biggl|
1_{R_i(x,\xi) <0}- \Phi\biggl(\frac{-R_i(x,\xi)}
{{\sqrt{{2\ell^2}/{N}}}|\zeta_i|}\biggr)\biggr|^2 \biggr]^{1/2}
\nonumber\\
&\leq& Me^{({2\ell^2}/{N})\|\zeta\|^2} \biggl[\mathbb{E}^{\xi
_i^-}_0\biggl|
1_{R_i(x,\xi) <0}- \Phi\biggl(\frac{-R_i(x,\xi)}
{{\sqrt{{2\ell^2}/{N}}}|\zeta_i|}\biggr)\biggr| \biggr]^{1/2},
\nonumber
\end{eqnarray}
where the last two observations follow from the computation done in
(\ref{eqnesteR})
and the fact that $|1_{R_i(x,\xi) <0}- \Phi(\frac
{-R_i(x,\xi)}
{\sqrt{{2\ell^2}/{N}}|\zeta_i|})| < 1$.

By applying Lemma~\ref{lemmodMillsratio}, with $\epsilon
= \frac{1}{\sqrt{2}\ell|\zeta_i|}$,
%
%
\begin{eqnarray}\label{eqnusingcuteadbd}
\biggl| 1_{R_i(x,\xi) <0}- \Phi\biggl({ \frac{-R_i(x,\xi)}
{\sqrt{{2\ell^2}/{N}}|\zeta_i|}}\biggr)\biggr| &=& 1-
\Phi\biggl({ \frac{|R_i(x,\xi)|}
{\sqrt{{2\ell^2}/{N}|\zeta_i|}}}\biggr)\nonumber\\
&=& 1-
\Phi\biggl({ \frac{|R_i(x,\xi)|\sqrt{N}}
{\sqrt{2}\ell|\zeta_i|}}\biggr)\\
&\leq&\bigl(1 + \sqrt{2}\ell
|\zeta_i|\bigr)\frac{1}{1 +|R_i(x,\xi)|\sqrt{N}} .
\nonumber
\end{eqnarray}
%
The right-hand side of the estimate~(\ref{eqnusingcuteadbd}) depends on
$i$ but we need estimates which are
independent of $i$. In the next lemma, we replace $R_i(x,\xi)$ by
$R(x,\xi)$ and control the extra error term. 
%
\begin{lemma}
%
%
\begin{equation}\label{eqnerrgamdo3}\qquad
\mathbb{E}^{\xi_i^-}_0\frac{1}{1 +|R_i(x,\xi)|\sqrt{N}} \leq
M (1 + |\zeta_i|) \biggl[\mathbb{E}^{\xi}_0\frac{1}
{(1 +|R(x,\xi)|\sqrt{N})^2}\biggr]^{1/2} .
\end{equation}
\end{lemma}
\begin{pf}
We write
%
%
\begin{eqnarray}\label{eqndo3ga1}\quad
\mathbb{E}^{\xi_i^-}_0 \frac{1}{1 +|R_i(x,\xi)|\sqrt{N}} &=&
\mathbb{E}^{\xi}_0\frac{1}{1 +|R_i(x,\xi)|\sqrt{N}} \nonumber\\
&=& \mathbb{E}^{\xi}_0\frac{1}{1 +|R(x,\xi)|\sqrt{N}}+ \gamma\\
&\leq&\biggl[\mathbb{E}^{\xi}_0\frac{1}
{(1 +|R(x,\xi)|\sqrt{N})^2}\biggr]^{1/2}
+ \gamma,\nonumber\\
\label{eqndo3ga2}
|\gamma| &\leq&\mathbb{E}^{\xi}_0\biggl|\frac{1}{1 +|R_i(x,\xi
)|\sqrt{N}} -
\frac{1}{1 +|R(x,\xi)|\sqrt{N}}\biggr|\nonumber\\
&\leq&\mathbb{E}^{\xi}_0\frac{\sqrt{2}\ell|\zeta_i| |\xi_i| +
{\ell^2}/{\sqrt{N}}{\xi_i}^2}
{(1 +|R_i(x,\xi)|\sqrt{N})(1 +|R(x,\xi)|\sqrt{N})}
\nonumber\\[-8pt]\\[-8pt]
&\leq&\mathbb
{E}_0^{\xi} \frac{\sqrt{2}\ell|\zeta_i| |\xi_i| + {\ell
^2}/{\sqrt{N}}{\xi_i}^2}
{(1 +|R(x,\xi)|\sqrt{N})}\nonumber\\[-2pt]
&\leq& M(|\zeta_i| + 1) \biggl[\mathbb{E}^{\xi}_0\frac{1}
{(1 +|R(x,\xi)|\sqrt{N})^2}\biggr]^{1/2},\nonumber
\end{eqnarray}
and the claim follows from~(\ref{eqndo3ga1}) and~(\ref{eqndo3ga2}).
\end{pf}
%
Now, by applying the estimates obtained in~(\ref{eqnestcomp}),
(\ref{eqnusingcuteadbd}) and~(\ref{eqnerrgamdo3}), we obtain
\[
|\omega_2(i)| \leq Me^{({2\ell^2}/{N})\|\zeta\|^2}(|\zeta_i| +
1)\biggl[\mathbb{E}^{\xi}_0\frac{1}
{(1 +|R(x,\xi)|\sqrt{N})^2}\biggr]^{1/4}
\]
and the proof is complete.
\end{pf}

The error estimate in $\omega_2$ has
$R(x,\xi)$ instead of $R_i(x,\xi)$. This bound can be achieved
because the terms
$R_i(x,\xi)$ for all $i \in\mathbb{N}$ have the same weak limit as
$R(x,\xi)$
and thus the additional error term due to the replacement of $R_i(x,\xi
)$ by
$R(x,\xi)$ in the expression 
can be controlled uniformly over $i$ for large~$N$.
\begin{lemma} \label{lemdo4}
\begin{eqnarray*}
\mathbb{E}_0^{\xi_i^-} e^{R_i(x,\xi) +
{\ell^2 {\zeta_i}^2}/{N}}1_{R_i(x,\xi) < 0} &=& \frac{\beta
}{2} + \omega_3(i) ,\\[-20pt]
\end{eqnarray*}
\begin{eqnarray*}
|\omega_3(i)| &\leq& M\frac{\zeta_i^2}{N} e^{\ell^2
{\|\zeta\|^2}/{N}} \\[-2pt]
&&{}+ M\Biggl(\frac{1+|\zeta_i| }{\sqrt{N}} +
\frac{1}{N^{3/2}}\sum_{j=1}^{N} |\zeta_j|^3 + \biggl|1 - \frac{\|\zeta\|
^2}{N}\biggr| \Biggr)^{1/2} .
\end{eqnarray*}
\end{lemma}
\begin{pf}
Set 
$g(y) \eqdef e^{y} 1_{y < 0}$.
We first need to estimate the following:
\[
\bigl|\mathbb{E}_0^{\xi}\bigl(g(R_i(x,\xi)) - g(Z_{\ell})\bigr)\bigr| .
\]
Notice that the function $g(\cdot)$ is not Lipschitz and therefore,
the Wasserstein bounds obtained earlier cannot be used directly.
However, we use the fact that the normal distribution has a~density
which is bounded above.
So by Lemma~\ref{lemwass2},~(\ref{eqnWasswiz}) and~(\ref{eqnkswassbound}),
\begin{eqnarray*}
\operatorname{KS}(R_i(x,\xi),Z_{\ell}) &\leq&2M\sqrt{\operatorname{Wass}(R_i(x,\xi
),Z_{\ell})} \\[-2pt]
&\leq& M\Biggl(\frac{1+|\zeta_i| }{\sqrt{N}} +
\frac{1}{N^{3/2}}\sum_{j=1}^{N} |\zeta_j|^3 + \biggl|1 - \frac{\|\zeta\|
^2}{N} \biggr|\Biggr)^{1/2} .\vadjust{\goodbreak}
\end{eqnarray*}
Since $g$ is positive on $(-\infty,0] $, for a~real valued
continuous random variable~$X$,
\[
\mathbb{E}(g(X)) = \int_{-\infty}^0 g'(t)\bigl(\mathbb{P}(X > t)
\bigr)\,dt - g(0) \mathbb{P}(X \geq0) .\vspace*{-2pt}
\]
Hence,
\begin{eqnarray*}
|\mathbb{E}_0^{\xi} g(R_i(x,\xi)) - \mathbb{E} g(Z_{\ell})| &\leq&
\biggl|\int_{-\infty}^0 g'(t)\bigl(\mathbb{P}\bigl(R_i(x,\xi) > t\bigr) -
\mathbb{P}(Z_\ell>t)\bigr)\,dt\biggr|\\[-1pt]
&&{} + g(0)\bigl|\mathbb{P}\bigl(R_i(x,\xi) \geq0\bigr) - \mathbb{P}(Z_\ell\geq0)\bigr|\\[-1pt]
&\leq&\operatorname{KS}(R_i(x,\xi),Z_\ell)\biggl( \int_{-\infty}^0g'(t)\,dt +
g(0)\biggr) \\[-1pt]
&\leq& M\operatorname{KS}(R_i(x,\xi),Z_\ell).\vspace*{-2pt}
\end{eqnarray*}

Hence, putting the above calculations together and noticing that\break
$\mathbb{E}(e^{Z_\ell}1_{Z_\ell< 0}) =
\beta/2 $, we have just shown that
\[
\biggl|\mathbb{E}_0^\xi\bigl(e^{R_i(x,\xi)}1_{R_i(x,\xi)<0}\bigr) - \frac
{\beta}{2}\biggr| \leq M\sqrt{\frac{1+|\zeta_i| }{\sqrt{N}} +
\frac{1}{N^{3/2}}\sum_{j=1}^{N} |\zeta_j|^3+ \biggl|1 - \frac{\|\zeta\|
^2}{N} \biggr|}.\vspace*{-2pt}
\]
Notice that
\begin{eqnarray*}
|\omega_3(i)| &\leq& \bigl|e^{\ell^2 {\zeta^2_i}/{N}}\mathbb
{E}_0^\xi\bigl(e^{R_i(x,\xi)}1_{R_i(x,\xi)<0}\bigr) - \beta/2\bigr| \\[-1pt]
&\leq& |e^{\ell^2 {\zeta^2_i}/{N}}-1| \bigl|\mathbb{E}_0^\xi
\bigl(e^{R_i(x,\xi)}1_{R_i(x,\xi)<0}\bigr)\bigr|\\[-1pt]
&&{}+
\bigl|\mathbb{E}_0^\xi\bigl(e^{R_i(x,\xi)}1_{R_i(x,\xi)<0}\bigr) - \beta/2\bigr| \\[-1pt]
&\leq& M\frac{\zeta_i^2}{N} e^{\ell^2 {\|\zeta\|^2}/{N}}
+\bigl|\mathbb{E}_0^\xi\bigl(e^{R_i(x,\xi)}1_{R_i(x,\xi)<0}\bigr) - \beta/2\bigr|,
\end{eqnarray*}
where the last bound follows from~(\ref{eqnesteR}),
proving the claimed error bound for $\omega_3(i)$.\vspace*{-2pt}
\end{pf}

For deriving the error bounds on $\omega_3$, we cannot
directly apply the Wasserstein bounds obtained in~(\ref{eqnWasswiz}),
because the function $ y \mapsto e^y1_{y<0}$ is not Lipschitz on
$\mathbb{R}$.
However, using~(\ref{eqnkswassbound}), the KS distance
between $R_i(x,\xi)$ and~$Z_\ell$ is bounded by the square
root of the Wasserstein
distance. Thus, using the
fact that $e^y1_{y<0}$ is bounded and positive, we bound
the expectation in Lemma~\ref{lemdo4} by the
KS distance.

Combining all the above estimates, we see that
%
%
\begin{equation}
N\mathbb{E}_0^\xi[x^{1}_i - x_i ] =
-\ell^2 \beta\bigl(P^N x + C \nabla\Psi(P^N x)\bigr)_i + r_i^{N} %
\end{equation}
with
%
%
\begin{equation}\label{eqnrNest}
|r_i^{N}| \leq|\omega_0(i)| + M \lambda_i \bigl( {\sqrt{N}}
|\omega_1(i)|
+ |\zeta_i| |\omega_2(i)|+ |\zeta_i| |\omega_3(i)|\bigr) .
\end{equation}
The following lemma gives the control over $r^N$ and completes
the proof of~(\ref{thmDD-drift}), Proposition~\ref{thmdrift-diffus}.\vadjust{\goodbreak}
\begin{lemma}\label{lemremest} For $s < \kappa-1/2$,
\[
\lim_{N \rightarrow\infty} \mathbb{E}^{\pi^N}\|r^{N}\|_s^2 = \lim
_{N \rightarrow\infty}
\mathbb{E}^{\pi^N}\sum_{i=1}^N
i^{2s} |r^{N}_i|^2 = 0.\vspace*{-2pt}
\]
\end{lemma}
\begin{pf}
By~(\ref{eqnrNest}), we have
$ |r_i^{N}| \leq|\omega_0(i)| + M \lambda_i ( {\sqrt{N}} |\omega
_1(i)| + |\zeta_i| |\omega_2(i)|
+ |\zeta_i| |\omega_3(i)|)$.
Therefore,
%
%
\begin{eqnarray}\label{eqnrdoremain}\qquad
&&
\mathbb{E}^{\pi^N}\sum_{i=1}^N i^{2s} |r^N_i|^2 \nonumber\\[-9pt]\\[-9pt]
&&\qquad\leq M \mathbb
{E}^{\pi^N}
\sum_{i=1}^N\bigl( i^{2s} |\omega_0(i)|^2+
i^{2s} \lambda_i^2\bigl(N \omega_1(i)^2 + {\zeta_i} ^2 \omega_2(i)^2 +
{\zeta_i} ^2 \omega_3(i)^2 \bigr) \bigr) .\nonumber\vspace*{-2pt}
\end{eqnarray}
%
Now we will evaluate each sum of the right-hand side of the above equation
and show that they converge to zero.
%
\begin{itemize}
\item Since $\sum_{i=1}^\infty\lambda_i^2 i^{2s} < \infty$,
%
%
\begin{equation}
\label{eqnestom0} \sum_{i=1}^N
\mathbb{E}^{\pi^N}i^{2s} |\omega_0(i)|^2 \leq M
\frac{1}{N}\sum_{i=1}^N i^{2s} \lambda_i^2 \leq M \frac{1}{N}
\sum_{i=1}^\infty\lambda_i^2 i^{2s} \rightarrow0.\vspace*{-2pt}
\end{equation}
\item By Lemmas~\ref{lemdo2} and~\ref{lemshiftmeas},
%
%
\begin{equation}\label{eqnestom1}\quad
N \mathbb{E}^{\pi^N}
\sum_{i=1}^N \lambda_i^2 i^{2s} |\omega_1(i)|^2 \leq M\frac{1}{N}
\sum_{i=1}^N\mathbb{E}^{\pi^N} \lambda_i^2 i^{2s} |\zeta_i|^4
e^{({2\ell^2}/{N})\|\zeta\|^2} \rightarrow0.\vspace*{-2pt}
\end{equation}
\item From Lemma~\ref{lemdo3} and Cauchy--Schwarz, we obtain
\begin{eqnarray*}
&&\sum_{i=1}^N \mathbb{E}^{\pi^N} \lambda_i^2 i^{2s} |\zeta_i|^2
|\omega_2(i)|^2\\[-2pt]
&&\qquad\leq M \biggl( \mathbb{E}^{\pi^N}\biggl[\mathbb{E}^{\xi}_0 \frac{1} {(1
+|R(x,\xi)|\sqrt{N})^2}\biggr]\biggr)^{1/2}\\[-2pt]
&&\qquad\quad{}\times\sum_{i=1}^N\bigl(\mathbb
{E}^{\pi^N}
e^{({8\ell^2}/{N})\|\zeta\|^2}\lambda_i^4
i^{4s} (|\zeta_i|^8+1)\bigr)^{1/2} .\vspace*{-2pt}
\end{eqnarray*}
Proceeding similarly as in Lemma~\ref{lemshiftmeas}, it follows that
\[
\sum_{i=1}^N\bigl(\mathbb{E}^{\pi^N} e^{({8\ell^2}/{N})\|\zeta\|
^2}\lambda_i^4 i^{4s} (|\zeta_i|^8+1)\bigr)^{1/2}\vspace*{-2pt}
\]
is bounded in~$N$.
Since, with $x \dist\pi_0$, $R(x,\xi)$ converges weakly to $Z_\ell$
as $N \rightarrow\infty$,
by the bounded convergence theorem we obtain
\[
\lim_{N\rightarrow\infty}\mathbb{E}^{\pi_0} \biggl[\mathbb
{E}^{\xi}_0 \frac{1}
{(1 +|R(x,\xi)|\sqrt{N})^2}\biggr] = 0\vadjust{\goodbreak}
\]
and thus, by Lemma~\ref{lemmchmeas},
\[
\lim_{N\rightarrow\infty}\mathbb{E}^{\pi^N} \biggl[\mathbb{E}^{\xi
}_0 \frac{1}
{(1 +|R(x,\xi)|\sqrt{N})^2}\biggr] = 0.
\]
Therefore, we deduce that
%
%
\begin{equation}
\label{eqnestom2}
\lim_{N\rightarrow\infty}\sum_{i=1}^N \mathbb{E}^{\pi^N} |\zeta
_i|^2 i^{2s} \lambda_i^2 |\omega_2(i)|^2 = 0.
\end{equation}
\item After some algebra we obtain from Lemma~\ref{lemdo4} that
\begin{eqnarray*}
\hspace*{-4pt}&&\mathbb{E}^{\pi^N}\sum_{i=1}^N \lambda_i^2 i^{2s} |\zeta_i |^2
|\omega_3(i)|^2 \\
\hspace*{-4pt}&&\qquad\leq M \frac{1}{N^2}
\sum_{i=1}^N\mathbb{E}^{\pi^N}
\lambda_i^2 i^{2s} |\zeta_i|^6 e^{2\ell^2 ({\|\zeta\|^2}/{N})}
+ M\frac{1}{\sqrt{N}} \mathbb{E}^{\pi^N}\sum_{i=1}^N\lambda_i^2
i^{2s} \zeta_i^2(1 +|\zeta_i| ) \\
\hspace*{-4pt}&&\qquad\quad{} + M \Biggl[\Biggl(\mathbb{E}^{\pi^N}\Biggl(\frac{1}{N^{3/2}}\sum_{j=1}^{N}
|\zeta_j|^3
\Biggr)^2 + \mathbb{E}^{\pi^N}\biggl|1 - \frac{\|\zeta\|^2}{N}\biggr|^2
\Biggr)^{1/2} \Biggr] \\
\hspace*{-4pt}&&\qquad\quad\hspace*{10pt}{}\times\sum_{i=1}^N
(\mathbb{E}^{\pi^N} \lambda_i^4 i^{4s} \zeta_i^4)^{1/2} .
\end{eqnarray*}
Similar to the previous calculations, using Lemma~\ref{lemshiftmeas},
it is quite
straightforward to verify that each of the four terms above converges
to $0$.
Thus we obtain
%
%
\begin{equation}\label{eqnestom3}
\lim_{N\rightarrow\infty}\sum_{i=1}^N \mathbb{E}^{\pi^N} \lambda
_i^2 i^{2s} |\zeta_i|^2 |\omega_3(i)|^2 = 0.
\end{equation}
\end{itemize}
Now the proof of Lemma~\ref{lemremest} follows from
(\ref{eqnrdoremain})--(\ref{eqnestom3}).
\end{pf}

This completes the proof of Proposition~\ref{thmdrift-diffus},
equation~(\ref{thmDD-drift}).

\subsection{\texorpdfstring{Rigorous estimates for the diffusion coefficient: Proof of Proposition~\protect\ref{thmdrift-diffus}, equation (\protect\ref{thmDD-diffus})}
{Rigorous estimates for the diffusion coefficient: Proof of Proposition 2.1, equation (2.15)}}\label{secdiffes}
Recall that for $1 \leq i,j \leq N$,
\[
N \mathbb{E}_0[(x_i^{1} - x_i^{0})(x_j^{1} - x_j^{0})]
=2\ell^2 \mathbb{E}_0^{\xi}\bigl[(C^{1/2}\xi)_i
(C^{1/2}\xi)_j \bigl(1 \wedge\exp
Q(x,\xi)\bigr)\bigr] .
\]
The following lemma quantifies the approximations made
in~(\ref{eqnvardiffbd}) and~(\ref{eqnrconvbet}).

\begin{lemma}\label{lemdrifvar}
\begin{eqnarray*}
\mathbb{E}_0^{\xi}\bigl[(C^{1/2}\xi)_i (C^{1/2}\xi)_j \bigl(1
\wedge
\exp Q(x,\xi)\bigr)\bigr]&=&\lambda_i\lambda_j
\delta_{ij} \mathbb{E}^{\xi^-_{ij}}\bigl[
\bigl(1 \wedge\exp R_{ij}(x,\xi)\bigr)\bigr] + \theta_{ij},\\
\mathbb{E}^{\xi^-_{ij}}\bigl[ \bigl(1 \wedge\exp
R_{ij}(x,\xi)\bigr)\bigr] &=& \beta+ \rho_{ij},
\end{eqnarray*}
where the error terms satisfy
%
%
\begin{eqnarray}\label{eqnrhoij}\quad
|\theta_{ij}| &\leq& M\lambda_{i}\lambda_j(1 + |\zeta_i| ^2+
|\zeta_j|^2)^{1/2}\frac{1}{\sqrt{N}}, \\
\label{eqnkapp}
|\rho_{ij}| &\leq& M\Biggl(\frac{1}{\sqrt{N}}(1 +|\zeta_i| + |\zeta_j|)
+ \frac{1}{N^{3/2}}\sum_{s=1}^{N}|\zeta_s|^3 + \biggl|1 - \frac{\|\zeta\|
^2}{N} \biggr|\Biggr) .
\end{eqnarray}
\end{lemma}
\begin{pf}
We first
derive the bound for $\theta$. Indeed,
\begin{eqnarray*}
|\theta_{ij}| &\leq& \mathbb{E}_0^{\xi}\bigl[\bigl|
(C^{1/2}\xi)_i (C^{1/2}\xi)_j \bigl(\bigl(1 \wedge e^{Q(x,\xi)}\bigr) -
\bigl(1 \wedge e^{R_{ij}(x,\xi)}\bigr)\bigr)\bigr|\bigr] \\
&\leq& M \lambda_i\lambda_j \mathbb{E}_0^{\xi}\bigl[\bigl|
\xi_i \xi_j \bigl(\bigl(1 \wedge e^{Q(x,\xi)}\bigr) -
\bigl(1 \wedge e^{R_{ij}(x,\xi)}\bigr)\bigr)\bigr|\bigr] .
\end{eqnarray*}
By the Cauchy--Schwarz inequality,
\begin{eqnarray*}
|\theta_{ij}|&\leq& M \lambda_i\lambda_j \bigl(\mathbb{E}_0^{\xi}
\bigl|\bigl(1\wedge e^{Q(x,\xi)}\bigr) -
\bigl(1 \wedge e^{R_{ij}(x,\xi)}\bigr)\bigr|\bigr)^{1/2} \\
&\leq& M \lambda_i\lambda_j \bigl(\mathbb{E}_0^{\xi} |Q(x,\xi) -
R_{ij}(x,\xi)|^2\bigr)^{1/2} .
\end{eqnarray*}
Using the estimate obtained in~(\ref{eqnRijest}),
\[
|\theta_{ij}|\leq M \lambda_i\lambda_j(1 + |\zeta_i|^2 + |\zeta_j|^2)^{1/2}\frac
{1}{\sqrt{N}}
\]
verifying~(\ref{eqnrhoij}).

Now we turn to verifying the error bound
in~(\ref{eqnkapp}). 
We need to bound
\[
\mathbb{E}_0^{\xi}\bigl(g(R_{ij}(x,\xi)) - g(Z_\ell)\bigr),
\]
where $g(y) \eqdef1 \wedge e^{y}$. Notice that $\mathbb{E}(g(Z_\ell
)) = \beta$.
Since $g(\cdot)$ is Lipschitz,
%
%
\begin{equation}\label{eqnsteinfunvar}
\bigl|\mathbb{E}_0^{\xi}\bigl(g(R_{ij}(x,\xi)) - g(Z_\ell)\bigr)\bigr| \leq M
\operatorname{Wass}(R_{ij}(x,\xi), Z_\ell) .
\end{equation}
%
calculation will yield that
\[
\operatorname{Wass}(R_{ij}(x,\xi),R(x,\xi)) \leq M (|\zeta_i| + |\zeta
_j|+1)\frac{1}{\sqrt{N}}.
\]
Therefore, by the triangle inequality and Lemma~\ref{lemwass2},
\[
\operatorname{Wass}(R_{ij}(x,\xi),Z_\ell) \leq M\Biggl(\frac{1}{\sqrt{N}}(1
+|\zeta_i| + |\zeta_j|) +
\frac{1}{N^{3/2}}\sum_{r=1}^{N}|\zeta_r|^3 + \biggl|1 - \frac{\|\zeta\|
^2}{N} \biggr|\Biggr) .
\]
Hence the estimate in
(\ref{eqnrhoij}) follows from the observation made in (\ref
{eqnsteinfunvar}).~%
\end{pf}

Putting together all the estimates produces
%
%
\begin{eqnarray}\label{eqnputtestdiff}
N \mathbb{E}_0[(x_i^{1} - x_i^{0})(x_j^{1} - x_j^{0})] &=&
2\ell^2\beta\lambda_i \lambda_j \delta_{ij} + E^{N}_{ij}
\quad\mbox{and}\nonumber\\[-8pt]\\[-8pt]
|E^{N}_{ij}| &\leq& M(|\theta_{ij}| +
\lambda_{i}\lambda_{j} \delta_{ij}|\rho_{ij}|) .\nonumber
\end{eqnarray}
Finally we estimate the error of $E^{N}_{ij}$.\vadjust{\goodbreak}
%
\begin{lemma} \label{lemestukN}
We have
\[
\lim_{N \rightarrow\infty} \sum_{i=1}^N\mathbb{E}^{\pi^N} |
\langle\phi_i,E^N \phi_j \rangle_s
| = 0,\qquad
\lim_{N \rightarrow\infty} \mathbb{E}^{\pi^N} | \langle\phi
_i,E^N \phi_j \rangle_s
| = 0
\]
for any pair of indices $i,j$.
\end{lemma}
\begin{pf}
From~(\ref{eqnputtestdiff}) we obtain that
%
%
\begin{eqnarray}\qquad
\label{eqntrbd}
\sum_{i=1}^N\mathbb{E}^{\pi^N} | \langle\phi_i,E^N \phi_i \rangle
_s |
&\leq&
M\Biggl(\sum_{i=1}^N \mathbb{E}^{\pi^N} i^{2s} |\theta_{ii}| + \sum
_{i=1}^N \lambda_i^2 i^{2s} \mathbb{E}^{\pi^N}|\rho_{ii}|\Biggr) , \\
\label{eqnlemkappa}
\sum_{i=1}^N\mathbb{E}^{\pi^N} i^{2s} | \theta_{ii}| &\leq&
M\sum_{i=1}^N\mathbb{E}^{\pi_0}|\theta_{ii}| i^{2s}\nonumber\\
&\leq& M \sum_{i=1}^N\mathbb{E}^{\pi_0}
\lambda^2_{i} i^{2s} (1 + |\zeta_i|^2)^{1/2}\frac{1}{\sqrt
{N}}\\
&\leq& M \sum_{i=1}^N\mathbb{E}^{\pi_0}
\lambda^2_{i} i^{2s}(1 + |\zeta_i|)\frac{1}{\sqrt{N}}
\rightarrow0\nonumber
\end{eqnarray}
due to the fact that $\sum_{i=1}^\infty\lambda_{i}^2 i^{2s} < \infty
$ and Lemma~\ref{lemshiftmeas}.
Now the second term of~(\ref{eqntrbd}),
\begin{eqnarray*}
&&\sum_{i=1}^N \lambda_i^2 i^{2s} \mathbb{E}^{\pi^N} |\rho_{ii}|\\
&&\qquad\leq M\mathbb{E}^{\pi_0}\sum_{i=1}^N \lambda_i^2 i^{2s}
\Biggl(\frac{1}{\sqrt{N}}(1 +|\zeta_i|) +
\frac{1}{N^{3/2}}\sum_{s=1}^{N}|\zeta_s|^3 + \biggl|1 - \frac{\|\zeta\|
^2}{N} \biggr|\Biggr) .
\end{eqnarray*}
The first term above goes to zero by~(\ref{eqnlemkappa}) and the last
term converges to zero by the same arguments used in Lemma \ref
{lemshiftmeas}. As mentioned in the proof
of the estimate for the term $\omega_3$ in Lemma~\ref{lemremest},
the sum $\mathbb{E}^{\pi^N}\frac{1}{N^{3/2}}\sum_{s=1}^{N}|\zeta
_s|^3$ goes to zero.
Therefore, we have shown that
\[
\lim_{N \rightarrow\infty} \sum_{i=1}^N\mathbb{E}^{\pi^N}|
\langle\phi_i,E^N \phi_i \rangle_s| =0,
\]
proving the first claim.
Finally, from~(\ref{eqnrhoij}) it immediately follows that
\[
\mathbb{E}^\pi|\langle\phi_i, E^N \phi_j \rangle_s| \leq\mathbb
{E}^\pi i^{s} j^{s} |\theta_{ij}| \rightarrow0,
\]
proving the second claim as well.
\end{pf}

Therefore, we have shown
\begin{eqnarray*}
N \mathbb{E}_0[(x_i^{1} - x_i^{0})(x_j^{1} - x_j^{0})]
&=& 2\ell^2\beta\langle\phi_i, C \phi_j\rangle+ E^{N},\\
\lim_{N \rightarrow\infty} \sum_{i=1}^N\mathbb{E}^{\pi^N} |
\langle\phi_i,E^N \phi_i \rangle
|&=& 0 .
\end{eqnarray*}
%
This finishes the proof of Proposition~\ref{thmdrift-diffus},
equation~(\ref{thmDD-diffus}).

\section*{Acknowledgments}
We thank Alex Thiery and an anonymous referee for their careful reading
and very insightful comments which significantly improved the clarity of
the presentation.


%

%
\printaddresses


\begin{thebibliography}{28}

\bibitem{Beda07}
%
\begin{barticle}[mr]
\bauthor{\bsnm{B{\'e}dard},~\bfnm{Myl{\`e}ne}\binits{M.}}
(\byear{2007}).
\btitle{Weak convergence of {M}etropolis algorithms for non-i.i.d. target
distributions}.
\bjournal{Ann. Appl. Probab.}
\bvolume{17}
\bpages{1222--1244}.
\bid{doi={10.1214/105051607000000096}, issn={1050-5164}, mr={2344305}}
\end{barticle}
%
\endbibitem

\bibitem{Beda09}
%
\begin{bmisc}[auto:STB|2011-03-03|12:04:44]
\bauthor{\bsnm{B{\'e}dard},~\bfnm{M.}\binits{M.}}
(\byear{2009}).
\bhowpublished{On the optimal scaling problem of Metropolis algorithms for
hierarchical target distributions. Preprint}.
\end{bmisc}
%
\endbibitem

\bibitem{Berg86}
%
\begin{barticle}[mr]
\bauthor{\bsnm{Berger},~\bfnm{Erich}\binits{E.}}
(\byear{1986}).
\btitle{Asymptotic behaviour of a~class of stochastic approximation
procedures}.
\bjournal{Probab. Theory Related Fields}
\bvolume{71}
\bpages{517--552}.
\bid{doi={10.1007/BF00699040}, issn={0178-8051}, mr={0833268}}
\end{barticle}
%
\endbibitem

\bibitem{Besketal0101}
%
\begin{barticle}[mr]
\bauthor{\bsnm{Beskos},~\bfnm{Alexandros}\binits{A.}},
\bauthor{\bsnm{Roberts},~\bfnm{Gareth}\binits{G.}} \AND
\bauthor{\bsnm{Stuart},~\bfnm{Andrew}\binits{A.}}
(\byear{2009}).
\btitle{Optimal scalings for local {M}etropolis--{H}astings chains on nonproduct
targets in high dimensions}.
\bjournal{Ann. Appl. Probab.}
\bvolume{19}
\bpages{863--898}.
\bid{doi={10.1214/08-AAP563}, issn={1050-5164}, mr={2537193}}
\end{barticle}
%
\endbibitem

\bibitem{Besketal08}
%
\begin{barticle}[mr]
\bauthor{\bsnm{Beskos},~\bfnm{Alexandros}\binits{A.}},
\bauthor{\bsnm{Roberts},~\bfnm{Gareth}\binits{G.}},
\bauthor{\bsnm{Stuart},~\bfnm{Andrew}\binits{A.}} \AND
\bauthor{\bsnm{Voss},~\bfnm{Jochen}\binits{J.}}
(\byear{2008}).
\btitle{M{CMC} methods for diffusion bridges}.
\bjournal{Stoch. Dyn.}
\bvolume{8}
\bpages{319--350}.
\bid{doi={10.1142/S0219493708002378}, issn={0219-4937}, mr={2444507}}
\end{barticle}
%
\endbibitem

\bibitem{BeskStua07}
%
\begin{bincollection}[auto:STB|2011-03-03|12:04:44]
\bauthor{\bsnm{Beskos},~\bfnm{A.}\binits{A.}} \AND
\bauthor{\bsnm{Stuart},~\bfnm{A.~M.}\binits{A.~M.}}
(\byear{2008}).
\btitle{MCMC methods for sampling function space}.
In \bbooktitle{ICIAM Invited Lecture 2007}
(\beditor{\bfnm{R.}\binits{R.}~\bsnm{Jeltsch}} \AND
\beditor{\bfnm{G.}\binits{G.}~\bsnm{Wanner}}, eds.).
\bpublisher{European Mathematical Society}, \baddress{Z\"{u}rich}.
\end{bincollection}
%
\endbibitem

\bibitem{BouVan09}
%
\begin{barticle}[mr]
\bauthor{\bsnm{Bou-Rabee},~\bfnm{Nawaf}\binits{N.}} \AND
\bauthor{\bsnm{Vanden-Eijnden},~\bfnm{Eric}\binits{E.}}
(\byear{2010}).
\btitle{Pathwise accuracy and ergodicity of Metropolized integrators for
{SDE}s}.
\bjournal{Comm. Pure Appl. Math.}
\bvolume{63}
\bpages{655--696}.
\bid{issn={0010-3640}, mr={2583309}}
\bptnote{check year}%
\end{barticle}
%
\endbibitem

\bibitem{Breyetal04}
%
\begin{barticle}[mr]
\bauthor{\bsnm{Breyer},~\bfnm{Laird~Arnault}\binits{L.~A.}},
\bauthor{\bsnm{Piccioni},~\bfnm{Mauro}\binits{M.}} \AND
\bauthor{\bsnm{Scarlatti},~\bfnm{Sergio}\binits{S.}}
(\byear{2004}).
\btitle{Optimal scaling of {M}A{L}A for nonlinear regression}.
\bjournal{Ann. Appl. Probab.}
\bvolume{14}
\bpages{1479--1505}.
\bid{doi={10.1214/105051604000000369}, issn={1050-5164}, mr={2071431}}
\end{barticle}
%
\endbibitem

\bibitem{BreyRobe00}
%
\begin{barticle}[mr]
\bauthor{\bsnm{Breyer},~\bfnm{L.~A.}\binits{L.~A.}} \AND
\bauthor{\bsnm{Roberts},~\bfnm{G.~O.}\binits{G.~O.}}
(\byear{2000}).
\btitle{From {M}etropolis to diffusions: {G}ibbs states and optimal scaling}.
\bjournal{Stochastic Process. Appl.}
\bvolume{90}
\bpages{181--206}.
\bid{doi={10.1016/S0304-4149(00)00041-7}, issn={0304-4149}, mr={1794535}}
\end{barticle}
%
\endbibitem

\bibitem{Chat07}
%
\begin{bmisc}[auto:STB|2011-03-03|12:04:44]
\bauthor{\bsnm{Chatterjee},~\bfnm{S.}\binits{S.}}
(\byear{2007}).
\bhowpublished{Stein's method. Lecture notes.
Available at \texttt{
\href{http://www.stat.berkeley.edu/\textasciitilde sourav/stat206Afall07.html}{http://}
\href{http://www.stat.berkeley.edu/\textasciitilde sourav/stat206Afall07.html}{www.stat.berkeley.edu/\textasciitilde sourav/stat206Afall07.html}}.}
\end{bmisc}
%
\endbibitem

\bibitem{ChenWhit98}
%
\begin{barticle}[mr]
\bauthor{\bsnm{Chen},~\bfnm{Xiaohong}\binits{X.}} \AND
\bauthor{\bsnm{White},~\bfnm{Halbert}\binits{H.}}
(\byear{1998}).
\btitle{Central limit and functional central limit theorems for
{H}ilbert-valued dependent heterogeneous arrays with applications}.
\bjournal{Econometric Theory}
\bvolume{14}
\bpages{260--284}.
\bid{doi={10.1017/S0266466698142056}, issn={0266-4666}, mr={1629340}}
\end{barticle}
%
\endbibitem

\bibitem{cottetal09}
%
\begin{barticle}[auto:STB|2011-03-03|12:04:44]
\bauthor{\bsnm{Cotter},~\bfnm{S.~L.}\binits{S.~L.}},
\bauthor{\bsnm{Dashti},~\bfnm{M.}\binits{M.}} \AND
\bauthor{\bsnm{Stuart},~\bfnm{A.~M.}\binits{A.~M.}}
(\byear{2010}).
\btitle{Approximation of Bayesian inverse problems}.
\bjournal{SIAM Journal of Numerical Analysis}
\bvolume{48}
\bpages{322--345}.
\end{barticle}
%
\endbibitem

\bibitem{DaprZaby92}
%
\begin{bbook}[mr]
\bauthor{\bsnm{Da~Prato},~\bfnm{Giuseppe}\binits{G.}} \AND
\bauthor{\bsnm{Zabczyk},~\bfnm{Jerzy}\binits{J.}}
(\byear{1992}).
\btitle{Stochastic Equations in Infinite Dimensions}.
\bseries{Encyclopedia of Mathematics and Its Applications}
\bvolume{44}.
\bpublisher{Cambridge Univ. Press}, \baddress{Cambridge}.
\bid{doi={10.1017/CBO9780511666223}, mr={1207136}}
\end{bbook}
%
\endbibitem

\bibitem{EthiKurt86}
%
\begin{bbook}[mr]
\bauthor{\bsnm{Ethier},~\bfnm{Stewart~N.}\binits{S.~N.}} \AND
\bauthor{\bsnm{Kurtz},~\bfnm{Thomas~G.}\binits{T.~G.}}
(\byear{1986}).
\btitle{Markov Processes: Characterization and Convergence}.
\bpublisher{Wiley}, \baddress{New York}.
\bid{doi={10.1002/9780470316658}, mr={0838085}}
\end{bbook}
%
\endbibitem

\bibitem{HairStuaVoss07}
%
\begin{barticle}[mr]
\bauthor{\bsnm{Hairer},~\bfnm{M.}\binits{M.}},
\bauthor{\bsnm{Stuart},~\bfnm{A.~M.}\binits{A.~M.}} \AND
\bauthor{\bsnm{Voss},~\bfnm{J.}\binits{J.}}
(\byear{2007}).
\btitle{Analysis of {SPDE}s arising in path sampling. {II}. {T}he nonlinear
case}.
\bjournal{Ann. Appl. Probab.}
\bvolume{17}
\bpages{1657--1706}.
\bid{doi={10.1214/07-AAP441}, issn={1050-5164}, mr={2358638}}
\end{barticle}
%
\endbibitem

\bibitem{HairStuaVoss10}
%
\begin{bmisc}[auto:STB|2011-03-03|12:04:44]
\bauthor{\bsnm{Hairer},~\bfnm{M.}\binits{M.}},
\bauthor{\bsnm{Stuart},~\bfnm{A.~M.}\binits{A.~M.}} \AND
\bauthor{\bsnm{Voss},~\bfnm{J.}\binits{J.}}
(\byear{2011}).
\bhowpublished{Signal processing problems on function space: Bayesian
formulation, stochastic PDEs and effective MCMC methods.
In \textit{The Oxford Handbook of Nonlinear Filtering}
(D. Crisan and B. Rozovsky, eds.).
Oxford Univ. Press, Oxford.}
\end{bmisc}
%
\endbibitem

\bibitem{Hairetal05}
%
\begin{barticle}[mr]
\bauthor{\bsnm{Hairer},~\bfnm{M.}\binits{M.}},
\bauthor{\bsnm{Stuart},~\bfnm{A.~M.}\binits{A.~M.}},
\bauthor{\bsnm{Voss},~\bfnm{J.}\binits{J.}} \AND
\bauthor{\bsnm{Wiberg},~\bfnm{P.}\binits{P.}}
(\byear{2005}).
\btitle{Analysis of {SPDE}s arising in path sampling. {I}. {T}he {G}aussian
case}.
\bjournal{Commun. Math. Sci.}
\bvolume{3}
\bpages{587--603}.
\bid{issn={1539-6746}, mr={2188686}}
\end{barticle}
%
\endbibitem

\bibitem{Hast70}
%
\begin{barticle}[auto:STB|2011-03-03|12:04:44]
\bauthor{\bsnm{Hastings},~\bfnm{W.~K.}\binits{W.~K.}}
(\byear{1970}).
\btitle{Monte Carlo sampling methods using Markov chains and their
applications}.
\bjournal{Biometrika}
\bvolume{57}
\bpages{97--109}.
\end{barticle}
%
\endbibitem

\bibitem{Liu08}
%
\begin{bbook}[mr]
\bauthor{\bsnm{Liu},~\bfnm{Jun~S.}\binits{J.~S.}}
(\byear{2008}).
\btitle{Monte {C}arlo Strategies in Scientific Computing}.
\bpublisher{Springer}, \baddress{New York}.
\bid{mr={2401592}}
\end{bbook}
%
\endbibitem

\bibitem{MaRock92}
%
\begin{bbook}[mr]
\bauthor{\bsnm{Ma},~\bfnm{Zhi~Ming}\binits{Z.~M.}} \AND
\bauthor{\bsnm{R{\"o}ckner},~\bfnm{Michael}\binits{M.}}
(\byear{1992}).
\btitle{Introduction to the Theory of (nonsymmetric) {D}irichlet Forms}.
\bpublisher{Springer}, \baddress{Berlin}.
\bid{mr={1214375}}
\end{bbook}
%
\endbibitem

\bibitem{Metretal53}
%
\begin{barticle}[auto:STB|2011-03-03|12:04:44]
\bauthor{\bsnm{Metropolis},~\bfnm{N.}\binits{N.}},
\bauthor{\bsnm{Rosenbluth},~\bfnm{A.~W.}\binits{A.~W.}},
\bauthor{\bsnm{Teller},~\bfnm{M.~N.}\binits{M.~N.}} \AND
\bauthor{\bsnm{Teller},~\bfnm{E.}\binits{E.}}
(\byear{1953}).
\btitle{Equations of state calculations by fast computing machines}.
\bjournal{J. Chem. Phys.}
\bvolume{21}
\bpages{1087--1092}.
\end{barticle}
%
\endbibitem

\bibitem{CaseRobe04}
%
\begin{bbook}[mr]
\bauthor{\bsnm{Robert},~\bfnm{Christian~P.}\binits{C.~P.}} \AND
\bauthor{\bsnm{Casella},~\bfnm{George}\binits{G.}}
(\byear{2004}).
\btitle{Monte {C}arlo Statistical Methods},
\bedition{2nd} ed.
\bpublisher{Springer}, \baddress{New York}.
\bid{mr={2080278}}
\end{bbook}
%
\endbibitem

\bibitem{Robeetal97}
%
\begin{barticle}[mr]
\bauthor{\bsnm{Roberts},~\bfnm{G.~O.}\binits{G.~O.}},
\bauthor{\bsnm{Gelman},~\bfnm{A.}\binits{A.}} \AND
\bauthor{\bsnm{Gilks},~\bfnm{W.~R.}\binits{W.~R.}}
(\byear{1997}).
\btitle{Weak convergence and optimal scaling of random walk {M}etropolis
algorithms}.
\bjournal{Ann. Appl. Probab.}
\bvolume{7}
\bpages{110--120}.
\bid{doi={10.1214/aoap/1034625254}, issn={1050-5164}, mr={1428751}}
\end{barticle}
%
\endbibitem

\bibitem{RobeRose98}
%
\begin{barticle}[mr]
\bauthor{\bsnm{Roberts},~\bfnm{Gareth~O.}\binits{G.~O.}} \AND
\bauthor{\bsnm{Rosenthal},~\bfnm{Jeffrey~S.}\binits{J.~S.}}
(\byear{1998}).
\btitle{Optimal scaling of discrete approximations to {L}angevin diffusions}.
\bjournal{J. R. Stat. Soc. Ser. B Stat. Methodol.}
\bvolume{60}
\bpages{255--268}.
\bid{doi={10.1111/1467-9868.00123}, issn={1369-7412}, mr={1625691}}
\end{barticle}
%
\endbibitem

\bibitem{RobeRose01}
%
\begin{barticle}[mr]
\bauthor{\bsnm{Roberts},~\bfnm{Gareth~O.}\binits{G.~O.}} \AND
\bauthor{\bsnm{Rosenthal},~\bfnm{Jeffrey~S.}\binits{J.~S.}}
(\byear{2001}).
\btitle{Optimal scaling for various {M}etropolis--{H}astings algorithms}.
\bjournal{Statist. Sci.}
\bvolume{16}
\bpages{351--367}.
\bid{doi={10.1214/ss/1015346320}, issn={0883-4237}, mr={1888450}}
\end{barticle}
%
\endbibitem

\bibitem{Stro93}
%
\begin{bbook}[mr]
\bauthor{\bsnm{Stroock},~\bfnm{Daniel~W.}\binits{D.~W.}}
(\byear{1993}).
\btitle{Probability Theory, an Analytic View}.
\bpublisher{Cambridge Univ. Press}, \baddress{Cambridge}.
\bid{mr={1267569}}
\end{bbook}
%
\endbibitem

\bibitem{Stua10}
%
\begin{barticle}[mr]
\bauthor{\bsnm{Stuart},~\bfnm{A.~M.}\binits{A.~M.}}
(\byear{2010}).
\btitle{Inverse problems: A {B}ayesian perspective}.
\bjournal{Acta Numer.}
\bvolume{19}
\bpages{451--559}.
\bid{doi={10.1017/S0962492910000061}, issn={0962-4929}, mr={2652785}}
\end{barticle}
%
\endbibitem

\end{thebibliography}
\end{document}